\documentclass[12pt,reqno]{amsart}
\usepackage{latexsym}
\usepackage{amssymb}
\usepackage{epsfig}
\usepackage{multirow}
\usepackage{caption}
\usepackage{subfigure}
\usepackage{graphicx}
\usepackage{epstopdf}
\usepackage{graphicx}
\usepackage{subfigure}
\usepackage{color}
\allowdisplaybreaks
\numberwithin{equation}{section}

\usepackage{geometry}
\usepackage{algorithm}
\usepackage{algorithmic}

\newtheorem{example}{\textbf{Example}}[section]

\newtheorem{thm}{Theorem}[section]

\newtheorem{lem}[thm]{Lemma}

\newtheorem{ass}{Assumption}[section]
\newtheorem{definition}{Definition}[section]

\theoremstyle{definition}

\theoremstyle{remark}
\newtheorem{rem}{Remark}[section]

\def\bq{\begin{equation}}
\def\eq{\end{equation}}
\def\br{\begin{eqnarray}}
\def\er{\end{eqnarray}}
\def\brr{\bq\begin{array}{rlll}}
\def\err{\end{array}\eq}

\makeatletter
\newcommand{\definetitlefootnote}[1]{%
  \newcommand\addtitlefootnote{%
    \makebox[0pt][l]{$^{*}$}%
    \footnote{\protect\@titlefootnotetext}
  }%
  \newcommand\@titlefootnotetext{\spaceskip=\z@skip $^{*}$#1}%
}
\makeatother

\begin{document}
\title[AFEM for Cahn--Hilliard equation]
{Recovery type a posteriori error estimation of an adaptive finite element method for Cahn--Hilliard equation\addtitlefootnote}
\author[Y. Chen,\quad Y. Huang,\quad N. Yi\quad and P. Yin]
{Yaoyao Chen$^{\dag}$, Yunqing Huang$^\ddag$, Nianyu Yi$^{\S}$ and Peimeng Yin$^{\P}$}
\thanks{$^\dag$ School of Mathematics and Statistics, Anhui Normal University, Wuhu 241000, P.R.China,
Email: cyy1012xtu@126.com (Y. Chen).}
\thanks{$^\ddag$ Key Laboratory of Intelligent Computing \& Information Processing of Ministry of Education,
School of Mathematics and Computational Science, Xiangtan University, Xiangtan 411105, Hunan, P.R.China,
Email: huangyq@xtu.edu.cn (Y. Huang).}
\thanks{$^\S$ Hunan Key Laboratory for Computation and Simulation in Science and Engineering,
School of Mathematics and Computational Science, Xiangtan University, Xiangtan 411105, Hunan, P.R.China,
Email: yinianyu@xtu.edu.cn (N. Yi).}
\thanks{$^\P$ Multiscale Methods and Dynamics Group, Computer Science and Mathematics Division, Oak Ridge National Laboratory,  Oak Ridge, TN 37831, USA,
Email: yinp@ornl.gov (P. Yin).}


\definetitlefootnote{Notice:  This manuscript has been authored in part by UT-Battelle, LLC, under contract DE-AC05-00OR22725 with the US Department of Energy (DOE). The US government retains and the publisher, by accepting the article for publication, acknowledges that the US government retains a nonexclusive, paid-up, irrevocable, worldwide license to publish or reproduce the published form of this manuscript, or allow others to do so, for US government purposes. DOE will provide public access to these results of federally sponsored research in accordance with the DOE Public Access Plan (http://energy.gov/downloads/doe-public-access-plan).}

\begin{abstract}
In this paper, we derive a novel recovery type a posteriori error estimation of the Crank-Nicolson finite element method for the Cahn--Hilliard equation.
To achieve this, we employ both the elliptic reconstruction technique and a time reconstruction technique based on three time-level approximations, resulting in an optimal a posteriori error estimator.
We propose a time-space adaptive algorithm that utilizes the derived a posteriori error estimator as error indicators. Numerical experiments are presented to validate the theoretical findings, including comparing with an adaptive finite element method based on a residual type a posteriori error estimator.
\end{abstract}

\date{\today}
\maketitle

\hskip\parindent
\section{Introduction} \label{secInt}
\setcounter{equation}{0}


In this paper, we are interested in 
an adaptive finite element method for the Cahn--Hilliard equation

\bq\label{1e1}
\left\{\begin{aligned}
u_{t}+\mathcal{A}\left(\varepsilon\mathcal{A} u+\frac{1}{\varepsilon}f(u)\right)&=0,\qquad  \text{in}\ \Omega\times(0,T],\\
\partial_{\mathbf{n}}u\mid_{\partial\Omega}&=0,  \qquad \text{on}\ \partial\Omega\times[0,T],\\
\partial_{\mathbf{n}}\left(\varepsilon\mathcal{A} u+\frac{1}{\varepsilon}f(u)\right)\mid_{\partial\Omega}&=0, \qquad \text{on}\ \partial\Omega\times[0,T],\\
u(x,0)&=u_{0},\quad \ \ \text{in}\ \Omega\times\{t=0\},
\end{aligned}
\right.
\eq
where $\Omega\subset R^{d}(d=2,3)$ is a bounded domain with Lipschitz boundary $\partial\Omega$, $\mathbf{n}$ is the unit outward normal to the boundary $\partial\Omega$, the operator $\mathcal{A}:=-\Delta$, and the interface width $\varepsilon>0$ is a small parameter compared with the characteristic length of the laboratory scale. The nonlinear function $f(u)=F^{'}(u)=u^{3}-u$ with $F(u)=\frac{1}{4}(u^{2}-1)^{2}$, which is a double well potential and drives the solution to two pure states $u=\pm 1$.  


The Cahn--Hilliard equation, which was introduced by Cahn and Hilliard in the late 1950s to describe the process of phase separation \cite{CH1958}, has become a fundamental model in engineering and materials science. It also plays an increasingly important role in many other fields \cite{BF1993,EF1987}. The Cahn--Hilliard equation can be expressed as the $H^{-1}$-gradient flow, given by $u_{t}=\delta_{u} E(u)$, where $\delta_{u} E(u)$ is the variational derivative of the total free energy functional
\[E(u)=\int_{\Omega}\left(\frac{\varepsilon}{2}|\nabla u|^{2}+\frac{1}{\varepsilon}F(u)\right)dx.\]
It is well-known that the Cahn--Hilliard equation (\ref{1e1}), subject to the prescribed boundary conditions, satisfies an energy dissipative law given by
\[\frac{d}{dt}E(u(t))=-(u_{t},u_{t})\leq 0.\]


Efficient and easy-to-implement numerical methods for the Cahn--Hilliard equation face several challenges, including the presence of high order derivatives, a nonlinear reaction term $f(u)$, and the smallness of the parameter $\varepsilon$. To overcome these challenges, many spatial discretizations have been studied, including finite difference methods \cite{CFWW2019}, finite element methods \cite{DCW2016,FP2004, JLFL2020, WKG2006}, discontinuous Galerkin methods \cite{FK07, FLX16, LY2021}, and spectral methods \cite{DN1991}. Strategies to address the nonlinearity include convex-splitting methods \cite{GWWY2016}, stabilization methods \cite{SY2010}, invariant energy quantization (IEQ) approach \cite{YZWS2017, LY2021}, and scalar variable auxiliary (SAV) approach \cite{HAX2019,SXY2018,SXY2019}. Numerical approximations of the Cahn--Hilliard equation have been extensively investigated, but efficient and accurate methods are still an active area of research.


The smallness of the parameter $\varepsilon$ in gradient flow models, including the Cahn--Hilliard equation and the Allen--Cahn equation, results in the interface layer phenomenon. To accurately simulate macroscopic processes described by these equations, it is necessary to use adaptive techniques to adjust the spatial mesh size and time step size according to the interface width $\varepsilon$.  
In recent years, some works on a posteriori error estimators and adaptive methods have been proposed. Feng and Wu \cite{FW2008} developed residual-type a posteriori error estimates for conforming and mixed finite element approximations of the Cahn--Hilliard equation. 
A superconvergent cluster recovery (SCR)-based a posteriori error estimation and a time-space adaptive finite element algorithm was proposed in \cite{CHY2019n} for the Allen--Cahn equation. 
The SCR method produces a superconvergent recovered gradient, which leads to an asymptotically exact SCR-based error estimator.
The primary focus of \cite{CHY2019n} was to design an adaptive algorithm based on the SCR-based error estimator, while the time adaptation of the error indicator was simply constructed based on  approximations on two time levels.

In this paper, we present a novel SCR-based recovery type a posteriori error estimator for the Crank-Nicolson finite element method applied to the Cahn--Hilliard equation. 
The a posterior error estimator is derived using both the elliptic reconstruction technique and the time reconstruction technique. Therefore, the a posterior error estimator constructed is of greater precision and efficiency.
The elliptic reconstruction technique involves separating the error between the finite element approximation and the exact solution into two categories: elliptic type and parabolic type. The key idea is to leverage pre-existing elliptic a posteriori estimators for the elliptic type error, while controlling the parabolic type error using parabolic energy estimates.  
In \cite{CHY2019}, a time reconstruction technique using approximations on two time levels were introduced for the Allen--Cahn equation, which allowed for the construction of a first-order a posteriori error estimator for time discretization. In this work, we utilize the time reconstruction technique involving approximations on three time levels \cite{LPP2009}, leading to a second-order a posteriori error estimator for time discretization.
We employ the derived a posteriori error estimator as error indicators and propose an efficient time-space adaptive algorithm to solve the Cahn--Hilliard equation. Our numerical results show that the proposed recovery type a posteriori error estimator is more effective than a residual type error estimator and a space-only adaptive algorithm. Furthermore, our results demonstrate that the use of time step adaptation is essential in achieving accurate numerical solutions for the Cahn--Hilliard equation.

The paper is organized as follows: In Section \ref{secSch}, we introduce the Crank-Nicolson finite element method for discretizing the Cahn--Hilliard equation, followed by an introduction of the elliptic reconstruction for a nonlinear elliptic problem and its properties. In Section \ref{secEst}, we derive an optimal a posteriori error estimation for the Cahn--Hilliard equation based on the elliptic reconstruction and time reconstruction techniques. Based on the derived error estimator, we propose a time-space adaptive algorithm. In Section \ref{secNum}, we present several numerical examples to verify the accuracy and effectiveness of the proposed error indicators and the corresponding time-space adaptive algorithm. We present concluding remarks in Section \ref{secCon}. Finally, in Appendix \ref{aprofthm}, we provide the proof of Theorem \ref{thm3e2}.


\hskip\parindent
\section{The discrete scheme and elliptic reconstruction} \label{secSch}
\setcounter{equation}{0}
For a bounded domain $\Omega\subset{R}^{d}$, we adopt the standard notations for the Sobolev space $W^{m,p}(\Omega)$ equipped with
the norm $\|\cdot\|_{m,p,\Omega}$ and the semi-norm $|\cdot|_{m,p,\Omega}$. If $p=2$, we set
$W^{m,p}(\Omega)=H^m(\Omega)$, $\|\cdot\|_{m,p,\Omega}=\|\cdot\|_{m,\Omega}$ and $|\cdot|_{m,p,\Omega}=|\cdot|_{m,\Omega}$. Further, if $m=2$, we take $\|\cdot\|=\|\cdot\|_{0,\Omega}$.

By introducing the chemical potential
\bq\label{1e2}
w:=\varepsilon\mathcal{A} u+\frac{1}{\varepsilon}f(u),
\eq
we can get the equivalent form of \eqref{1e1},
\bq\label{1e3}
\left\{\begin{aligned}
u_{t}+\mathcal{A}w&=0,\qquad  \text{in}\ \Omega\times(0,T],\\
\partial_{\mathbf{n}}w\mid_{\partial\Omega}&=0, \qquad \text{on}\ \partial\Omega\times[0,T],\\
\varepsilon\mathcal{A} u+\frac{1}{\varepsilon}f(u)-w&=0,\qquad  \text{in}\ \Omega\times(0,T],\\
\partial_{\mathbf{n}}u\mid_{\partial\Omega}&=0,\qquad \text{on}\ \partial\Omega\times[0,T],\\
u(x,0)&=u_{0},\quad\ \ \text{in}\ \Omega\times\{t=0\}.
\end{aligned}
\right.
\eq

\subsection{The  Crank-Nicolson Finite Element Scheme}
For the homogeneous Neumann boundary conditions, the problem \eqref{1e3} is understood in the following weak form: find $(u,w)\in H^{1}(\Omega)\times H^{1}(\Omega)$ such that
\bq\label{1e4}
\left\{\begin{aligned}
(u_{t},v)+(\nabla w,\nabla v)&=0,\qquad  \forall v\in H^{1}(\Omega),\\
\varepsilon(\nabla u,\nabla \varphi)+\frac{1}{\varepsilon}\big(f(u),\varphi\big)-(w,\varphi)&=0,\qquad \forall \varphi\in H^{1}(\Omega),\\
u(\cdot,0)&=u_{0}.
\end{aligned}
\right.
\eq
Let $\mathcal{T}_{h}$ be a shape regular triangulation of $\Omega$, and $V_{h}$ be the corresponding finite element space, which is defined as
\[
V_{h}:=\left\{ v\in H^{1}(\Omega),v|_{K}\in P_{1}(K),\forall\ K\in\mathcal{T}_{h}\right\},
\]
where $P_{1}(K)$ denotes the set of linear polynomials defined in $K$. 
The semi-discrete finite element scheme of \eqref{1e3} reads: find $(u_{h},w_{h})\in V_{h}\times V_{h}$ such that
\bq\label{1e5}
\left\{\begin{aligned}
(u_{h,t},v_{h})+(\nabla w_{h},\nabla v_{h})&=0,\qquad  \forall v_{h}\in V_h,\\
\varepsilon(\nabla u_{h},\nabla \varphi_{h})+\frac{1}{\varepsilon}\big(f(u_{h}),\varphi_{h}\big)-(w_{h},\varphi_{h})&=0,\qquad \forall \varphi_{h}\in V_h,\\
 (u_{h}(x,0)-u_{0},\phi_h)&=0,\qquad \forall \phi_{h}\in V_h.
\end{aligned}
\right.
\eq
Generally, we rewrite the scheme \eqref{1e5} in its pointwise form
\bq\label{1e6}
\left\{\begin{aligned}
u_{h,t}+Aw_{h}&= 0,\\
\varepsilon Au_{h}+\frac{1}{\varepsilon}P f(u_{h})-w_{h}&=0,\\
 u_{h}(x,0)=u_{h}^{0}:&=Pu_{0},
\end{aligned}
\right.
\eq
where the finite-dimensional space operator $A:V_{h}\rightarrow V_{h}$ is the discrete Laplacian defined, through the Riesz representation in $V_{h}$, by
\[\langle Av,\Phi\rangle=a(v,\Phi),\qquad \forall\Phi\in V_{h},\]
and $P:L^{2}(\Omega)\rightarrow V_{h}$ is the $L^{2}(\Omega)$-projection operator such that, for each $v\in L^{2}(\Omega)$, we have
\[\langle Pv,\Phi\rangle=\langle v,\Phi\rangle,\qquad\forall \Phi\in V_{h}.\]

We divide the time interval $[0,T]$ into a partition of $N$ consecutive adjacent subintervals whose endpoints are denoted by $0=t_{0}<t_{1}<\cdots<t_{N}=T$, the $n$-th time interval $I_{n}:=[t_{n-1},t_{n}]$ and the corresponding time step is defined as $\tau_{n}:=t_{n}-t_{n-1}$.
The Crank-Nicolson finite element is to find a sequence of function $(u_{h}^{n},w_{h}^{n})\in V_{h}^{n}\times V_{h}^{n}$ such that, for each $n=1,2,\ldots,N$, 
\bq\label{1e7}
\left\{\begin{aligned}
\left(\frac{u_{h}^{n}-u_{h}^{n-1}}{\tau_{n}},v_{h}\right)+\frac{1}{2}\left(\nabla w_{h}^{n}+\nabla w_{h}^{n-1},\nabla v_{h}\right)&= 0,\qquad \forall v_{h} \in V_{h}^{n},\\
\frac{\varepsilon}{2}\left(\nabla u_{h}^{n}+\nabla u_{h}^{n-1},\nabla \varphi_{h}\right)+\frac{1}{\varepsilon}\left(\frac{f(u_{h}^{n})+f(u_{h}^{n-1})}{2},\varphi_{h}\right)\\-\frac{1}{2}\left(w_{h}^{n}+w_{h}^{n-1}, \varphi_{h}\right)&=0,\qquad \forall \varphi_{h} \in V_{h}^{n},\\
\qquad \qquad\qquad \qquad \qquad u_{h}(x,0)&=u_h^{0}.
\end{aligned}
\right.
\eq

Similarly to the semi-discrete scheme, the fully discrete scheme can be written in a pointwise form as follows
\bq\label{1e8}
\left\{\begin{aligned}
\frac{u_{h}^{n}-u_{h}^{n-1}}{\tau_{n}}+\frac{1}{2}\left(A^{n}w_{h}^{n}+A^{n-1}w_{h}^{n-1}\right)&= 0,\\
\frac{\varepsilon}{2}\left(A^{n}u_{h}^{n}+A^{n-1}u_{h}^{n-1}\right)+\frac{P^{n}f(u_{h}^{n})+P^{n-1}f(u_{h}^{n-1})}{2\varepsilon}
\\-\frac{1}{2}\left(w_{h}^{n}+w_{h}^{n-1}\right)&=0,\\
\qquad \qquad\qquad \qquad \qquad u_{h}(x,0)&=u_h^{0},
\end{aligned}
\right.
\eq
where $A^{n}: V_{h}^n\rightarrow V_{h}^n$ is defined as the discrete Laplacian and $P^{n}: L^{2}(\Omega)\rightarrow V_{h}^n$ represents the $L^{2}(\Omega)$-projection operator.

\subsection{Elliptic Reconstruction}
The nonlinear elliptic problem corresponding to a steady state of the nonlinear evolution equation \eqref{1e1} is taken as follows: given $g\in L^{2}(\Omega)$, $r\in L^{2}(\Omega)$, find $(\mu,\nu)\in H^{1}(\Omega)\times H^{1}(\Omega)$ such that
\bq\label{1e9}
\left\{\begin{aligned}
\mathcal{A}\nu+\nu&=g,\quad\ \ \, \text{in}\ \Omega,\\
\varepsilon\mathcal{A}\mu+\frac{1}{\varepsilon}h(\mu)-\nu&=r,\qquad \text{in}\ \Omega,\\
\nabla \mu\cdot \mathbf{n}=0,\,\nabla \nu\cdot \mathbf{n}&=0, \qquad \text{on}\ \partial\Omega,
\end{aligned}
\right.
\eq
with $h(\mu):=\mu^{3}$. The weak form of the elliptic problem \eqref{1e9} reads: find $(\mu,\nu)\in H^{1}(\Omega)\times H^{1}(\Omega)$ such that
\begin{flalign}
(\nabla \nu,\nabla v)+(\nu, v)&=\langle g,v\rangle,\qquad \forall v\in H^{1}(\Omega),\label{1e10a}\\
\varepsilon(\nabla \mu,\nabla \varphi)+\frac{1}{\varepsilon}\left(h(\mu),\varphi\right)-(\nu,\varphi)&=\langle r,\varphi\rangle,\qquad \forall \varphi\in H^{1}(\Omega).\label{1e10b}
\end{flalign}

\begin{rem}
The well-posedness of the variational problem \eqref{1e10a}-\eqref{1e10b} can be derived as follows. Owing to the variational problem \eqref{1e10a} is the Euler-Lagrange equation of the functional
\bq
J(\nu)=\frac{1}{2}\int_{\Omega}|\nabla \nu|^{2}+\frac{1}{2}\int_{\Omega}\nu^{2}-\int_{\Omega}g \nu,
\eq
taking the derivative of the functional $J(\nu)$, it holds that
\bq
\left(\frac{\delta J(\nu)}{\delta \nu},v\right)=\left(\nabla \nu,\nabla v\right)+\left(\nu,v\right)-\left(g,v\right)=0, \quad \forall v\in H^{1}(\Omega).
\eq
Notice that $J(\nu)$ is a convex functional, then the uniqueness of the solution for scheme \eqref{1e10a} is proved. As for the variational problem \eqref{1e10b}, it is the Euler-Lagrange equation of the functional
\bq
H(\mu)=\frac{\varepsilon}{2}\int_{\Omega}|\nabla \mu|^{2}+\frac{1}{4\varepsilon}\int_{\Omega}\mu^{4}-\int_{\Omega}s \mu,
\eq
where $\int_{\Omega}s \mu:=(\nu,\mu)+\langle r,\mu\rangle$.
Similarly, taking the derivative of the functional $H(\mu)$, it has
\bq
\left(\frac{\delta H(\mu)}{\delta \mu},\varphi\right)=\varepsilon\left(\nabla \mu,\nabla \varphi\right)+\frac{1}{\varepsilon}\left(h(\mu),\varphi\right)-\left(s,\varphi\right)=0,\quad \forall \varphi\in H^{1}(\Omega).
\eq
Due to $H(\mu)$ be a convex functional, then the uniqueness of the solution for scheme \eqref{1e10b} is proved.
\end{rem}
The finite element discretization of the elliptic problem \eqref{1e9} reads: find $(\mu_{h},\nu_{h})\in V_{h}\times V_{h}$ such that
\bq\label{1e11}
\left\{\begin{aligned}
(\nabla \nu_{h},\nabla v_{h})+(\nu_{h}, v_{h})&=\langle g_{h},v_{h}\rangle,\qquad \forall v_{h}\in V_{h},\\
\varepsilon(\nabla \mu_{h},\nabla \varphi_{h})+\frac{1}{\varepsilon}\left(h(\mu_{h}),\varphi_{h}\right)-(\nu_{h},\varphi_{h})&=\langle r_{h},\varphi_{h}\rangle,\qquad \forall \varphi_{h}\in V_{h}.
\end{aligned}
\right.
\eq
\begin{definition}\label{def2e1} (Gradient recovery a posteriori estimator function)
For the nonlinear elliptic problem \eqref{1e9}, we define the gradient recovery a posteriori estimator functional
\bq\label{1e12}
\mathcal{E}_{v}:=\mathcal{E}[v,H^{1}(\Omega),V_{h}]:=\|Gv-\nabla v\|,\qquad \forall v\in V_{h},
\eq
where $G$ is a gradient recovery operator.
\end{definition}

\begin{rem} As in \cite{LP2012}, we utilize $H^1(\Omega)$ to estimate the elliptic a posteriori estimation for the gradient recovery a posteriori estimator functional $\mathcal{E}_v$. However, it's worth noting that there are alternative methods to compute upper and lower bounds for the error in other functional spaces, such as $L^2(\Omega)$ and $L^{\infty}(\Omega)$.
\end{rem}

Gradient recovery is a post-processing technique that has gained widespread popularity in the engineering community for its robustness as an a posteriori error estimator, its superconvergence of the recovered derivatives, and its efficiency in implementation. It involves reconstructing gradient approximations from finite element solutions to obtain improved solutions. The practical use of the recovery technique is not only to enhance the quality of the approximation but also to construct a posteriori error estimators in adaptive computation.
The gradient of the finite element approximation for the Lagrange element provides a discontinuous approximation to the true gradient. Various techniques have been proposed to recover the gradient, including averaging \cite{BS1977,HJY2012}, local or global projections \cite{HTW2002,HLY2012}, postprocessing interpolation \cite{LY1996,Y2008}, the superconvergent patch recovery (SPR) \cite{ZZ1992}, the polynomial preserving recovery (PPR) \cite{ZN2005} and the superconvergent cluster recovery (SCR) \cite{HY2010}.

\begin{ass}\label{ass2e1} (Elliptic a posteriori error estimators)
Assume that $(\mu,\nu)$, $(\mu_{h},\nu_{h})$ are the exact solution and numerical solution of nonlinear elliptic problem \eqref{1e9}, respectively, $\mathcal{E}$ defined as Definition \ref{def2e1}, there exists constants $C_{0}$ and $C_{1}$, such that the following bounds hold
\bq\label{1e13}
\begin{aligned}
\|\nabla(\mu_{h}-\mu)\|&\leq C_{0}\mathcal{E}_{\mu},\\
\|\nabla(\nu_{h}-\nu)\|&\leq C_{1}\mathcal{E}_{\nu}.
\end{aligned}
\eq
\end{ass}
In \cite{HZ2016}, He and Zhou derived both a priori and a posteriori finite element error estimates for the following semilinear elliptic problems
\bq\label{1e13a}
\left\{\begin{aligned}
-\Delta u+b(x,u)&=0,\qquad  \text{in}\ \Omega\times(0,T],\\
u&=0,  \qquad \text{on}\ \partial\Omega\times[0,T],
\end{aligned}
\right.
\eq
and if the nonlinear term $b$ satisfies 
\bq\label{1e13ap}
\sup_{x\in \bar{\Omega}}\left|b(x,y)-b(x,y_{0})+\frac{\partial b}{\partial y}(x,y_{0})(y_{0}-y) \right|\lesssim (1+\max\{|y|^{s},|y_{0}|^{s}\})|y-y_{0}|^{q}, \ \forall y,\ y_{0}\in R 
\eq
with $q\in (1,2],\ s\in [0,5-q]$, then it has the following $L^{2}$ promotional property. 
 
\begin{lem}\label{ass2e10} \cite{HZ2016} ($L^{2}$ promotional property)
If $h_{0}\ll 1$, $h\in (0,h_{0}]$, then
\bq\label{2ea}\begin{aligned}
\|\mu_{h}-\mu\|&\leq C_{0}h \|\mu_{h}-\mu\|_{1,\Omega},\\
\|\nu_{h}-\nu\|&\leq C_{1}h \|\nu_{h}-\nu\|_{1,\Omega},
\end{aligned}\eq
here $C_{0}, C_{1}$ are constants and $h=\max\{h_{K},K\in\mathcal{T}_{h}\}$.
\end{lem}
\begin{rem}\label{rem1e1}
In this paper, for the nonlinear elliptic problem \eqref{1e9}, the nonlinear term $b(u):=h(u)=u^3$ satisfies
\[b(y)-b(y_0)-b'(y_0)(y-y_0)=(2y_0+y)(y-y_0)^2 ,\]
which is consistent with the condition \eqref{1e13ap} as $q=2$.
Thus, in view of Lemma \ref{ass2e10}, we have
\bq\label{2eb}
\begin{aligned}
\|\mu_{h}-\mu\|^{2}}_{1,\Omega&=\|\mu_{h}-\mu\|^{2}+\|\nabla(\mu_{h}-\mu)\|^{2}\\
&\leq C_{0}^{2}h^{2} \|\mu_{h}-\mu\|^{2}_{1,\Omega}+\|\nabla(\mu_{h}-\mu)\|^{2},\\
\|\nu_{h}-\nu\|^{2}}_{1,\Omega&=\|\nu_{h}-\nu\|^{2}+\|\nabla(\nu_{h}-\nu)\|^{2}\\
&\leq C_{1}^{2}h^{2} \|\nu_{h}-\nu\|^{2}_{1,\Omega}+\|\nabla(\nu_{h}-\nu)\|^{2},
\end{aligned}
\eq
thus, if $h$ is small enough, it holds that
\bq\label{2ed}\begin{aligned}
\|\mu_{h}-\mu\|_{1,\Omega}\leq C_{0}\|\nabla(\mu_{h}-\mu)\|\leq C_{0}\mathcal{E}_{\mu},\\
\|\nu_{h}-\nu\|_{1,\Omega}\leq C_{1}\|\nabla(\nu_{h}-\nu)\|\leq C_{1}\mathcal{E}_{\nu}.
\end{aligned}\eq
\end{rem}

To link the Cahn--Hilliard equation and the elliptic recovered gradient estimates, we utilize the elliptic reconstruction technique.

\begin{definition}\label{def2e2} (Elliptic reconstruction) For $1\leq n \leq N$ with the discrete elliptic operator $A^{n}$ defined as \eqref{1e8}, we define the corresponding elliptic reconstruction operator $R^{n}:V_{h}^{n}\rightarrow H^{1}(\Omega)$, for each $(\chi,\vartheta)\in V_{h}^{n}\times V_{h}^{n}$, by solving for the elliptic problem
\bq\label{1e14}
\left\{\begin{aligned}
\mathcal{A}R^{n}\vartheta+R^{n}\vartheta&=A^{n}\vartheta+\vartheta,\\
\varepsilon\mathcal{A}R^{n}\chi+\frac{1}{\varepsilon}h(R^{n}\chi)-R^{n}\vartheta&=\varepsilon A^{n}\chi+\frac{1}{\varepsilon}P^{n}h(\chi)-\vartheta,\\
\end{aligned}
\right.
\eq
which can be written in weak form as
\bq\label{1e15}
\left\{\begin{aligned}
a(R^{n}\vartheta,v)+(R^{n}\vartheta,v)=&\langle A^{n}\vartheta,v\rangle+(\vartheta,v),\ \forall v\in H^{1}(\Omega),\\
\varepsilon a(R^{n}\chi,\varphi)+\frac{1}{\varepsilon}\left(h(R^{n}\chi),\varphi\right)-(R^{n}\vartheta,\varphi)=&\varepsilon\langle A^{n}\chi,\varphi\rangle+\frac{1}{\varepsilon}\langle h(\chi),\varphi\rangle\\&-\langle\vartheta,\varphi\rangle,\ \forall\varphi\in H^{1}(\Omega).
\end{aligned}
\right.
\eq
\end{definition}


By the Definition \ref{def2e2}, it is obviously that $(R^{n}u_{h}^{n},R^{n}w_{h}^{n})$, $(u_{h}^{n},w_{h}^{n})$ are the exact solution and numerical solution of \eqref{1e14}, respectively. According to  Assumption \ref{ass2e1}, we have
\begin{flalign}
\|\nabla(u_{h}^{n}-R^{n}u_{h}^{n})\|&\leq C_{0}\mathcal{E}_{u}^{n},\label{1e16a}\\
\|\nabla(w_{h}^{n}-R^{n}w_{h}^{n})\|&\leq C_{1}\mathcal{E}_{w}^{n},\label{1e16b}
\end{flalign}
where $\mathcal{E}_{u}^{n},\ \mathcal{E}_{w}^{n}$ are defined  following Definition \ref{def2e1}, respectively, by
\begin{flalign}
\mathcal{E}_{u}^{n}&:=\|G^{n}u_{h}^{n}-\nabla u_{h}^{n}\|,\ \qquad \forall u_{h}^{n}\in V_{h}^{n},\label{1e16c}\\
\mathcal{E}_{w}^{n}&:=\|G^{n} w_{h}^{n}-\nabla w_{h}^{n}\|,\qquad \forall w_{h}^{n}\in V_{h}^{n},\label{1e16d}
\end{flalign}
with $G^{n}:=G^{V_{h}^{n}}$.

\hskip\parindent
\section{A posteriori error estimation and adaptive algorithm}\label{secEst}
\setcounter{equation}{0}
In this section, we derive a recovery type a posteriori error estimation for the Cahn--Hilliard equation based on the elliptic reconstruction and time reconstruction techniques, and a time-space adaptive algorithm is also developed based on the proposed a posteriori error estimation.
\subsection{A posteriori error estimation}
The discrete solution is sequence of finite element functions $u_h^{n} \in V_h^{n}$ defined at each discrete time $t_{n}, 1\leq n \leq N$.
Define the piecewise quadratic extension \cite{LPP2009}
\bq\label{2e1}
\begin{aligned}
u_{h}(t):&=\frac{t-t_{n-1}}{\tau_{n}}u_{h}^{n}+\frac{t_{n}-t}{\tau_{n}}u_{h}^{n-1}+\frac{1}{2}(t-t_{n-1})(t-t_{n})\partial_{n}^{2}u_{h},\ t\in I_{n},\ 1\leq n\leq N,\\
w_{h}(t):&=\frac{t-t_{n-1}}{\tau_{n}}w_{h}^{n}+\frac{t_{n}-t}{\tau_{n}}w_{h}^{n-1}+\frac{1}{2}(t-t_{n-1})(t-t_{n})\partial_{n}^{2}w_{h},\ t\in I_{n},\ 1\leq n\leq N,
\end{aligned}
\eq
where the term $\partial_{n}^{2}\nu_{h}$ is defined as
\bq\label{2e2}
\begin{aligned}
\partial_{n}^{2}\nu_{h}:=\frac{\frac{\nu_{h}^{n}-\nu_{h}^{n-1}}{\tau_{n}}-\frac{\nu_{h}^{n-1}-\nu_{h}^{n-2}}{\tau_{n-1}}}{\frac{\tau_{n}+\tau_{n-1}}{2}}
\end{aligned}
\eq
with
$\nu_{h}^{-1}=\nu_{h}^{0}$ as $n=1$.

Then we also define
\bq\label{2e3}
p^{n}:=R^{n}u_{h}^{n},\qquad q^{n}:=R^{n}w_{h}^{n},\qquad \ n=0,1,2,\ldots,N,\\
\eq
and denote this sequence's piecewise quadratic reconstruction in time by $p(t)$ and $q(t)$, that is,
\bq\label{2e4}
\begin{aligned}
p(t):&=\frac{t-t_{n-1}}{\tau_{n}}p^{n}+\frac{t_{n}-t}{\tau_{n}}p^{n-1}+\frac{1}{2}(t-t_{n-1})(t-t_{n})\partial_{n}^{2}p,\quad t\in I_{n},\ 1\leq n\leq N,\\
q(t):&=\frac{t-t_{n-1}}{\tau_{n}}q^{n}+\frac{t_{n}-t}{\tau_{n}}q^{n-1}+\frac{1}{2}(t-t_{n-1})(t-t_{n})\partial_{n}^{2}q,\quad  t\in I_{n},\ 1\leq n\leq N.
\end{aligned}
\eq
The corresponding fully discrete error is defined by
\bq\label{2e5}
\begin{aligned}
e_{u}:&=u_{h}(t)-u(t),\\
e_{w}:&=w_{h}(t)-w(t),
\end{aligned}
\eq
and can be split, using the elliptic reconstruction $p(t)$ and $q(t)$, as follows
\bq\label{2e6}
\begin{aligned}
e_{u}&=(p(t)-u(t))-(p(t)-u_{h}(t)):=\rho_{u}-\epsilon_{u},\\
e_{w}&=(q(t)-w(t))-(q(t)-w_{h}(t)):=\rho_{w}-\epsilon_{w}.
\end{aligned}
\eq
For terms in (\ref{2e6}), the following result holds.
\begin{thm}\label{thm3e1} (Parabolic error identity) For each $n=1,2,\ldots,N$ and each $t\in (t_{n-1},t_{n}]$, it holds that
\bq\label{2e9}
\begin{aligned}
\partial_{t}e_u+\mathcal{A}\rho_{w}
&=\frac{A^{n}w_{h}^{n}-A^{n-1} w_{h}^{n-1}}{2}+\mathcal{A}(q(t)-q^{n})+w_{h}^{n}-R^{n}w_{h}^{n}+(t-t_{n-\frac{1}{2}})\partial_{n}^{2}u_{h},\\
\varepsilon\mathcal{A}\rho_{u}-\rho_{w}
&=\varepsilon\frac{A^{n}u_{h}^{n}-A^{n-1} u_{h}^{n-1}}{2}+\frac{P^{n}f(u_{h}^{n})-P^{n-1}f(u_{h}^{n-1})}{2\varepsilon}-\frac{w_{h}^{n}-w_{h}^{n-1}}{2}\\
&\quad+\varepsilon\mathcal{A}(p(t)-p^{n})-(q(t)-q^{n})+\frac{f(u)-f(p^{n})}{\varepsilon}+\frac{1}{\varepsilon}\left(u_{h}^{n}-p^{n}\right),
\end{aligned}
\eq
where $A^n$ and $P^n$ are defined in (\ref{1e8}), respectively.
\end{thm}

\begin{proof}
For $n=1,2,\ldots,N$ and $t\ \in\ (t_{n-1},t_{n}]$, by the definition of $u_{h}^{n}$, we have
$$\partial_{t}u_{h}=\frac{u_{h}^{n}-u_{h}^{n-1}}{\tau_{n}}+(t-t_{n-\frac{1}{2}})\partial_{n}^{2}u_{h},$$
and using the fully discrete scheme \eqref{1e8}, we obtain
\[
\begin{aligned}
\partial_{t}u_{h}+\mathcal{A}q^{n}+q^{n}
=&\frac{u_{h}^{n}-u_{h}^{n-1}}{\tau_{n}}+A^{n}w_{h}^{n}+w_{h}^{n}+(t-t_{n-\frac{1}{2}})\partial_{n}^{2}u_{h}\\
=&\frac{u_{h}^{n}-u_{h}^{n-1}}{\tau_{n}}+\frac{A^{n}w_{h}^{n}+A^{n-1} w_{h}^{n-1}}{2}+\frac{A^{n}w_{h}^{n}-A^{n-1} w_{h}^{n-1}}{2}+w_{h}^{n}+(t-t_{n-\frac{1}{2}})\partial_{n}^{2}u_{h}\\
=&\frac{A^{n}w_{h}^{n}-A^{n-1} w_{h}^{n-1}}{2}+w_{h}^{n}+(t-t_{n-\frac{1}{2}})\partial_{n}^{2}u_{h},
\end{aligned}
\]
\[
\begin{aligned}
\varepsilon\mathcal{A}p^{n}&+\frac{1}{\varepsilon}h(p^{n})-q^{n}
=\varepsilon A^{n}u_{h}^{n}+\frac{1}{\varepsilon}P^{n}h(u_{h}^{n})-w_{h}^{n}\\
&=\varepsilon\frac{A^{n}u_{h}^{n}+A^{n-1} u_{h}^{n-1}}{2}+\varepsilon\frac{A^{n}u_{h}^{n}-A^{n-1} u_{h}^{n-1}}{2}+\frac{P^{n}h(u_{h}^{n})+P^{n-1}f(u_{h}^{n-1})}{2\varepsilon}\\
&\quad+\frac{P^{n}h(u_{h}^{n})-P^{n-1}f(u_{h}^{n-1})}{2\varepsilon}-\frac{w_{h}^{n}+w_{h}^{n-1}}{2}-\frac{w_{h}^{n}-w_{h}^{n-1}}{2}\\
&=\varepsilon\frac{A^{n}u_{h}^{n}-A^{n-1} u_{h}^{n-1}}{2}+\frac{P^{n}f(u_{h}^{n})
-P^{n-1}f(u_{h}^{n-1})}{2\varepsilon}-\frac{w_{h}^{n}-w_{h}^{n-1}}{2}+\frac{1}{\varepsilon}u_{h}^{n}.
\end{aligned}\]
Hence
\[
\begin{aligned}
\partial_{t}u_{h}+\mathcal{A}q(t)
&=\frac{A^{n}w_{h}^{n}-A^{n-1} w_{h}^{n-1}}{2}+\mathcal{A}(q(t)-q^{n})+w_{h}^{n}-R^{n}w_{h}^{n}+(t-t_{n-\frac{1}{2}})\partial_{n}^{2}u_{h},\\
\varepsilon\mathcal{A}p(t)+\frac{1}{\varepsilon}h(p^{n})-q(t)
&=\varepsilon\frac{A^{n}u_{h}^{n}-A^{n-1} u_{h}^{n-1}}{2}+\frac{P^{n}f(u_{h}^{n})-P^{n-1}f(u_{h}^{n-1})}{2\varepsilon}-\frac{w_{h}^{n}-w_{h}^{n-1}}{2}\\
&\quad+\varepsilon\mathcal{A}(p(t)-p^{n})-(q(t)-q^{n})+\frac{1}{\varepsilon}u_{h}^{n},
\end{aligned}
\]
and subtracting \eqref{1e3} from the above formula, we get
\[
\begin{aligned}
\partial_{t}e_u+\mathcal{A}\rho_{w}
&=\frac{A^{n}w_{h}^{n}-A^{n-1} w_{h}^{n-1}}{2}+\mathcal{A}(q(t)-q^{n})+w_{h}^{n}-R^{n}w_{h}^{n}+(t-t_{n-\frac{1}{2}})\partial_{n}^{2}u_{h},\\
\varepsilon\mathcal{A}\rho_{u}-\rho_{w}
&=\varepsilon\frac{A^{n}u_{h}^{n}-A^{n-1} u_{h}^{n-1}}{2}+\frac{P^{n}f(u_{h}^{n})-P^{n-1}f(u_{h}^{n-1})}{2\varepsilon}-\frac{w_{h}^{n}-w_{h}^{n-1}}{2}\\
&\quad+\varepsilon\mathcal{A}(p(t)-p^{n})-(q(t)-q^{n})+\frac{f(u)-f(p^{n})}{\varepsilon}+\frac{1}{\varepsilon}\left(u_{h}^{n}-p^{n}\right).
\end{aligned}
\]
\end{proof}

Then we have the following result.

\begin{thm}\label{thm3e2}
Let $(u_{h}^{n}, w_h^{n})_{n\in[0:N]}$ be the fully discrete solution, defined at each discrete time $t_{n}$, its piecewise linear extension $u_{h},\ w_{h}$ defined as \eqref{2e1}, and let $u,\ w$ be the exact solution of the model problem \eqref{1e3}. 
Assume that $\overline{\Lambda}_{CH}\in L^{1}(0,T)$ is a function such that for almost every $t\in (0,T)$, we have
\[-\overline{\Lambda}_{CH}(t)\leq-\Lambda_{CH}(t):=\inf_{v\in \dot{V}\backslash\{0\}}\frac{\varepsilon\|\nabla v\|^{2}+\varepsilon^{-1}(f'(u_{h})v,v)}{\|\nabla\Delta^{-1}v\|^{2}},\]
and set
\[
\begin{aligned}
a(t):=&\left(1+\frac{5}{2\varepsilon^{2}}+2\left(1-\varepsilon\right)\overline{\Lambda}_{CH}(t)\right),\\
\mu_{g}:=&\sup\limits_{t\in (0,T)}\|\tilde{f}(u_{h})\|_{L^{\infty}(\Omega)}.
\end{aligned}
\]
Define
\[\eta^{2}:=\|\nabla\Delta^{-1}e_{u}^{0}\|^{2}+\sum_{n=1}^{N}4\widetilde{\mathcal{E}_{u}^{n}}^{2}+\sum_{n=1}^{N}\left(\eta_{0}^{2}+\eta_{1}^{2}\right)\tau_{n},\]
and assume
\[
\eta^{2}\leq\frac{\varepsilon^{2/\sigma}}{(2\mu_{g}C_{S}(1+T))^{1/\sigma}}\left(8\exp\left(\int_{0}^{T}a(t)dt\right)\right)^{-1-\frac{1}{\sigma}},
\]
then
\bq\label{2e10}
\begin{aligned}
\sup_{t\in[0,T]}\|\nabla\Delta^{-1}e_{u}\|^{2}+\int_{0}^{T}\frac{\varepsilon^{2}}{2}\|\nabla e_{u}\|^{2}dt\leq8\eta^{2}\exp\left(\int_{0}^{T}a(t)dt\right),
\end{aligned}
\eq
where
\begin{align*}
\eta_{0}:=&\gamma_{w}^{n}+\delta_{w}^{n}+\eta_{w}^{n}+\beta_{u}^{n};\\
\eta_{1}:=&\gamma_{u}^{n}+\xi_{u}^{n}+\beta_{w}^{n}+\theta_{u}^{n}+\delta_{u}^{n}+\alpha_{u}^{n}+\zeta_{u}^{n};\\
\widetilde{\mathcal{E}_{u}^{n}}^{2}:=&\frac{C_{0}^{2}}{3}\tau_{n}\left((\mathcal{E}_{u}^{n})^{2}+(\mathcal{E}_{u}^{n-1})^{2}+\mathcal{E}_{u}^{n}\mathcal{E}_{u}^{n-1}\right)\\
&+C_{0}^{2}\frac{\tau_{n}^{2}\tau_{n-1}\Big(\mathcal{E}_{u}^{n}+\mathcal{E}_{u}^{n-1}\Big)+\tau_{n}^{3}\Big(\mathcal{E}_{u-1}^{n}+\mathcal{E}_{u}^{n-2}\Big)}{6\tau_{n-1}(\tau_{n}+\tau_{n-1})}\Big(\mathcal{E}_{u}^{n}+\mathcal{E}_{u}^{n-1}\Big)\nonumber\\
&+C_{0}^{2}\tau_{n}^{3}\frac{\Big(\tau_{n-1}\Big(\mathcal{E}_{u}^{n}+\mathcal{E}_{u}^{n-1}\Big)+\tau_{n}\Big(\mathcal{E}_{u-1}^{n}+\mathcal{E}_{u}^{n-2}\Big)\Big)^{2}}{30\tau_{n-1}^{2}(\tau_{n}+\tau_{n-1})^{2}};\\
\gamma_{w}^{n}:=&\left\|\frac{A^{n}w_{h}^{n}-A^{n-1} w_{h}^{n-1}}{2}\right\|_{-1};\\
\beta_{u}^{n}:=&\left\|\frac{\tau_{n}^{2}}{8}\cdot\partial_{n}^{2}u_{h}\right\|_{-1};\\
\eta_{w}^{n}:=&\left\|(A^{n-1}w_{h}^{n-1}+w_{h}^{n-1})-(A^{n}w_{h}^{n}+w_{h}^{n})\right\|_{-1}+\left\|\frac{\tau_{n}^{2}}{8}\partial_{n}^{2}(Aw_{h}+w_{h})\right\|_{-1};\\
\delta_{w}^{n}:=&\left\|w_{h}^{n}-w_{h}^{n-1}\right\|_{-1}+\left\|\frac{\tau_{n}^{2}}{8}\partial_{n}^{2}w_{h}\right\|_{-1};\\
\delta_{u}^{n}:=&\left\|\frac{u_{h}^{n}-u_{h}^{n-1}}{\varepsilon}\right\|;\\
\gamma_{u}^{n}:=&\varepsilon\left\|\frac{A^{n}u_{h}^{n}-A^{n-1} u_{h}^{n-1}}{2}\right\|;\\
\xi_{u}^{n}:=&\left\|\frac{P^{n}f(u_{h}^{n})-P^{n-1}f(u_{h}^{n-1})}{2\varepsilon}\right\|;\\
\beta_{w}^{n}:=&\left\|\frac{w_{h}^{n}-w_{h}^{n-1}}{2}\right\|;\\
\theta_{u}^{n}:=&2\left\|\left(\varepsilon A^{n-1}u_{h}^{n-1}+\frac{1}{\varepsilon}h(u_{h}^{n-1})-w_{h}^{n-1}\right)-\left(\varepsilon A^{n}u_{h}^{n}+\frac{1}{\varepsilon}h(u_{h}^{n})-w_{h}^{n}\right)\right\|\\&\quad+\Big\|\Big(\varepsilon A^{n-1}u_{h}^{n-1}+\frac{1}{\varepsilon}h(u_{h}^{n-1})-w_{h}^{n-1}\Big)-\Big(\varepsilon A^{n-2}u_{h}^{n-2}+\frac{1}{\varepsilon}h(u_{h}^{n-2})-w_{h}^{n-2}\Big) \Big\|;\\
\alpha_{u}^{n}:=&\frac{1}{\varepsilon}C\left(\mathcal{E}_{u}^{n}+\mathcal{E}_{u}^{n-1}+\mathcal{E}_{u}^{n-2}\right);\\
\zeta_{u}^{n}:=&\left\|\frac{3(u_{h}^{n})^{2}u_{h}^{n-1}-2(u_{h}^{n})^{3}-(u_{h}^{n-1})^{3}}{\varepsilon}\right\|+\left\|\frac{3u_{h}^{n}(u_{h}^{n-1})^{2}-2(u_{h}^{n-1})^{3}-(u_{h}^{n})^{3}}{\varepsilon} \right\|\\
&\quad+\left\|\frac{\big(3(u_{h}^{n})^{2}\partial_{n}^{2}u_{h}-3u_{h}^{n}u_{h}^{n-1}\partial_{n}^{2}u_{h}\big)\tau_{n}^{2}}{8\varepsilon} \right\|+\left\|\frac{\big(3(u_{h}^{n-1})^{2}\partial_{n}^{2}u_{h}-3u_{h}^{n}u_{h}^{n-1}\partial_{n}^{2}u_{h}\big)\tau_{n}^{2}}{8\varepsilon} \right\|\\
&\quad+\left\|\frac{\big(3u_{h}^{n}u_{h}^{n-1}\partial_{n}^{2}u_{h}-\frac{\frac{(u_{h}^{n})^{3}-(u_{h}^{n-1})^{3}}{\tau_{n}}-\frac{(u_{h}^{n-1})^{3}-(u_{h}^{n-2})^{3}}{\tau_{n-1}}}{\frac{\tau_{n}+\tau_{n-1}}{2}}\big)\tau_{n}^{2}}{8\varepsilon} \right\|\\
&\quad+\left\|\frac{\big(3u_{h}^{n}(\partial_{n}^{2}u_{h})^{2}-3u_{h}^{n-1}(\partial_{n}^{2}u_{h})^{2}\big)\tau_{n}^{4}}{64\varepsilon} \right\|+\left\|\frac{\big(3u_{h}^{n-1}(\partial_{n}^{2}u_{h})^{2}\big)\tau_{n}^{4}}{64\varepsilon}\right\|+\left\|\frac{(\partial_{n}^{2}u_{h})^{3}\tau_{n}^{6}}{512\varepsilon} \right\|,
\end{align*}
here $C_{S},C$ are constants, which are independent of mesh size, and $\mathcal{E}_{u}^{n}:=\mathcal{E}[u_{h}^{n}]$ is defined as \eqref{1e16c}.
\end{thm}
The proof of this theorem is provided in Appendix \ref{aprofthm}.

\begin{rem}\label{rem3e01}
The a posteriori error estimator in Theorem \ref{thm3e2} can be divided into two categories. The terms $\gamma_{w}^{n},\ \beta_{u}^{n},\ \eta_{w}^{n},\ \delta_{w}^{n},\ \delta_{u}^{n},\ \gamma_{u}^{n},\ \xi_{u}^{n},\ \beta_{w}^{n},\ \theta_{u}^{n},\ \zeta_{u}^{n}$ are viewed as the a posteriori error
indicators for time discretization, the terms $\widetilde{\mathcal{E}_{u}^{n}}$ and $\alpha_{u}^{n}$ are the spatial discretization error indicators.
\end{rem}


\subsection{Adaptive Algorithm}
In view of the a posteriori error estimator of Theorem \ref{thm3e2}, we design the algorithms for time-step size control and spatial adaptation in this part.

We adjust the time-step size in view of the error equidistribution strategy, which means that the time discretization error should be evenly distributed to each time interval $(t_{n-1},t_{n}],\ n=1,2,\ldots,N$.  Let $TOL_{time}$ be the tolerance allowed for the part of the a posteriori error estimator in \eqref{2e10} related to the time discretization, that is,
\bq\label{2e27b}
\sum_{n=1}^{N}\tau_{n}(\gamma_{w}^{n}+\beta_{u}^{n}+\eta_{w}^{n}+\delta_{w}^{n}+\delta_{u}^{n}+\gamma_{u}^{n}+\xi_{u}^{n}+\beta_{w}^{n}+\theta_{u}^{n}+\zeta_{u}^{n})^{2}\leq TOL_{time}.
\eq
Generally, we can achieve \eqref{2e27b} by adjusting the time-step size $\tau_{n}$ so as to have the following relations
\bq\label{2e28}
\eta^{n}_{time}:=\gamma_{w}^{n}+\beta_{u}^{n}+\eta_{w}^{n}+\delta_{w}^{n}+\delta_{u}^{n}+\gamma_{u}^{n}+\xi_{u}^{n}+\beta_{w}^{n}+\theta_{u}^{n}+\zeta_{u}^{n}\leq \sqrt{TOL_{time}/T}:=TOL_{t}.
\eq
We summarize the procedure of time-step size control in Algorithm \ref{algtimadp}.

\begin{algorithm}
\caption{Time-step size control}
\begin{algorithmic}[1]\label{algtimadp}
\STATE Given tolerances $TOL_{t}$, $TOL_{t,m}:=\sqrt{TOL_{time,min}/T}$, parameters $\delta_{1}\in (0,1)$, $\delta_{2}>1$;\;

\STATE Set $\tau_{n}:=\tau_{n-1}$, $t_{n}:=t_{n-1}+\tau_{n}$;\;

\STATE Solve the discrete problem and compute the time error estimator $\eta^{n}_{time}$;\;

\WHILE{$\eta^{n}_{time}>TOL_{t}$ or $\eta^{n}_{time}<TOL_{t,m}$}

\IF{$\eta^{n}_{time}>TOL_{t}$}
\STATE Set $\tau_{n}:=\delta_{1}\cdot\tau_{n}$ and $t_{n}:=t_{n-1}+\tau_{n}$;\;
\ELSE
\STATE Set $\tau_{n}:=\delta_{2}\cdot\tau_{n}$ and $t_{n}:=t_{n-1}+\tau_{n}$;\;
\ENDIF

\STATE Solve the discrete problem and compute the time error estimator $\eta^{n}_{time}$;\;

\ENDWHILE
  \end{algorithmic}
\end{algorithm}

Let $TOL_{space}$ be the tolerance allowed for the part of the a posteriori error estimator in \eqref{2e10} related to the spatial discretization. For the recovery type error estimator, we adopt the SCR-based error estimator. The SCR gradient recovery method was proposed by Huang and Yi in \cite{HY2010,Y2011}, it can produce a superconvergent recovered gradient, which in turn provides the SCR-based error estimator that is asymptotically exact. Similar to time discretization, we aim to achieve the following relation at each time step $n$,
\bq\label{2e29}
\eta^{n}_{space}:=\widetilde{\mathcal{E}_{u}^{n}}+\alpha_{u}^{n}\leq \sqrt{TOL_{space}/T}:=TOL_{s}.
\eq

Given the refinement and the coarsening parameters $TOL_{r}$, $TOL_{c}$, respectively, we adopt the following Maximum mark strategy to mark the elements for refinement or coarsening. Set
\bq\label{2e31}
\begin{aligned}
\eta^{n}_{K}:=\left\|G^{n} u_{h}^{n}-\nabla u_{h}^{n}\right\|_K, \
\eta_{ \max}^{n}:=\max\{\eta^{n}_{K},K\in\mathcal{T}^{n}_{h} \},
\end{aligned}
\eq
choose the elements $\{K: \eta^{n}_{K}> TOL_{r}\times\eta_{\max}^{n}\}$ for refinement, and 
choose the elements $\{K: \eta^{n}_{K}< TOL_{c}\times\eta_{\max}^{n}\}$ for coarsening.

In view of the error indicators above, we design the following time-space adaptive algorithm for Cahn--Hilliard equation, which is outlined in Algorithm \ref{algadp}.

  \begin{algorithm}
  \caption{Time-space adaptive algorithm for the Cahn--Hilliard equation}
\begin{algorithmic}[1]\label{algadp}
\STATE Given $TOL_{t}$, $TOL_{t,m}$, $TOL_{s}$, $TOL_{i}$, $\delta_{1}\in (0,1)$, $\delta_{2}>1$;\;

\STATE Given the initial time step $\tau_{0}$, initial mesh $\mathcal{T}^{0}_{h}$, and initial solution $u_{h}^{0}$;\;

\STATE Set $n=0$, $t_{0}=0$, $E(u_{h}^{-1})=0$;\;

\STATE Compute the initial error estimator $\eta^{0}_{initial}=\|u_{0}-u_{h}^{0}\|$;\;

\STATE Refine $\mathcal{T}^{0}_{h}$ to get a mesh such that $\eta^{0}_{initial}\leq TOL_{i}$;\;

\STATE Compute the energy $E(u_{h}^{0})$;\;

\WHILE{$E(u_{h}^{n})-E(u_{h}^{n-1})>TOL_e$}

\STATE Set $n:=n+1$, $\mathcal{T}^{n}_{h}:=\mathcal{T}^{n-1}_{h}$, $\tau_{n}:=\tau_{n-1}$, $t_{n}:=t_{n-1}+\tau_{n}$;\;

\STATE Solve the discrete problem and compute the time error estimator $\eta^{n}_{time}$;\;

\WHILE{$\eta^{n}_{time}>TOL_{t}$ or $\eta^{n}_{time}<TOL_{t,m}$}

\IF{$\eta^{n}_{time}>TOL_{t}$}

\STATE Set $\tau_{n}:=\delta_{1}\cdot\tau_{n}$ and $t_{n}:=t_{n-1}+\tau_{n}$;\;
\ELSE
\STATE Set $\tau_{n}:=\delta_{2}\cdot\tau_{n}$ and $t_{n}:=t_{n-1}+\tau_{n}$;\;
\ENDIF
\STATE Solve the discrete problem and compute the time error estimator $\eta^{n}_{time}$;\;
\ENDWHILE
\STATE Compute the space error estimator $\eta^{n}_{space}$, $\eta^{n}_{K}$ and $\eta_{\max}^{n}$;\;
\WHILE{$\eta^{n}_{space}> TOL_{s}$}
\STATE Mark elements for refinement;\; 

\STATE Refine mesh $\mathcal{T}^{n}_{h}$ to generate a new mesh $\mathcal{T}^{n}_{h}$;\;

\STATE Solve the discrete problem for $u_{h}^{n}$ on the new mesh $\mathcal{T}^{n}_{h}$ using data $u_{h}^{n-1}$;\;

\STATE Compute the time error estimator $\eta^{n}_{time}$;\;

\WHILE{$\eta^{n}_{time}>TOL_{t}$ or $\eta^{n}_{time}<TOL_{t,m}$}

\IF{$\eta^{n}_{time}>TOL_{t}$}

\STATE Set $\tau_{n}:=\delta_{1}\cdot\tau_{n}$ and $t_{n}:=t_{n-1}+\tau_{n}$;\;

\ELSE

\STATE Set $\tau_{n}:=\delta_{2}\cdot\tau_{n}$ and $t_{n}:=t_{n-1}+\tau_{n}$;\;

\ENDIF

\STATE Solve the discrete problem and compute the time error estimator $\eta^{n}_{time}$;\;

\ENDWHILE

\STATE Compute the space error estimator $\eta^{n}_{space}$, $\eta^{n}_{K}$ and $\eta_{\max}^{n}$;\;

\ENDWHILE

%
%
%
%

\STATE Compute the energy $E(u_{h}^{n})$ ;\;

\STATE Mark elements for coarsen and coarsen $\mathcal{T}^{n}_{h}$ producing a modified mesh $\mathcal{T}^{n}_{h}$;\;

\ENDWHILE

  \end{algorithmic}
\end{algorithm}

\hskip\parindent
\section{Numerical examples} \label{secNum}
\setcounter{equation}{0}

In this section, we present three examples to demonstrate the reliability and effectiveness of the proposed adaptive algorithm based on the a posteriori error estimator of Theorem \ref{thm3e2}. In Example\ \ref{exm1}, we investigate the main part of the space and time discretization error indicators numerically. In Example\ \ref{cexm2}, we focus on illustrating the efficiency of the a posteriori error estimator based on the recovery type and the necessity of time-space adaptation by comparing them with the residual type and space-only adaptation, respectively. We provide the corresponding numerical results, including the discrete energy history, the change in the number of nodes and time steps, the numerical solutions, adaptive meshes, and CPU time, to support our conclusions. For the last example, we apply the proposed time-space adaptive algorithm to the three-dimensional Cahn--Hilliard equation.

In all examples, we take the parameters
\[\delta_{1}=\frac{1}{2},\qquad \delta_{2}=2,\]
and the remaining parameters will be specified in each example.

\begin{example}\label{exm1} Consider the Cahn--Hilliard equation \eqref{1e3} with the initial condition
\[u_{0}(x,y)=\tanh\Big(\big((x-0.3)^{2}+y^{2}-0.25^{2}\big)/\varepsilon\Big)\tanh\Big(\big((x+0.3)^{2}+y^{2}-0.3^{2}\big)/\varepsilon\Big),\]
where $\Omega=[-1,1]^{2}$ and the parameters $\varepsilon=0.01$, $TOL_{t}=50$, $TOL_{t,m}=5$, $TOL_{s}=10$,  $TOL_{i}=0.002$.
\end{example}


\begin{figure}[!htbp]
$\begin{array}{cccc}
\includegraphics[width=3.8cm,height=3.0cm]{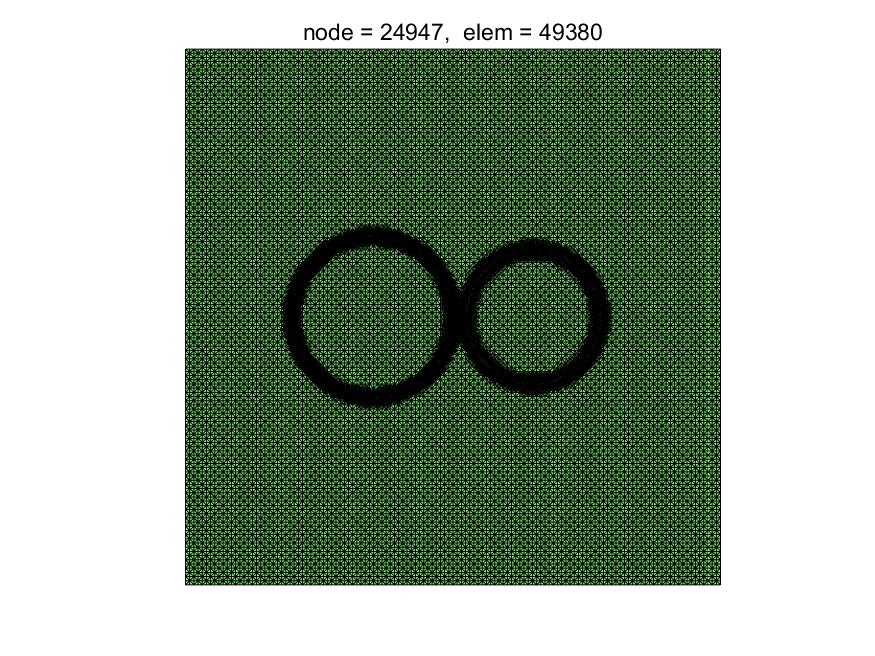}&
\includegraphics[width=3.8cm,height=3.0cm]{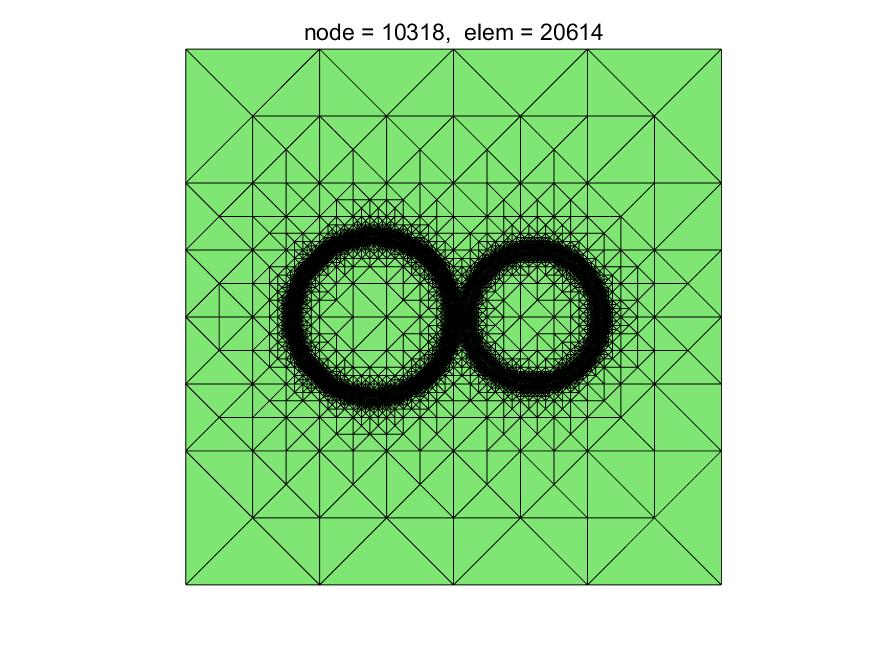}&
\includegraphics[width=3.8cm,height=3.0cm]{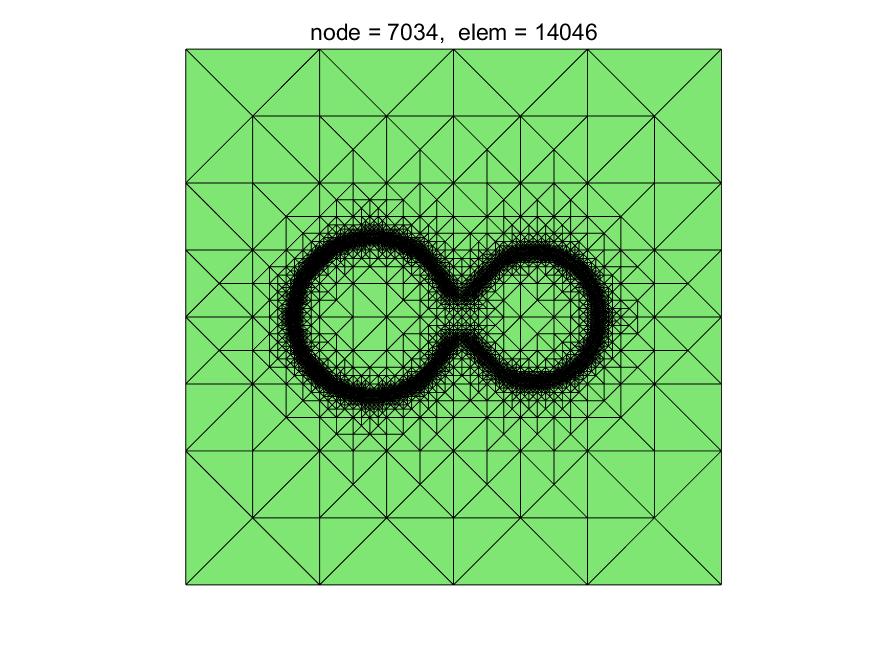}&
\includegraphics[width=3.8cm,height=3.0cm]{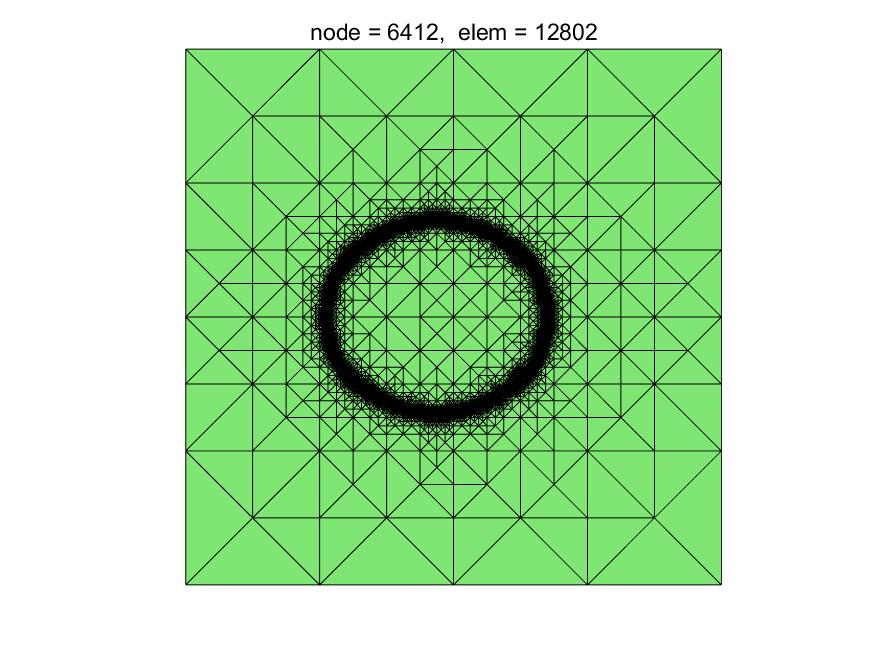}\\
\includegraphics[width=3.8cm,height=3.0cm]{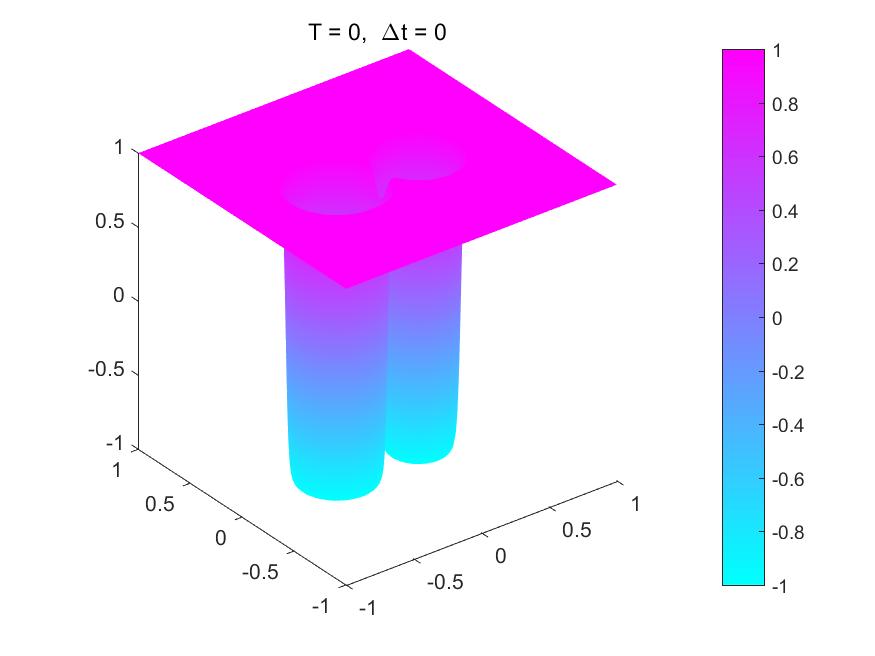}&
\includegraphics[width=3.8cm,height=3.0cm]{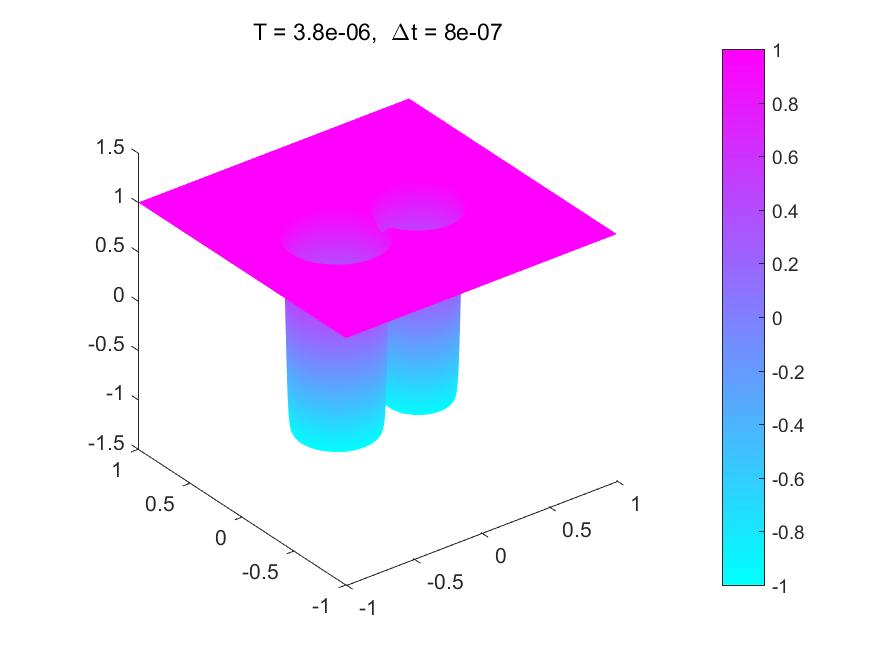}&
\includegraphics[width=3.8cm,height=3.0cm]{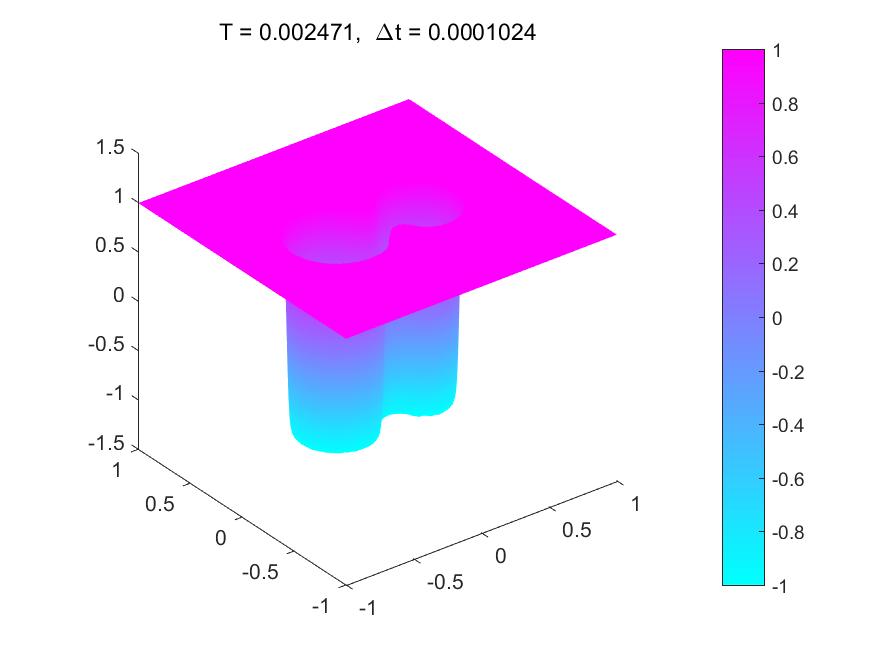}&
\includegraphics[width=3.8cm,height=3.0cm]{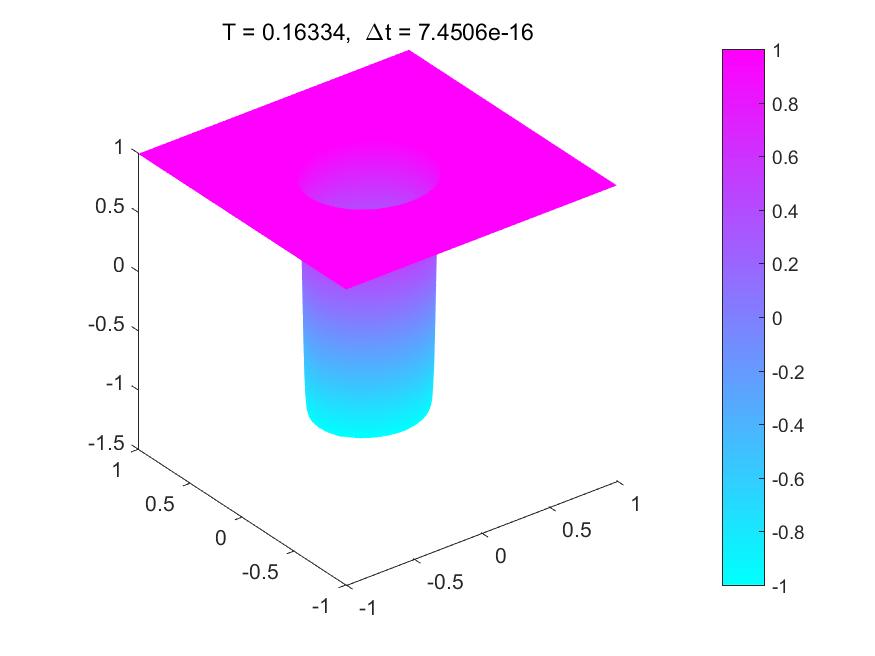}
\end{array}$
\caption{Example\ \ref{exm1}, First line: adaptive meshes; Second line: snapshots of numerical solutions for $u$.}\label{exp0u}
\end{figure}

We apply the proposed time-space adaptive algorithm to solve the Cahn--Hilliard equation, the numerical solutions of $u$ and the corresponding adaptive meshes are shown in Figure \ref{exp0u}, respectively. From the pictures, we can see that the meshes follow the zeros level set of $u$ as it moves.

\begin{figure}[!htbp]
$\begin{array}{ccc}
\includegraphics[width=5.5cm,height=4.5cm]{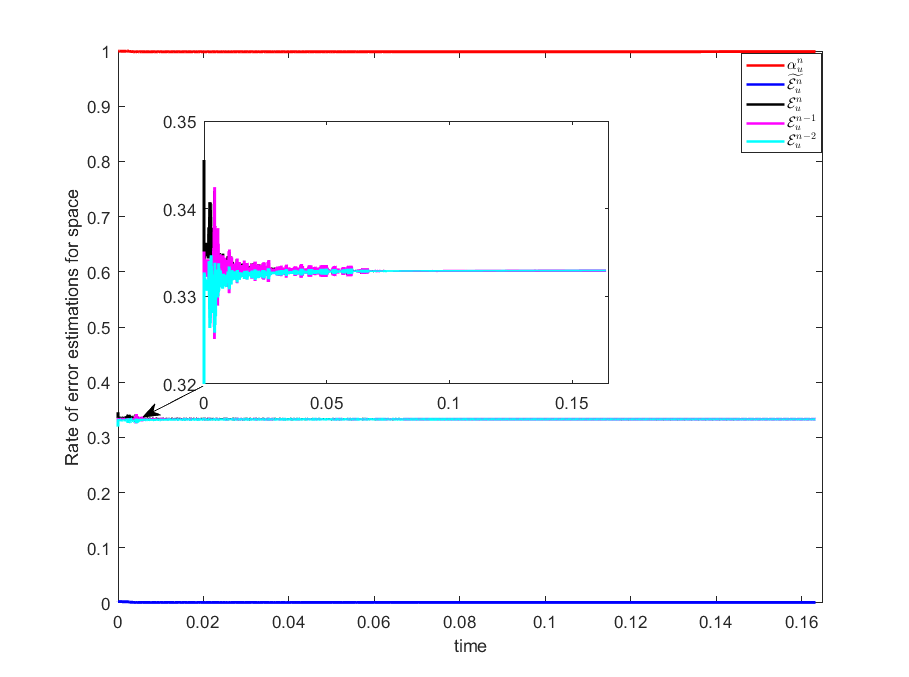}
\includegraphics[width=5.5cm,height=4.5cm]{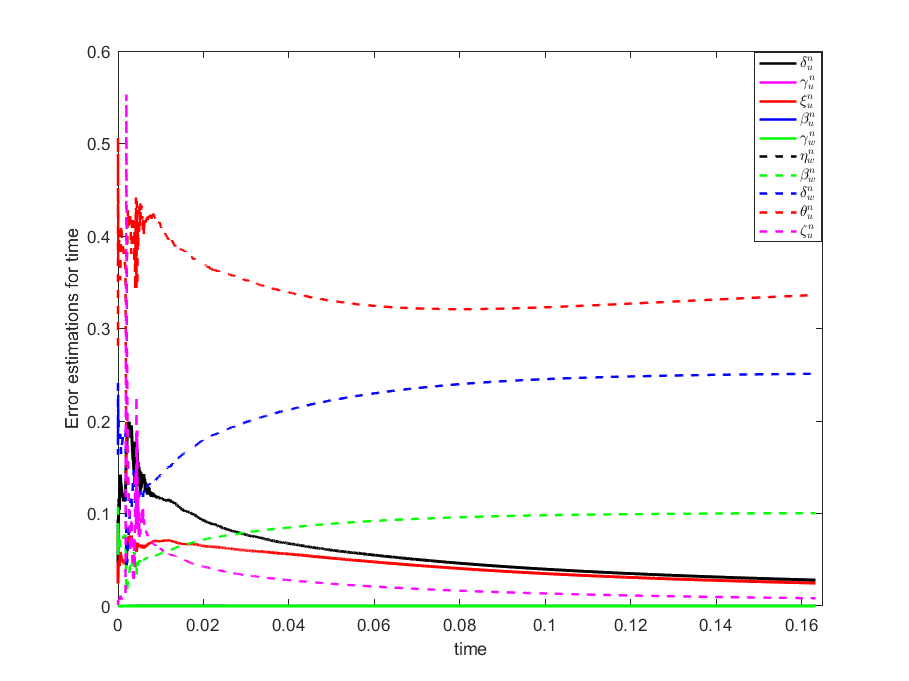}
\includegraphics[width=5.5cm,height=4.5cm]{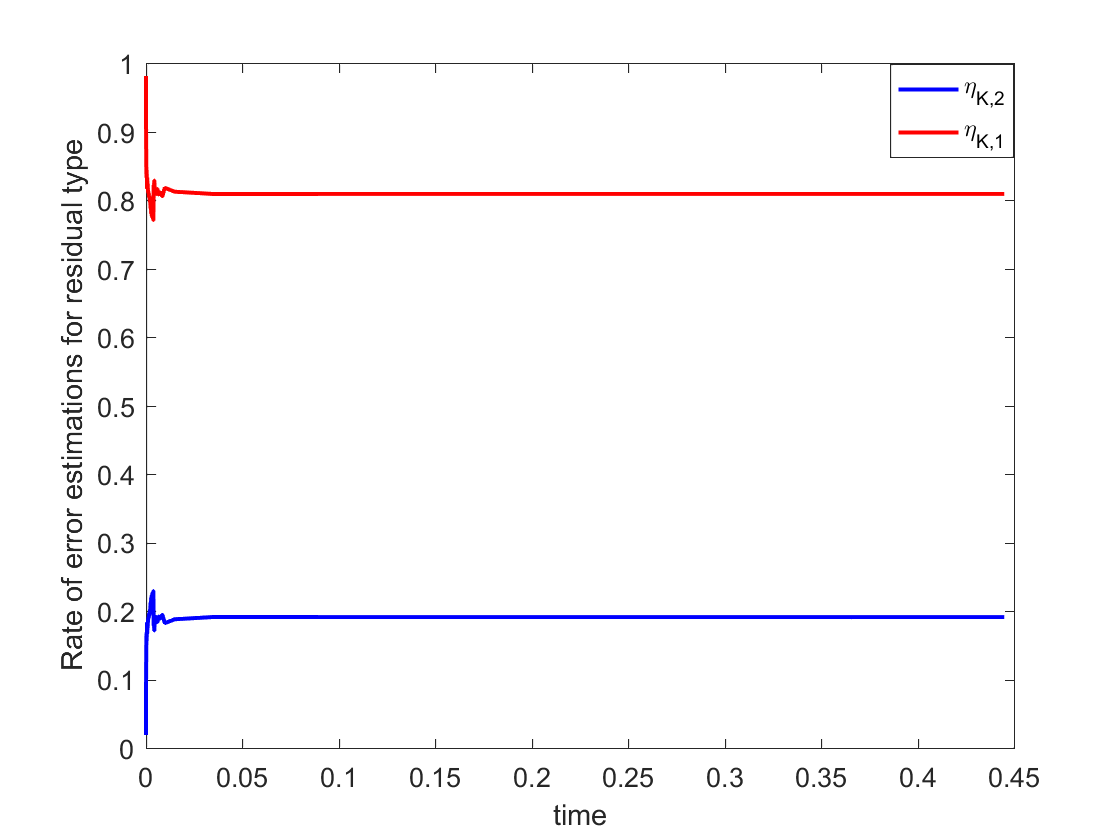}
\end{array}$
\caption{Example\ \ref{exm1}. Left: Error indicators for spatial discretization of recovery type; Middle: Error indicators for time discretization of recovery type; Right: Error indicators for spatial discretization of residual type.}\label{exp1u}
\end{figure}
%

From Theorem \ref{thm3e2}, the proposed error estimator contains twelve terms, 
\[\begin{aligned}
\eta^{n}_{time}&=\gamma_{w}^{n}+\beta_{u}^{n}+\eta_{w}^{n}+\delta_{w}^{n}+\delta_{u}^{n}+\gamma_{u}^{n}+\xi_{u}^{n}+\beta_{w}^{n}+\theta_{u}^{n}+\zeta_{u}^{n},\\
\eta^{n}_{space}&=\widetilde{\mathcal{E}_{u}^{n}}+\alpha_{u}^{n}.
\end{aligned}\]
We numerically investigate which terms are the main part of the time and space discretization error indicators. We also test the performance of the residual type error estimator provided in \cite{FW2008}, in which the local error estimators are defined by
\bq\label{4e1}
\begin{aligned}
\eta_{K,j}(t)=h_{K}\|R_{K,j}\|_{L^{2}(K)}+\sum_{\tau\in\partial K}\Big(\frac{1}{2}h_{\tau}\|J_{\tau,j}\|_{L^{2}(\tau)}^{2}\Big)^{\frac{1}{2}}, \quad j=1,2,
\end{aligned}
\eq
with the element residual
\bq\label{4e2}
\begin{aligned}
R_{K,1}&=u_{h,t}|_{K}+\mathcal{A}( w_{h}(t)|_{K}),\\
R_{K,2}&=\mathcal{A}( u_{h}(t)|_{K})+\frac{1}{\varepsilon^{2}}f(u_{h}(t)|_{K})-\frac{1}{\varepsilon}w_{h}(t)|_{K},
\end{aligned}
\eq
and the residual jumps across $\tau$
\bq\label{4e3}
\begin{aligned}
J_{\tau,1}(t)&=\Big(\nabla w_{h}(t)|_{K_{1}}-\nabla w_{h}(t)|_{K_{2}}\Big)\cdot \mathbf{n},\\
J_{\tau,2}(t)&=\Big(\nabla u_{h}(t)|_{K_{1}}-\nabla u_{h}(t)|_{K_{2}}\Big)\cdot \mathbf{n},
\end{aligned}
\eq
here $\mathbf{n}$ is the unit normal vector to $\tau$ pointing from $K_1$ to $K_2$. The corresponding total 
spatial discretization error estimator is taken as 
\bq\label{4e5}
\eta(t)=\Big(\sum_{K\in\mathcal{T}_{h}}\big(\eta_{K,1}^2(t)+\eta_{K,2}^2(t)\big)\Big)^\frac{1}{2}.
\eq

Figure \ref{exp1u} plots each parts of the error indicators. It shows that: i) for the recovery type error indicator, the time discretization error estiamtor $\eta_{time}^n$ is dominated by $\theta_u^n$, and the space discretization error estimator $\eta_{space}^n$ is dominated by $\mathcal{E}_{u}^{n}$; ii) $\eta_{K,1}(t^n)$ is the main part of the residual type error indicator $\eta(t^n)$. In the following examples, we adopt $\theta_{u}^{n}$ as the time discretization error indicator, and $\mathcal{E}_{u}^{n}$ or $\eta_{K,1}(t^n)$ as the spatial discretization error indicator, respectively.
 
\begin{example}\label{cexm2} Consider the model equation \eqref{1e3} with the parameters $\Omega=[-1,1]^{2},\, \varepsilon=0.01$,  $TOL_{t}=50$, $TOL_{t,m}=5$, $TOL_{s}=4$,  $TOL_{i}=0.002$ and the initial condition
\[\begin{aligned}
u_{0}(x,y)=&\tanh\Big(\big((x-0.3)^{2}+y^{2}-0.2^{2}\big)/\varepsilon\Big)\tanh\Big(\big((x+0.3)^{2}+y^{2}-0.2^{2}\big)/\varepsilon\Big)\times\\
&\tanh\Big(\big(x^{2}+(y-0.3)^{2}-0.2^{2}\big)/\varepsilon\Big)\tanh\Big(\big(x^{2}+(y+0.3)^{2}-0.2^{2}\big)/\varepsilon\Big).
\end{aligned}\]
\end{example}

\begin{figure}[!htbp]
$\begin{array}{ccc}
\includegraphics[width=5.5cm,height=4.5cm]{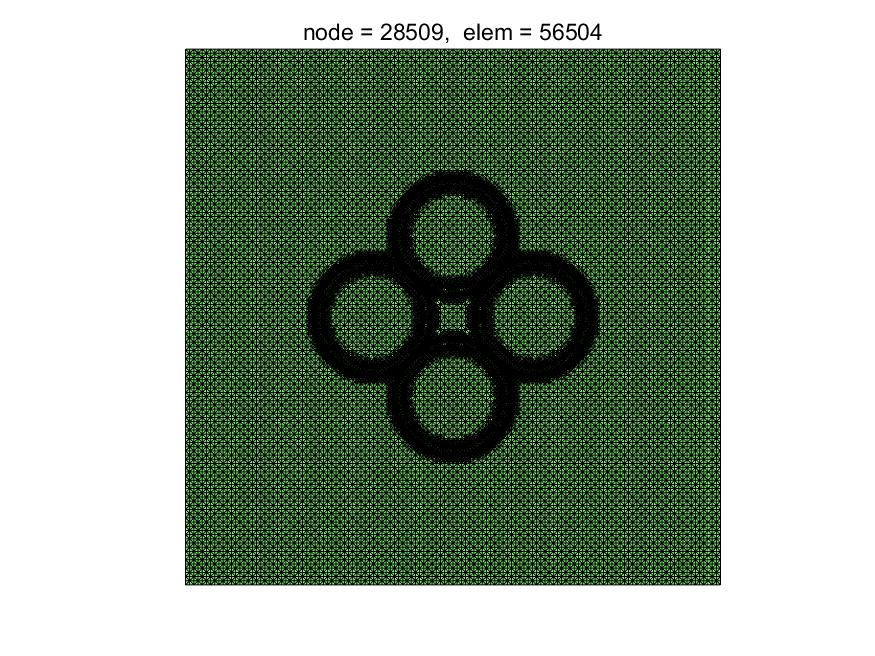}&
\includegraphics[width=5.5cm,height=4.5cm]{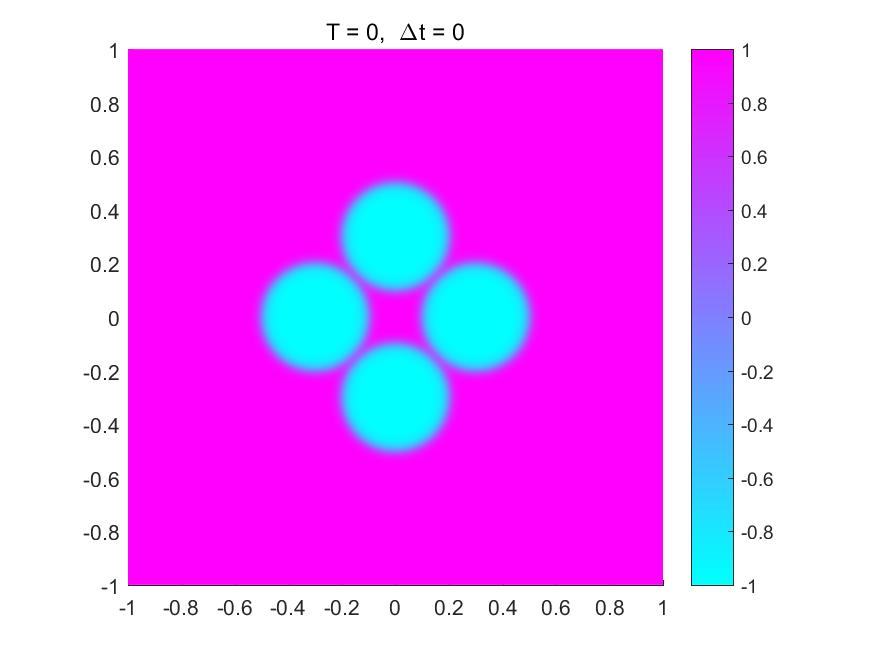}&
\includegraphics[width=5.5cm,height=4.5cm]{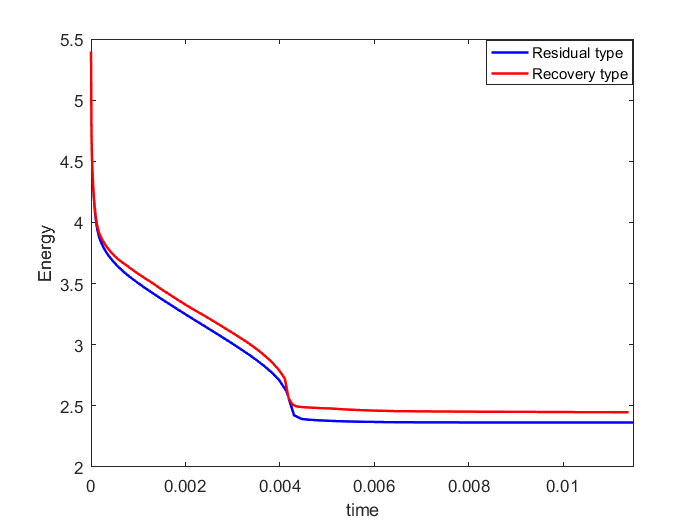}
\end{array}$\vspace{-0.2cm}
\caption{$\mathbf{Example\ \ref{cexm2}}$, Left: initial mesh; Middle: the contour plot of $u_{0}$; Right: discrete energy.}\label{Cexp2u}
\end{figure}
\begin{figure}[!htbp]
$\begin{array}{cccc}
\includegraphics[width=3.8cm,height=3.0cm]{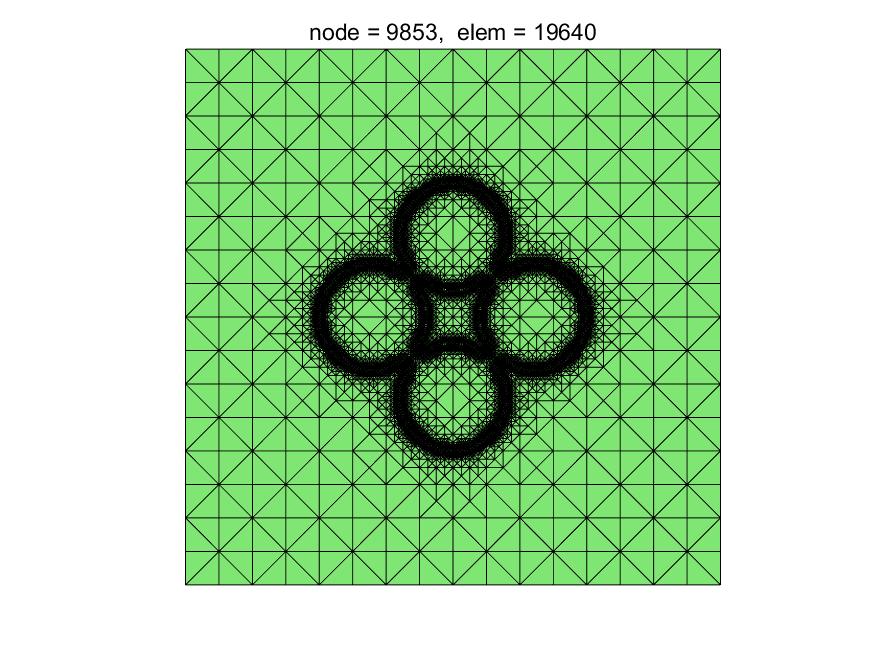}&
\includegraphics[width=3.8cm,height=3.0cm]{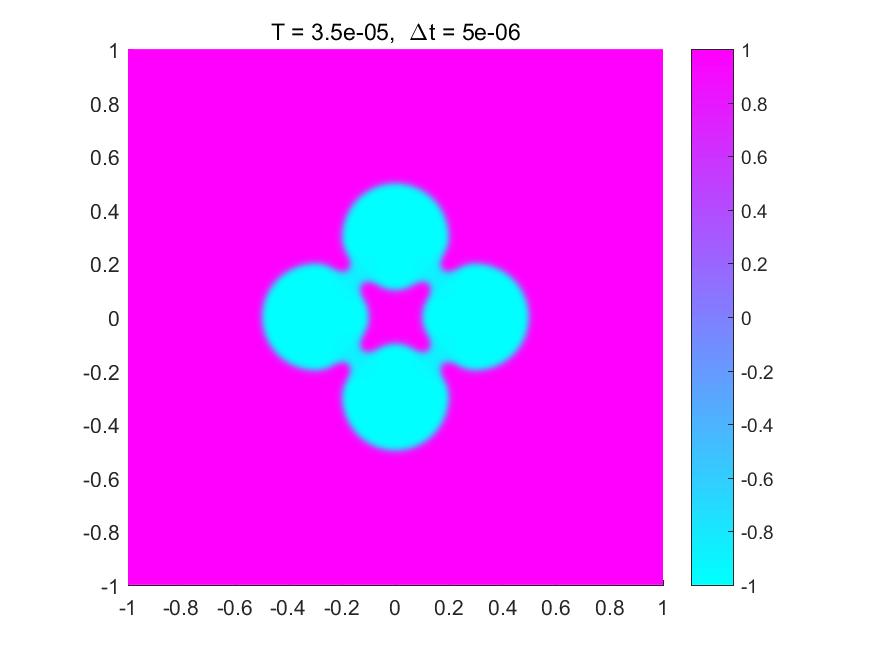}&
\includegraphics[width=3.8cm,height=3.0cm]{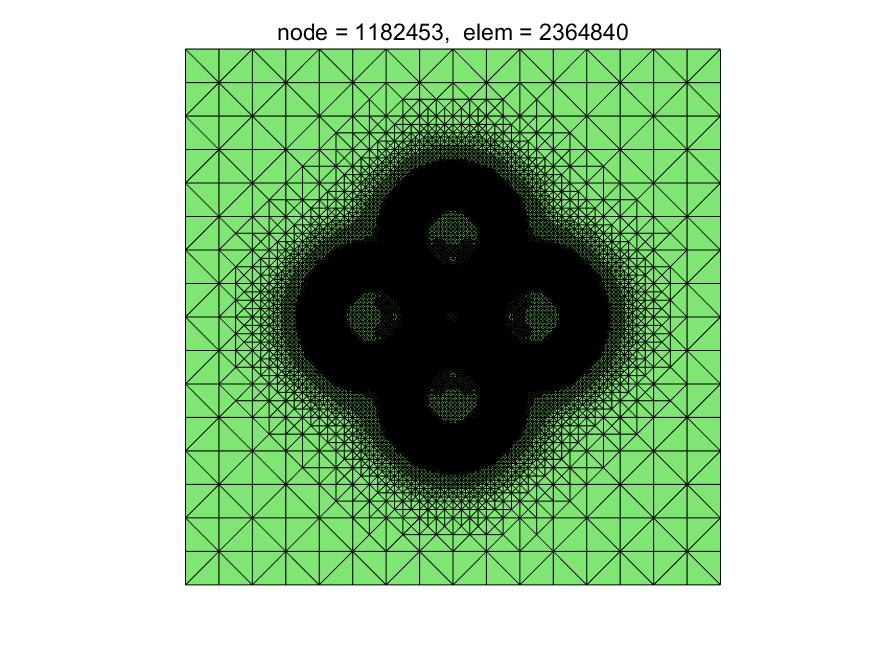}&
\includegraphics[width=3.8cm,height=3.0cm]{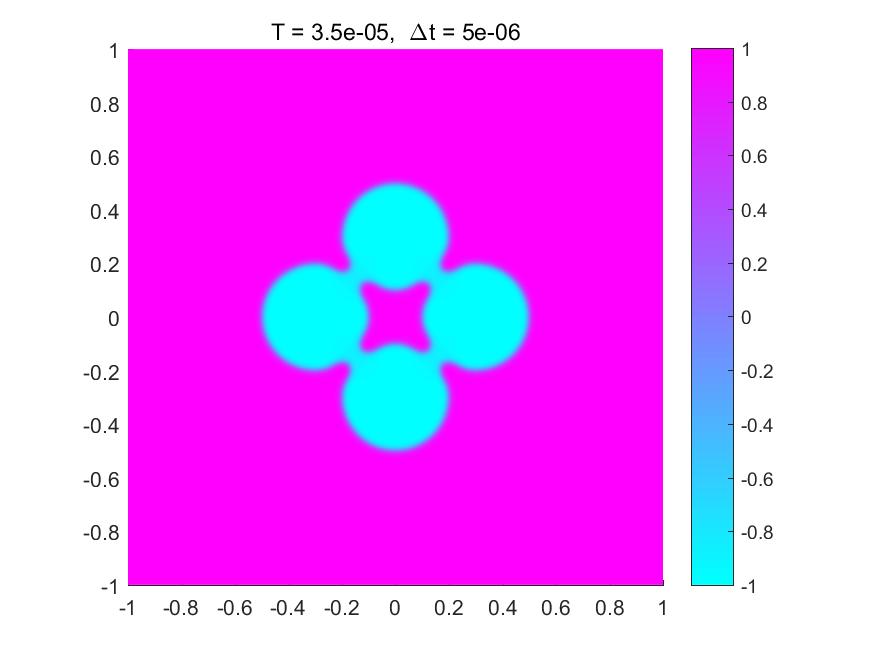}\\
\includegraphics[width=3.8cm,height=3.0cm]{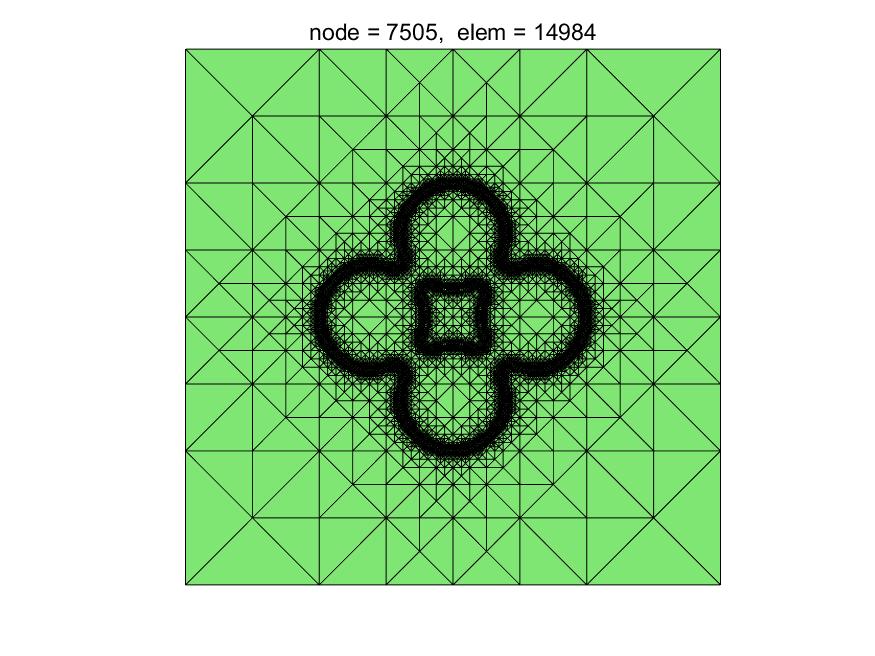}&
\includegraphics[width=3.8cm,height=3.0cm]{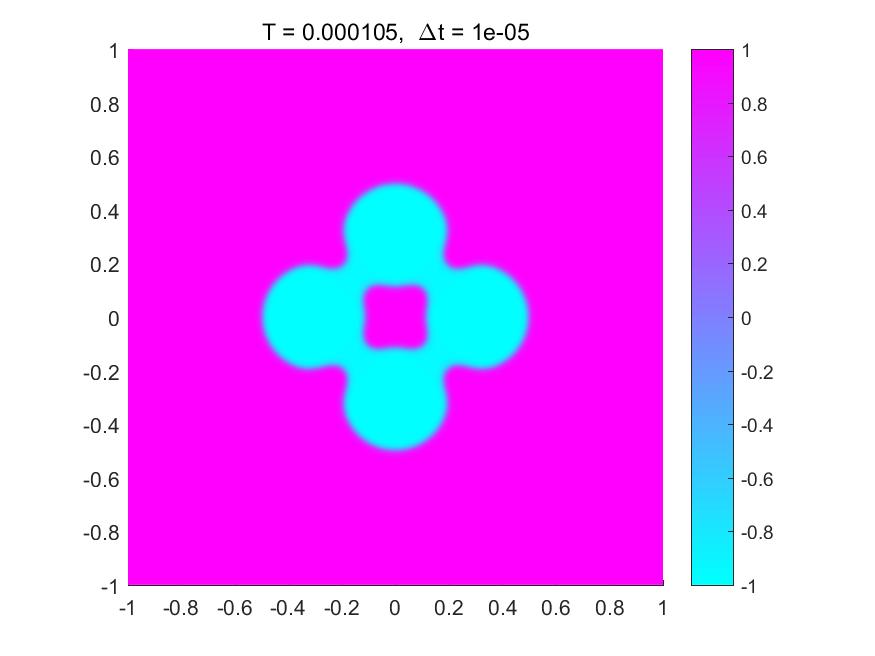}&
\includegraphics[width=3.8cm,height=3.0cm]{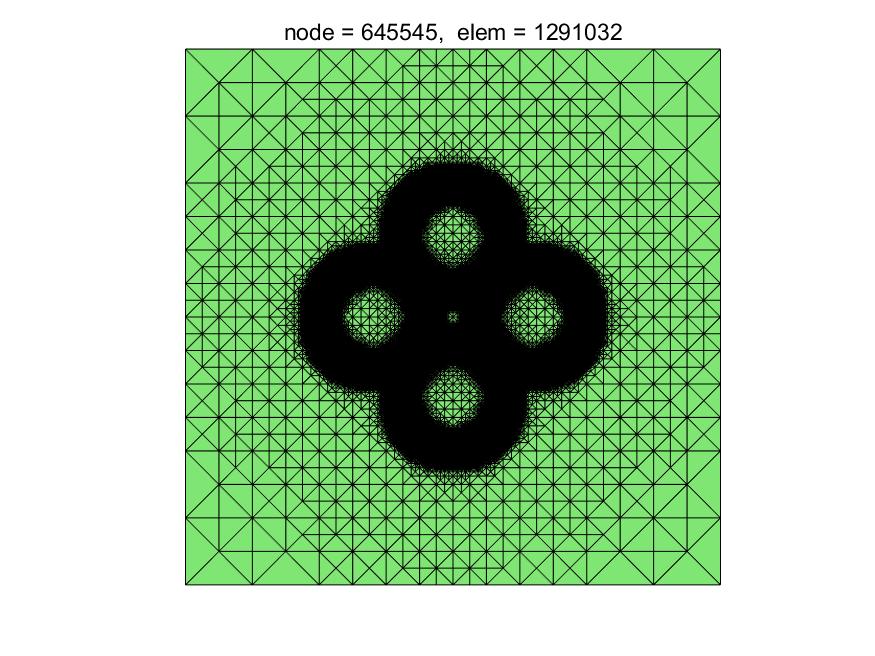}&
\includegraphics[width=3.8cm,height=3.0cm]{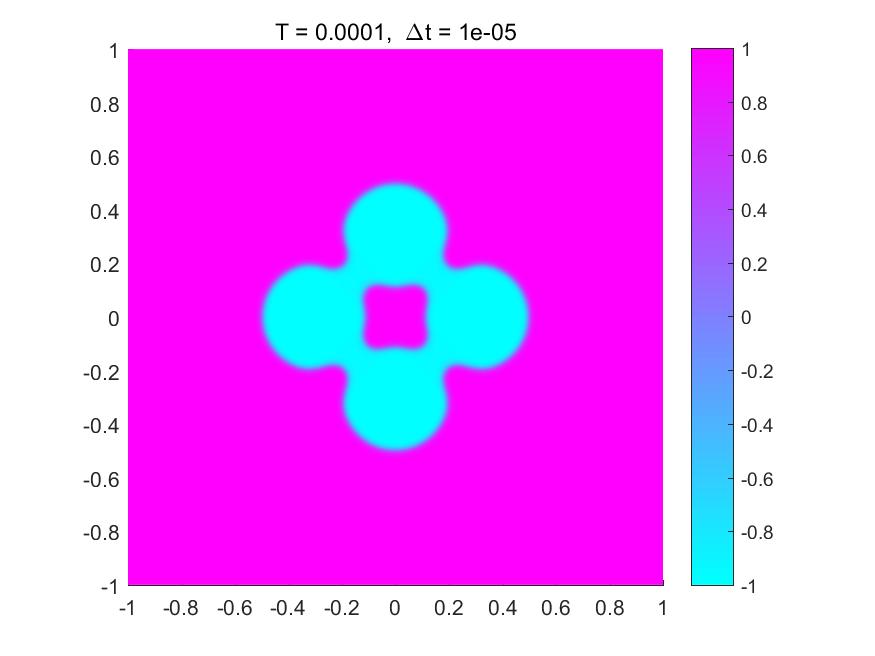}\\
\includegraphics[width=3.8cm,height=3.0cm]{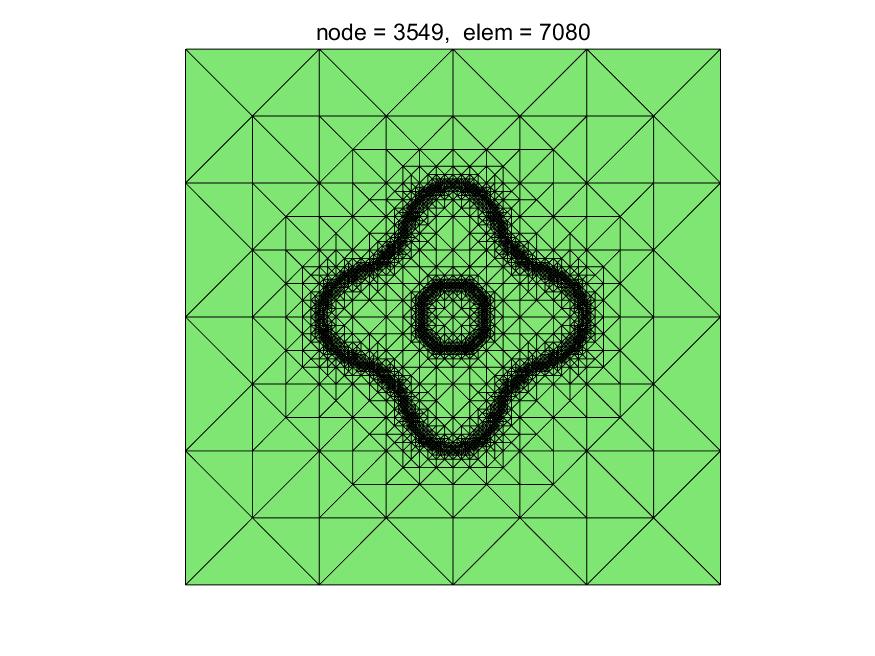}&
\includegraphics[width=3.8cm,height=3.0cm]{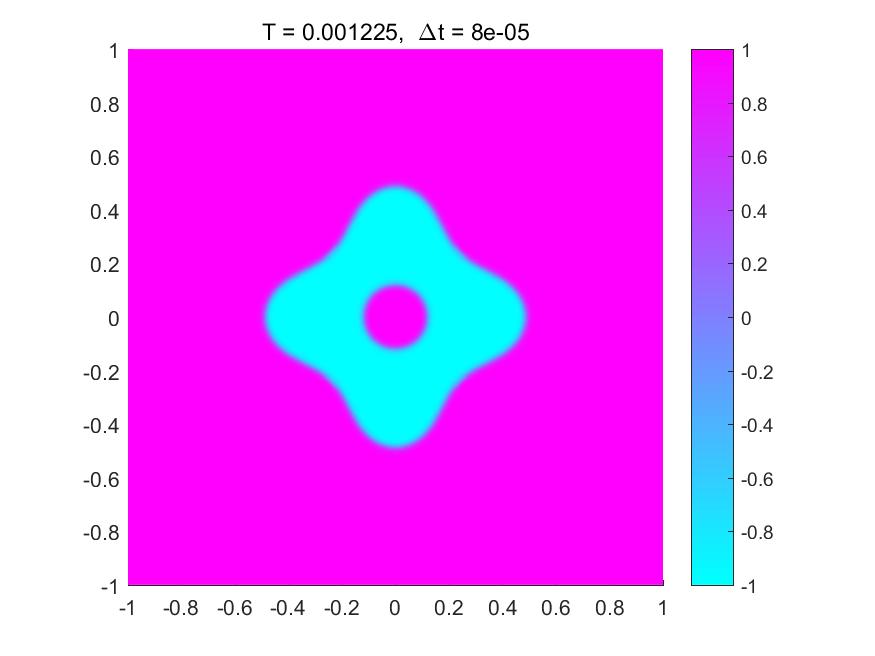}&
\includegraphics[width=3.8cm,height=3.0cm]{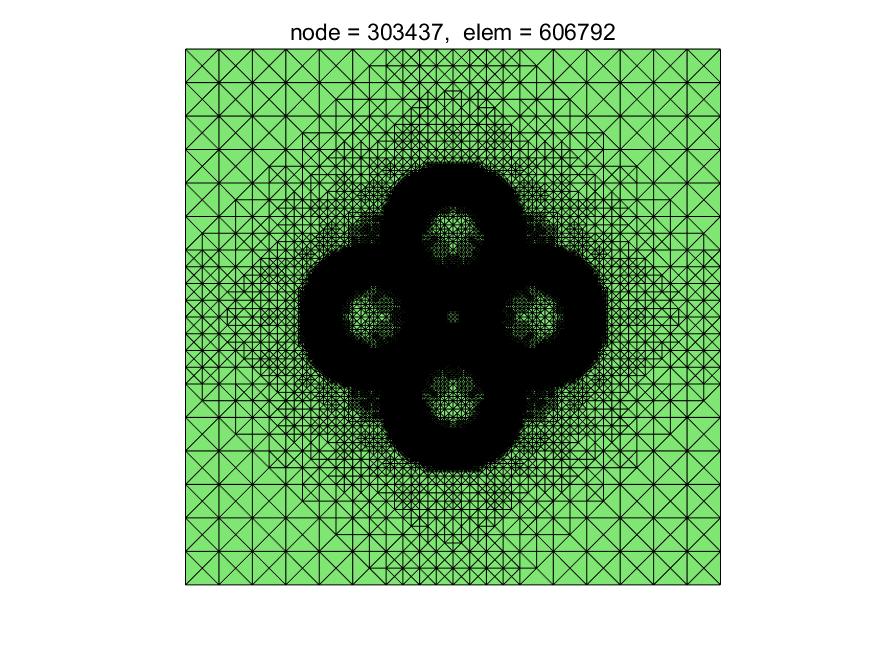}&
\includegraphics[width=3.8cm,height=3.0cm]{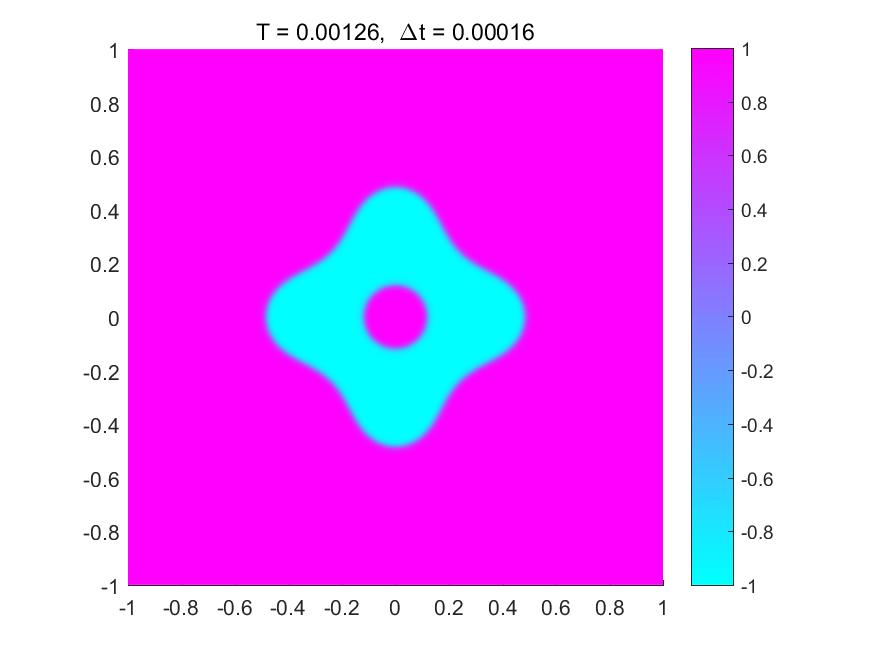}\\
\includegraphics[width=3.8cm,height=3.0cm]{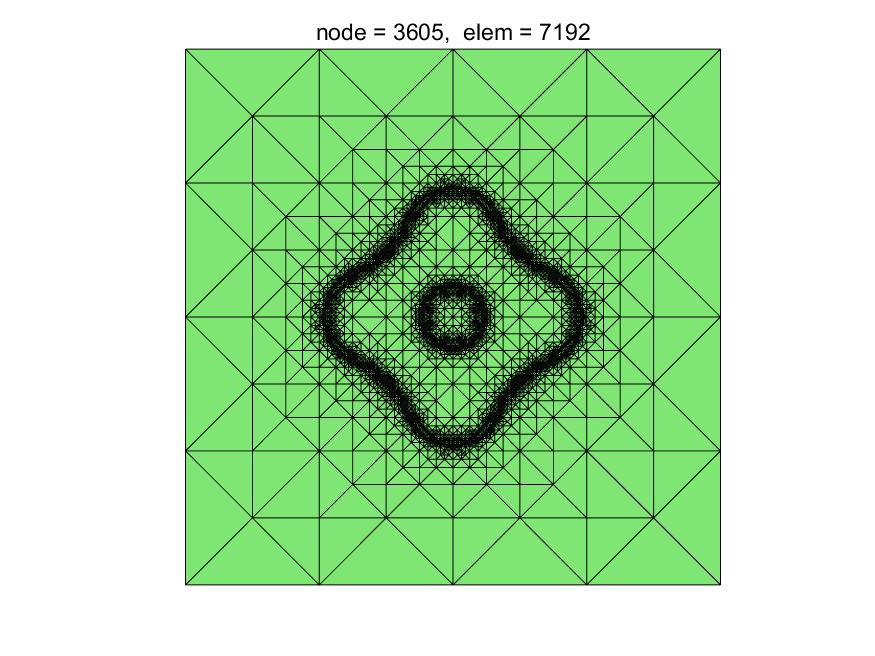}&
\includegraphics[width=3.8cm,height=3.0cm]{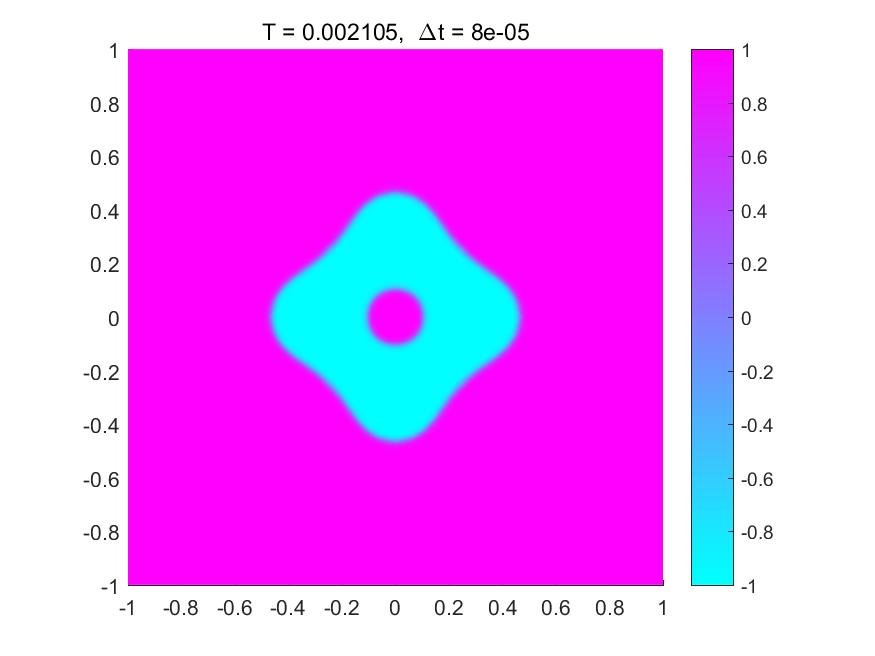}&
\includegraphics[width=3.8cm,height=3.0cm]{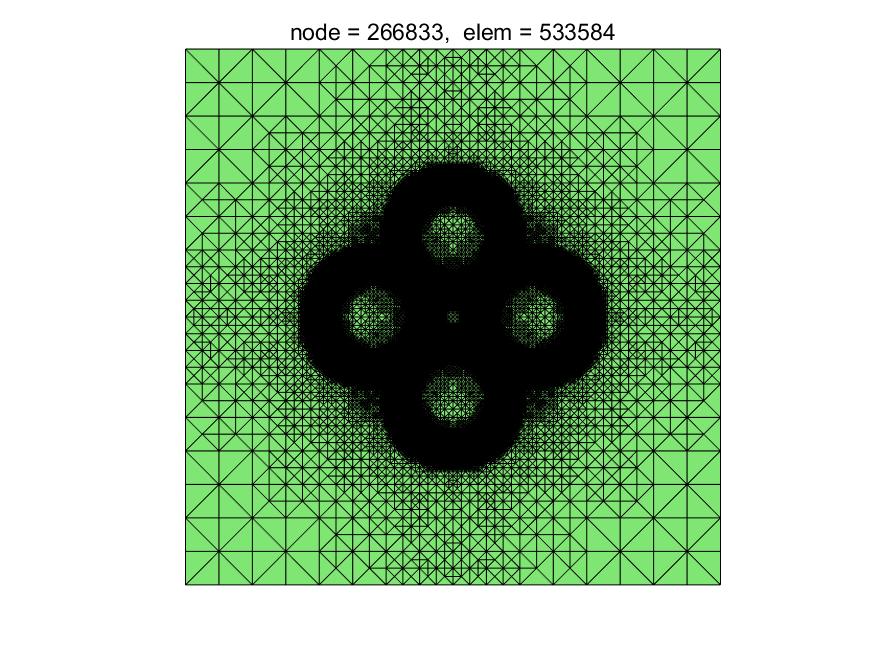}&
\includegraphics[width=3.8cm,height=3.0cm]{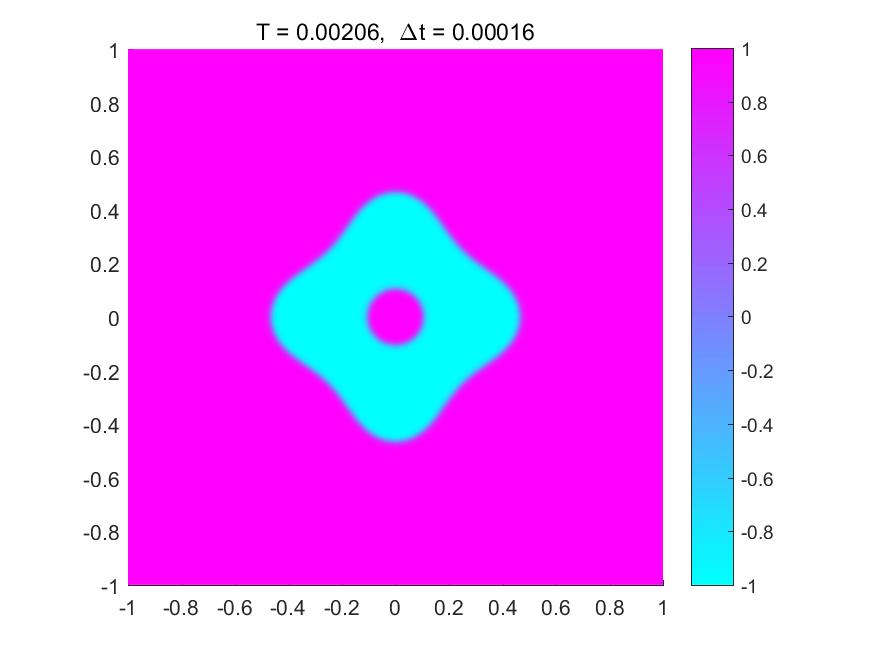}\\
\includegraphics[width=3.8cm,height=3.0cm]{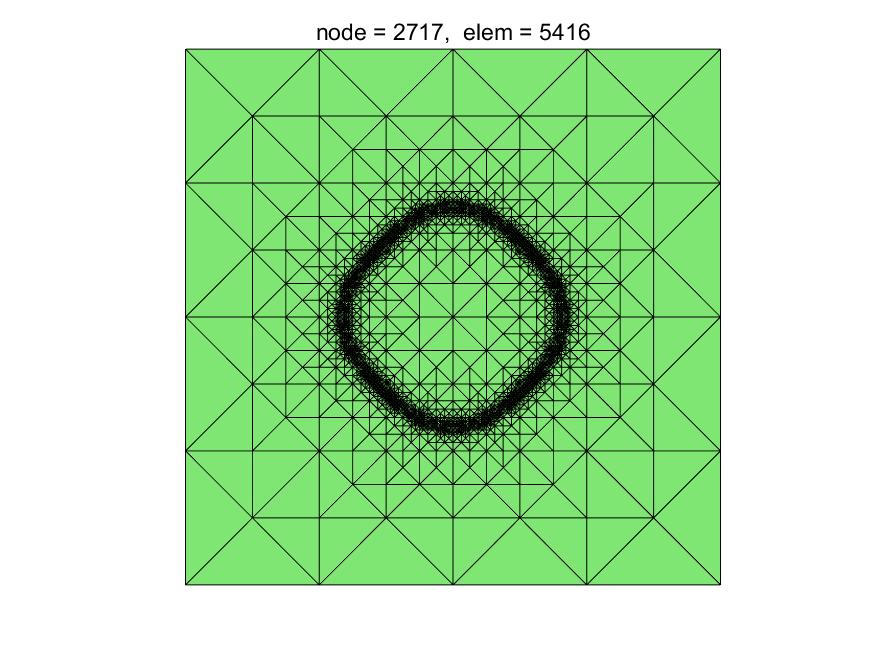}&
\includegraphics[width=3.8cm,height=3.0cm]{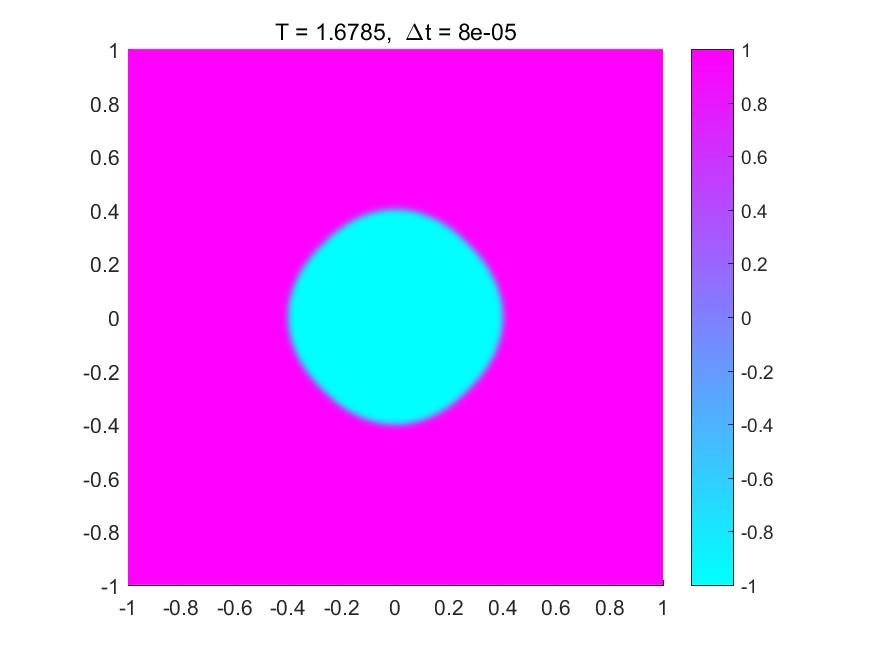}&
\includegraphics[width=3.8cm,height=3.0cm]{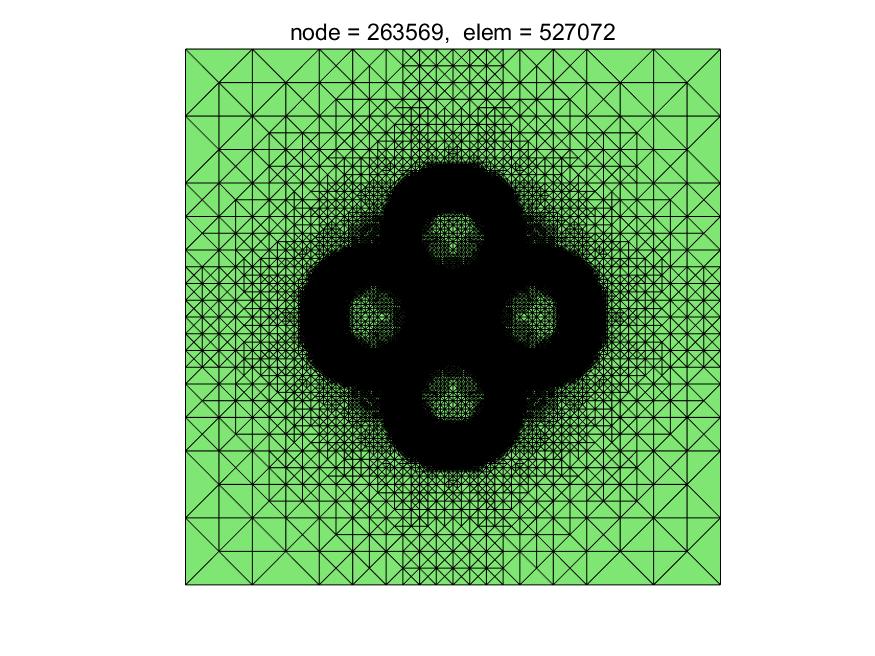}&
\includegraphics[width=3.8cm,height=3.0cm]{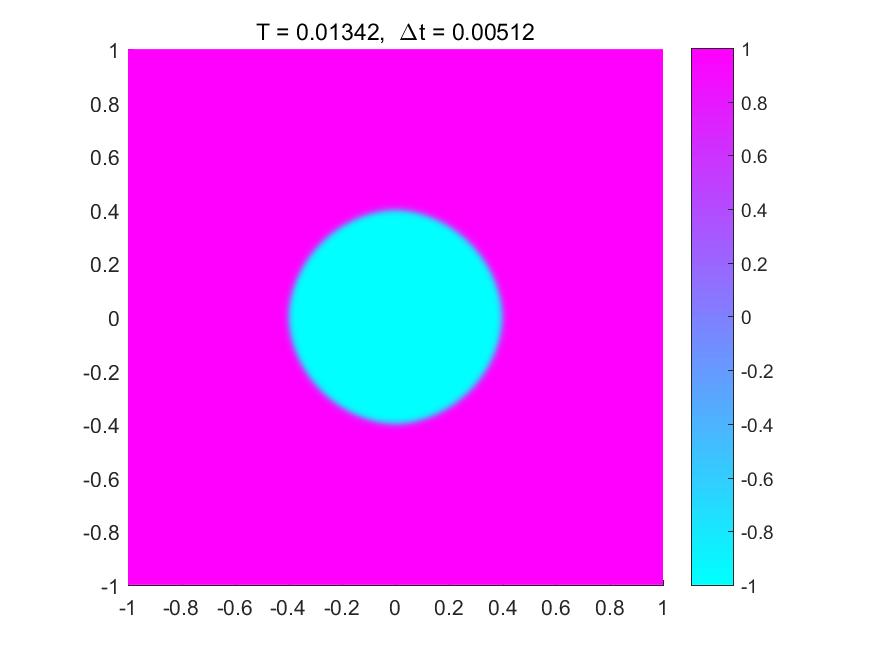}
\end{array}$
\caption{Example\ \ref{cexm2}, adaptive meshes and snapshots of numerical solutions; First and second column: recovery type; Third and fourth column:
residual type.}\label{exp2u}
\end{figure}

\begin{figure}[!htbp]
$\begin{array}{cc}
\includegraphics[width=5.5cm,height=4.5cm]{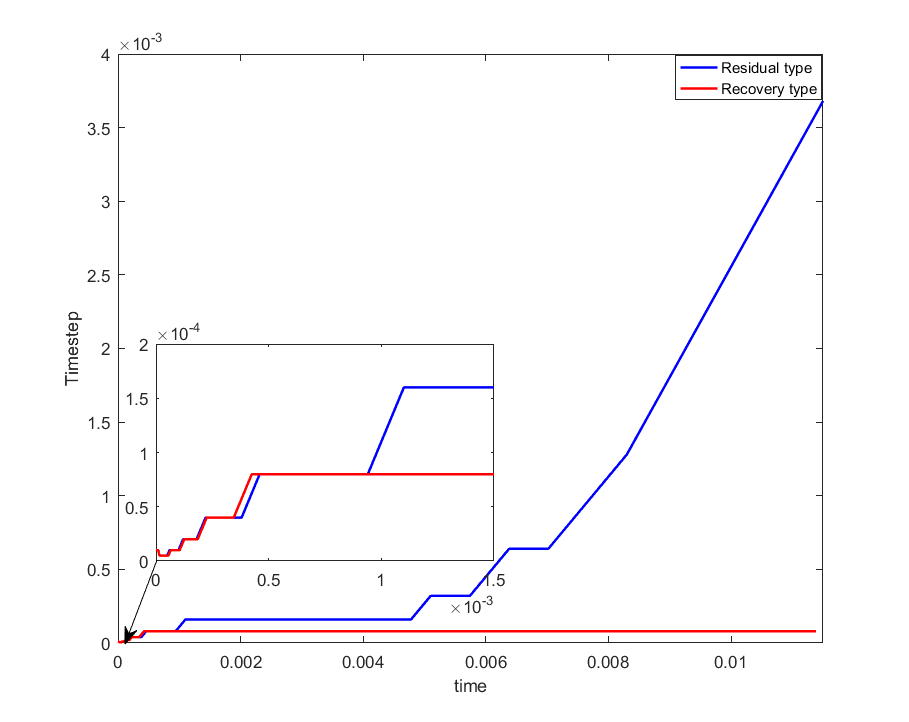}&
\includegraphics[width=5.5cm,height=4.5cm]{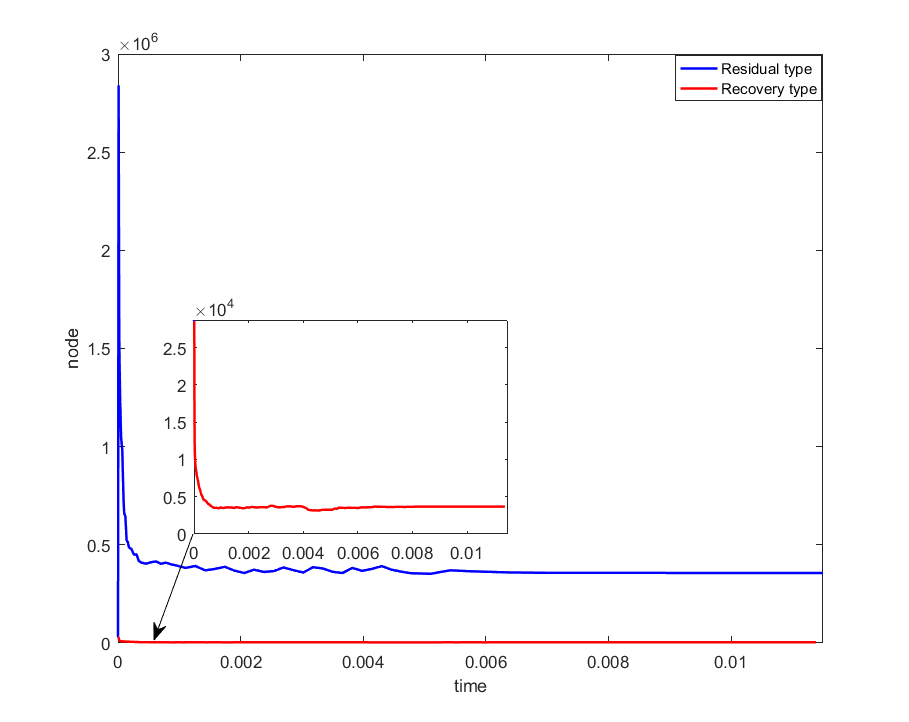}
\end{array}$\vspace{-0.2cm}
\caption{Example\ \ref{cexm2}, Left: time-steps; Right: number of nodes.}\label{Cexp2un}
\end{figure}

In this example, we compare the recovery type a posteriori error estimator with the residual type. Figure \ref{Cexp2u} displays the initial mesh, contour plot of the initial numerical solution, and discrete energy history for the two spatial error estimators based on the proposed time-space adaptive algorithm. We can see clearly that the energy decreases over time.
Figure \ref{exp2u} shows the sequences of adaptive meshes and contour plots of the corresponding approximate solutions produced by the time-space adaptive algorithm guided by the recovery and residual type error indicators for the spatial discretization, respectively. The adaptive meshes match the numerical solutions of Algorithm \ref{algadp} based on the recovery type error indicator better than the residual type. The corresponding time-step and number of nodes are also displayed in Figure \ref{Cexp2un}. We observe that as the time-step grows, the degree of freedom based on the recovery type is much less than the residual type, indicating that the recovery type a posteriori error estimation is clearly superior to the residual type.

\begin{table}[!ht]
\centering
\caption{Example\ \ref{cexm2}\ (T=0.01), CPU time for two kinds of types by using time-space adaptive algorithm and space-only adaptation, respectively (11th Gen Intel(R) Core(TM) i5-1135G7 @ 2.40GHz 2.42GHz).}\label{Table1:S1}
\begin{tabular}{| c | c | c | c | c |}
\hline
CPU time   & time-space adaptation & space-only adaptation                         \\  \hline
Recovery type   &  389s    &       1095s                   \\  \hline
Residual type     &  106975s    &       -                   \\  \hline
\end{tabular}
\end{table}


Furthermore, we evaluate the efficiency of the adaptive algorithm with time and space adaptation. Table \ref{Table1:S1} reports the corresponding CPU time. We observe that: i) the time-space adaptive method based on our proposed recovery type error estimator is significantly more efficient than the adaptive method based on the residual type error indicator; ii) the time-space adaptation is more efficient than the adaptive method with space-only adaptation.

\begin{example}\label{cexm3} In the last example, we consider the three dimensional Cahn--Hilliard equation \eqref{1e3} with the following initial condition
\[\begin{aligned}
u_{0}(x,y,z)=\varepsilon\cos(1.5\pi x)\cos(1.5\pi y)\big(\sin(\pi z)+\sin(2\pi z)\big),
\end{aligned}\]
where $\Omega=[-1,1]^{3}$ and the parameters $\varepsilon=0.05$, $TOL_{t}=20$, $TOL_{t,m}=1$, $TOL_{s}=1.5$,  $TOL_{i}=8e-5$.
\end{example}

\begin{figure}[!htbp]
$\begin{array}{ccc}
\includegraphics[width=5.5cm,height=4.5cm]{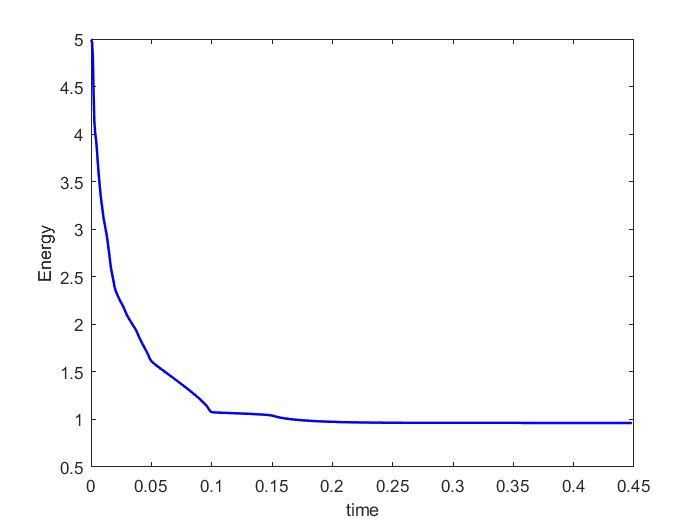}&
\includegraphics[width=5.5cm,height=4.5cm]{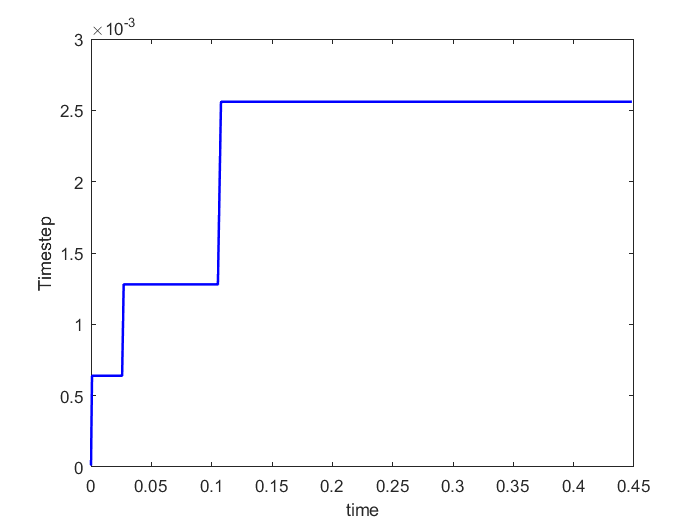}&
\includegraphics[width=5.5cm,height=4.5cm]{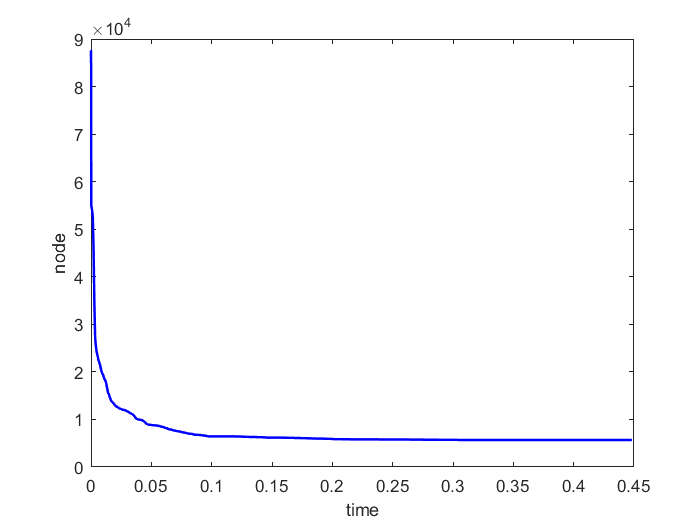}
\end{array}$\vspace{-0.2cm}
\caption{$\mathbf{Example\ \ref{cexm3}}$, Left: discrete energy; Middle: time-steps; Right: number of nodes.}\label{Cexp3u}
\end{figure}

\begin{figure}[!htbp]
$\begin{array}{cccc}
\includegraphics[width=3.8cm,height=3.0cm]{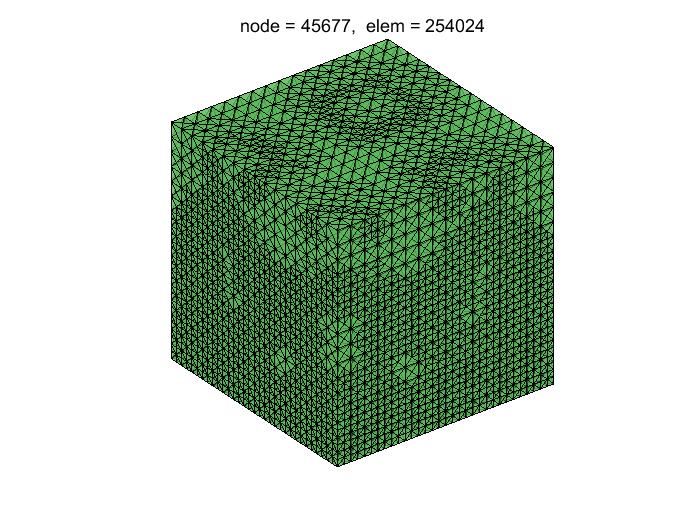}&
\includegraphics[width=3.8cm,height=3.0cm]{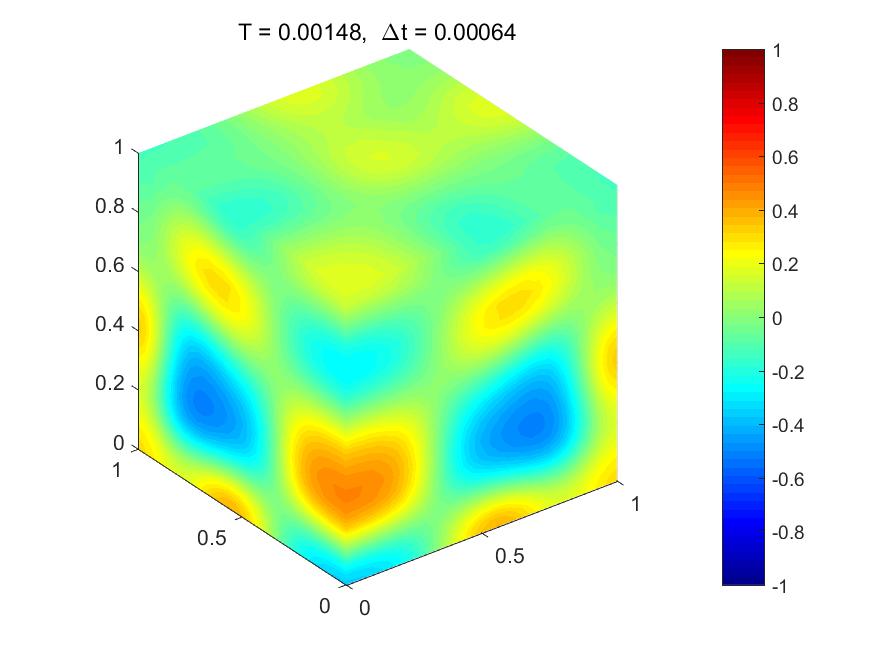}&
\includegraphics[width=3.8cm,height=3.0cm]{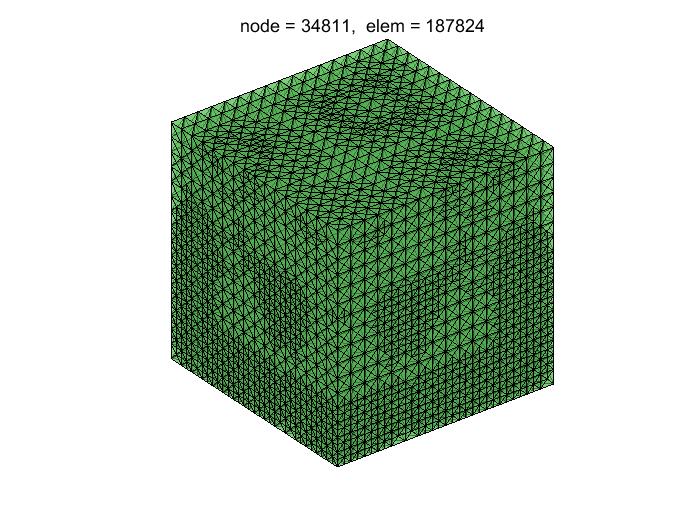}&
\includegraphics[width=3.8cm,height=3.0cm]{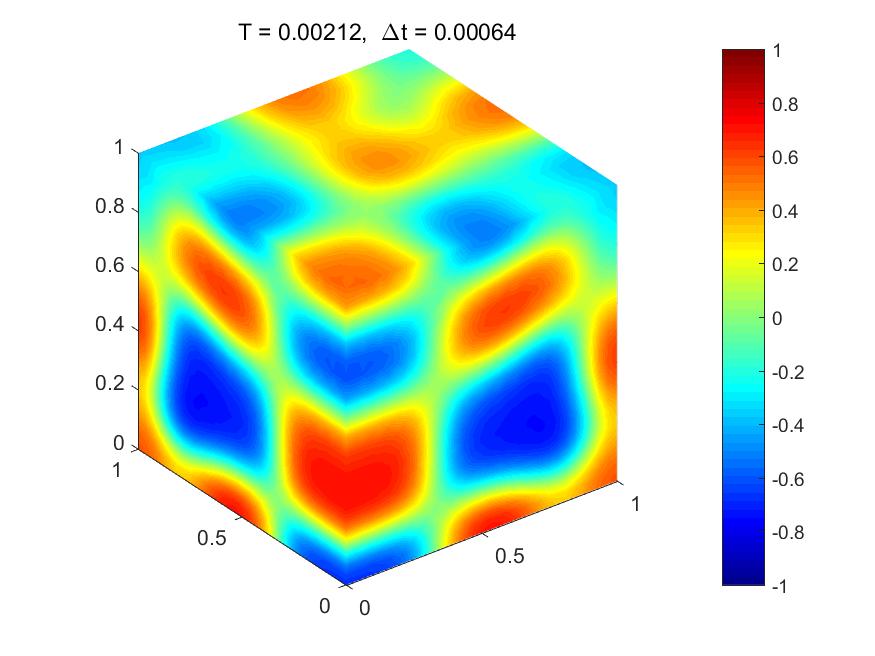}\\
\includegraphics[width=3.8cm,height=3.0cm]{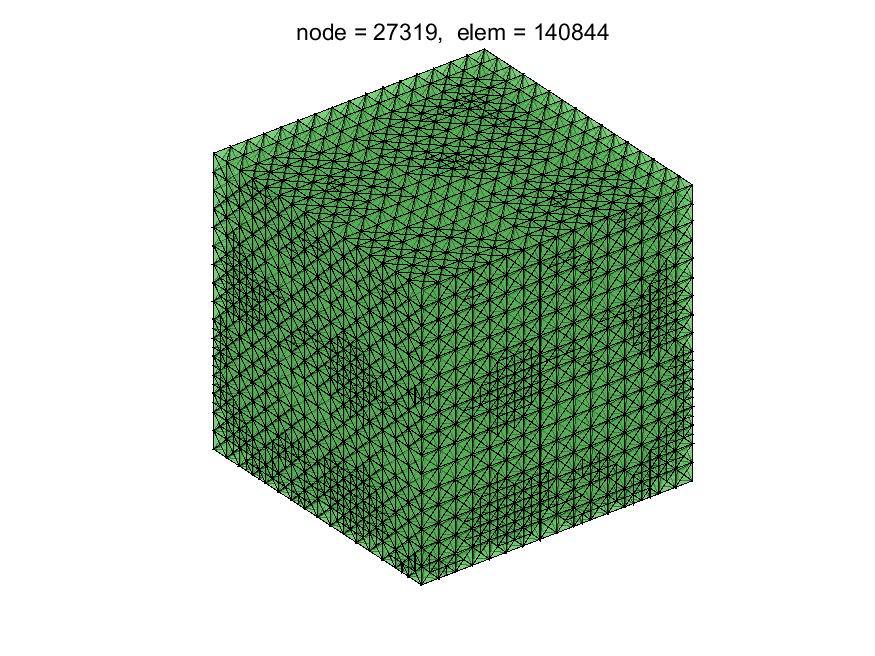}&
\includegraphics[width=3.8cm,height=3.0cm]{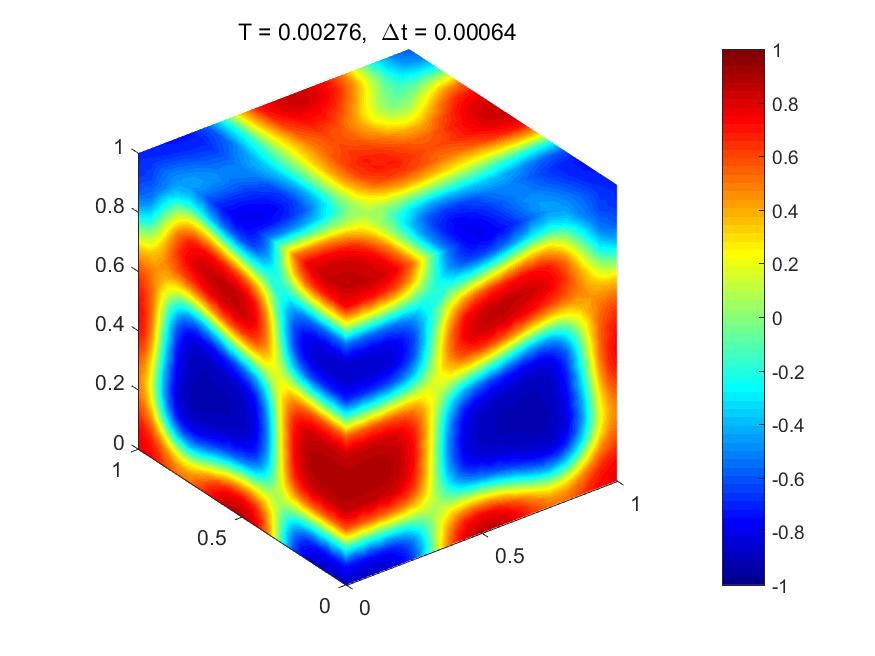}&
\includegraphics[width=3.8cm,height=3.0cm]{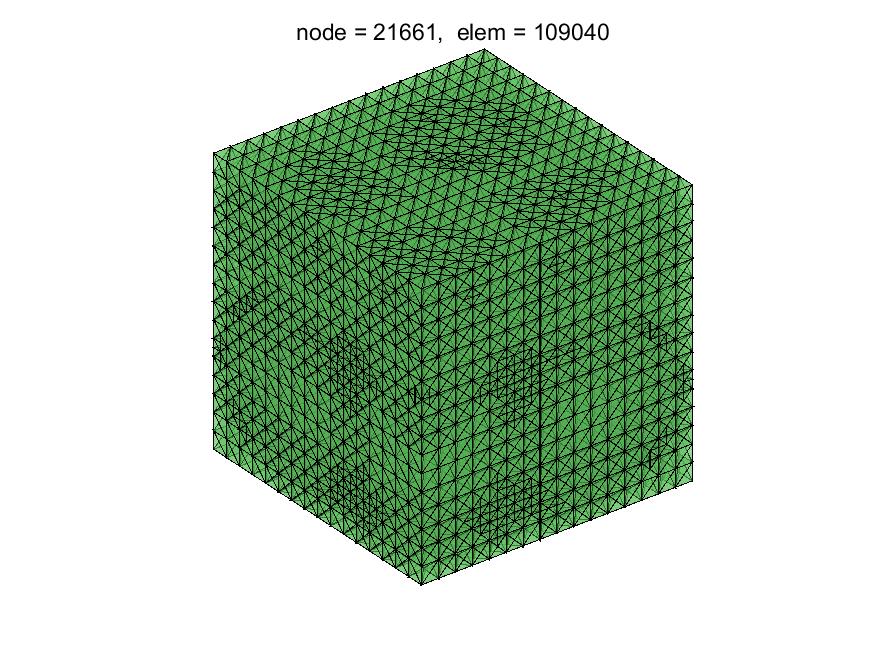}&
\includegraphics[width=3.8cm,height=3.0cm]{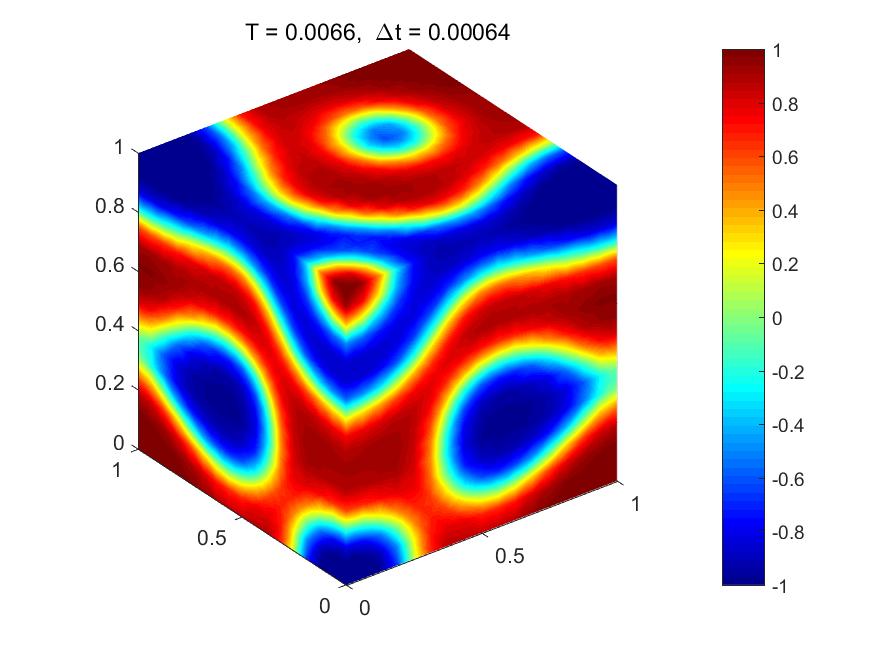}\\
\includegraphics[width=3.8cm,height=3.0cm]{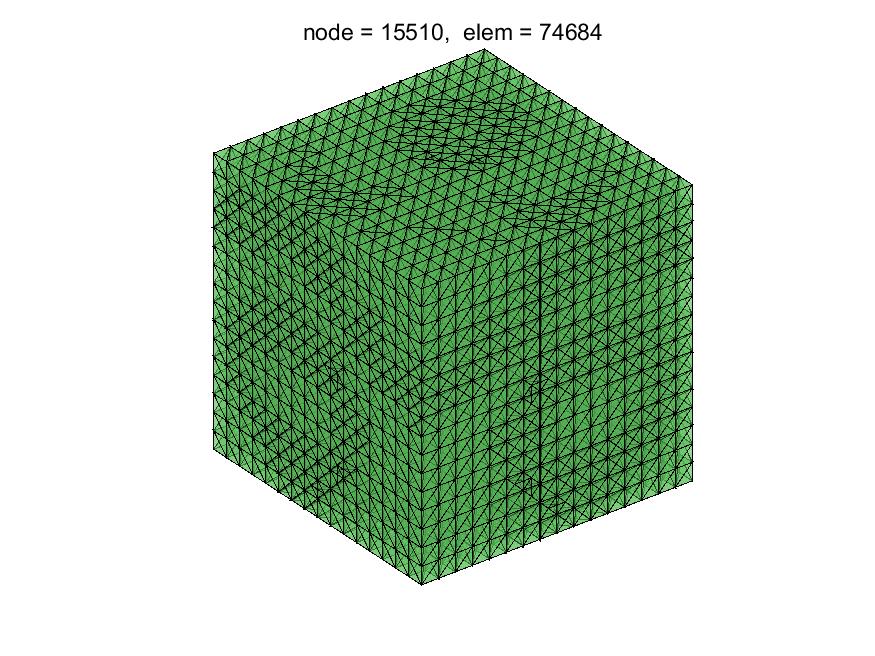}&
\includegraphics[width=3.8cm,height=3.0cm]{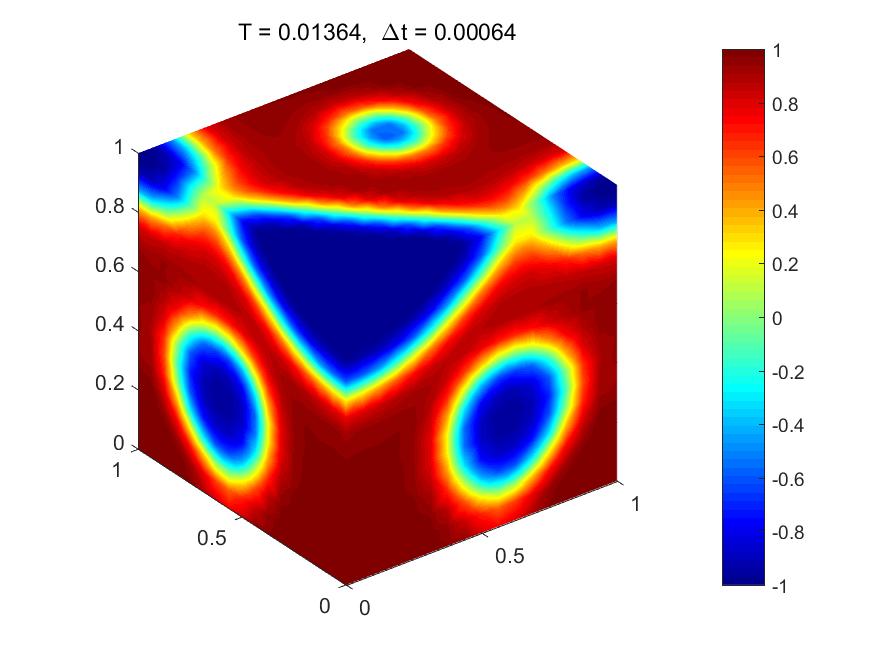}&
\includegraphics[width=3.8cm,height=3.0cm]{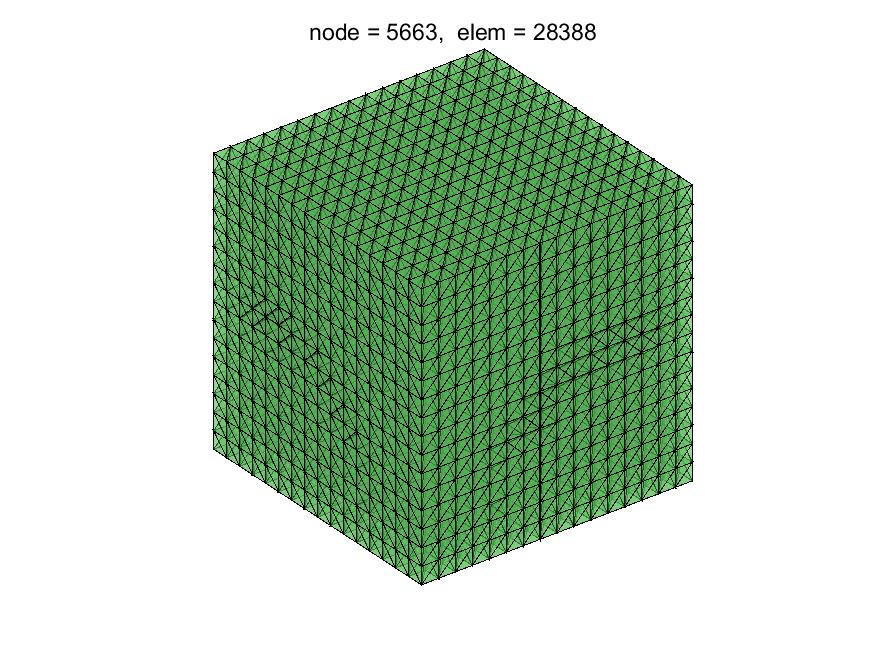}&
\includegraphics[width=3.8cm,height=3.0cm]{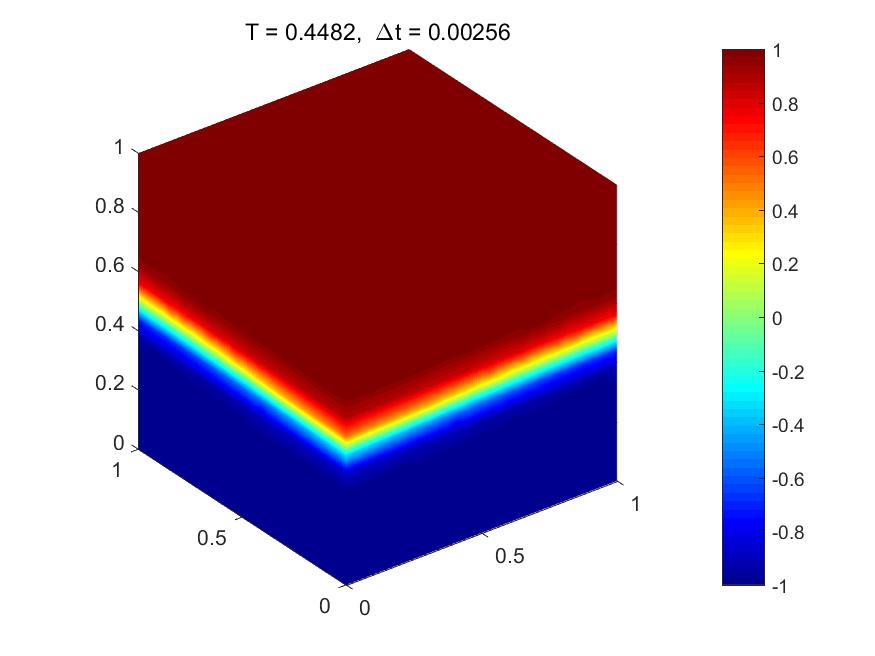}
\end{array}$
\caption{Example\ \ref{cexm3}, adaptive meshes and snapshots of numerical solutions.}\label{exp3u}
\end{figure}


Figure \ref{Cexp3u} displays the contour plots of the discrete energy history, time steps, and the change in the number of nodes with time. It is evident that the energy and the number of nodes both decrease over time, and the time steps change with time. In Figure \ref{exp3u}, we show the sequence of adaptive meshes and contour plots of the corresponding approximate solutions. We observe that the meshes adapt around the zero level set, which confirms the effectiveness of the derived a posteriori error estimation and adaptive algorithm for the three-dimensional Cahn--Hilliard equation.

\hskip\parindent
\section{Conclusions}\label{secCon} 
\setcounter{equation}{0}

In this paper, we derived a novel SCR-based recovery type a posteriori error estimator for the Crank-Nicolson finite element method applied to the Cahn--Hilliard equation. 
The derivation of the error estimator utilized the elliptic reconstruction technique and the time reconstruction technique, which was based on approximations on three time levels and led to a second order error estimator for time discretization.  Based on the derived a posteriori error estimator, we designed an efficient time-space adaptive algorithm. The numerical results indicated that the recovery-type a posteriori error estimator and the time-space adaptive strategy could greatly improve the efficiency of the adaptive algorithm for the Cahn--Hilliard equation. 
Notably, our proposed time-space adaptive finite element method outperformed the adaptive finite element method based on residual-type a posteriori error estimators, as well as the space-only adaptive finite element method. These results demonstrate the superior efficiency of our method in accurately solving the Chan--Hilliard equation at hand.



\hskip\parindent

\appendix

\section{Proof of Theorem \ref{thm3e2}}\label{aprofthm}

In this section, we present the proof of the Theorem \ref{thm3e2}. To begin with, we recall the following results.

\begin{lem}\label{lem3e1}\cite{BM2011}
Let $\dot{V}:=\left\{\phi\in H^{1}(\Omega),\bar{\phi}:=\frac{1}{|\Omega|}\int_{\Omega}\phi dx=0\right\},$ there exists $C_{I}>0$ such that for all $\phi\in \dot{V}$ if $d=2$ and for all $\phi\in \dot{V}\cap L^{\infty}(\Omega)$ if $d=3$, we have
\bq\label{propA1}
\|\phi\|_{L^{3}(\Omega)}^{3}\leq C_{I}\|\phi\|_{L^{\infty}(\Omega)}^{1-\sigma}\|\nabla\Delta^{-1}\phi\|^{\sigma}\|\nabla\phi\|^{2},
\eq
where $\sigma=1$ if $d=2$ and $\sigma=\frac{4}{5}$ if $d=3$.
\end{lem}

\begin{lem}\label{lem3e2}\cite{BMO2011} (Generalized Gronwall's Lemma) 
\ Suppose that the nonnegative functions $y_{1}\in C([0,T])$, $y_{2},\,y_{3}\in L^{1}(0,T)$, $a\in L^{\infty}(0,T)$, and the real number $A\geq 0$ satisfy
 \[y_{1}(t)+\int_{0}^{t}y_{2}(s)ds\leq A+\int_{0}^{t}a(s)y_{1}(s)ds+\int_{0}^{t}y_{3}(s)ds\]
for all $t\in [0,T]$. Assume that for $B\geq 0$, $\beta\geq 0$ and every $t\in [0,T]$, we have
\[\int_{0}^{t}y_{3}(s)ds\leq B\sup_{s\in [0,t]}y_{1}^{\beta}(s)\int_{0}^{t}\left(y_{1}(s)+y_{2}(s)\right)ds.\]
Setting $E:=exp\left(\int_{0}^{T}a(s)ds\right)$ and assume that $8AE\leq\left(8B(1+T)E\right)^{-1/\beta}$, then we obtain
\[\sup_{t\in [0,T]}y_{1}(t)+\int_{0}^{T}y_{2}(s)ds\leq 8Aexp\left(\int_{0}^{T}a(s)ds\right).\]
\end{lem}


Now, we are ready to present the proof of Theorem \ref{thm3e2}.
\begin{proof}
According to \eqref{2e9}, we have that
\begin{flalign}
\partial_{t}e_u+\mathcal{A}e_{w}
=&-\mathcal{A}\epsilon_{w}+\frac{A^{n}w_{h}^{n}-A^{n-1} w_{h}^{n-1}}{2}+\mathcal{A}(q(t)-q^{n})+w_{h}^{n}-R^{n}w_{h}^{n}\nonumber\\&+(t-t_{n-\frac{1}{2}})\partial_{n}^{2}u_{h},\label{2e12a}\\
\varepsilon\mathcal{A}e_{u}-e_{w}
=&-\varepsilon\mathcal{A}\epsilon_{u}+\epsilon_{w}+\varepsilon\frac{A^{n}u_{h}^{n}-A^{n-1}u_{h}^{n-1}}{2}+\frac{P^{n}f(u_{h}^{n})-P^{n-1}f(u_{h}^{n-1})}{2\varepsilon}-\frac{w_{h}^{n}-w_{h}^{n-1}}{2}\nonumber\\
&+\varepsilon\mathcal{A}(p(t)-p^{n})-(q(t)-q^{n})+\frac{f(u)-f(p^{n})}{\varepsilon}+\frac{1}{\varepsilon}\left(u_{h}^{n}-p^{n}\right).\label{2e12b}
\end{flalign}
To make the conclusion clean, we separate the remaining of this proof into eleven steps.

Step $1$: Multiplying both sides of \eqref{2e12a} by $-\Delta^{-1}e_{u}$ and \eqref{2e12b} by $e_{u}$, respectively, then adding the resulting equations, we obtain
\begin{align*}
\frac{1}{2}\frac{d}{dt}&\|\nabla\Delta^{-1}e_{u}\|^{2}+\varepsilon\|\nabla e_{u}\|^{2}
=\left(\frac{A^{n}w_{h}^{n}-A^{n-1} w_{h}^{n-1}}{2},-\Delta^{-1}e_{u}\right)\\
&+\left(\mathcal{A}(q(t)-q^{n}),-\Delta^{-1}e_{u}\right)+\left(w_{h}^{n}-R^{n}w_{h}^{n},-\Delta^{-1}e_{u}\right)+\left((t-t_{n-\frac{1}{2}})\partial_{n}^{2}u_{h},-\Delta^{-1}e_{u}\right)\\
&-\varepsilon a\left(\epsilon_{u},e_{u}\right)+\varepsilon\left(\frac{A^{n}u_{h}^{n}-A^{n-1} u_{h}^{n-1}}{2},e_{u}\right)+\frac{1}{\varepsilon}\left(u_{h}^{n}-p^{n},e_{u}\right)\\
&+\left(\frac{P^{n}f(u_{h}^{n})-P^{n-1}f(u_{h}^{n-1})}{2\varepsilon},e_{u}\right)+\left(-\frac{w_{h}^{n}-w_{h}^{n-1}}{2},e_{u}\right)\\
&+\left(\left(\varepsilon\mathcal{A}\left(p(t)-p^{n}\right)-\frac{1}{\varepsilon}\left(f(p^{n})\frac{t_{n}-t}{\tau_{n}}-f(p^{n-1})\frac{t_{n}-t}{\tau_{n}}\right)-(q(t)-q^{n})\right),e_{u}\right)\\
&+\left(\frac{1}{\varepsilon}\left(f(u_{h}^{n})\frac{t-t_{n-1}}{\tau_{n}}+f(u_{h}^{n-1})\frac{t_{n}-t}{\tau_{n}}-f(p^{n})+f(p^{n})\frac{t_{n}-t}{\tau_{n}}-f(p^{n-1})\frac{t_{n}-t}{\tau_{n}}\right),e_{u}\right)\\
&+\left(\frac{1}{\varepsilon}\left(f(u_{h})-f(u_{h}^{n})\frac{t-t_{n-1}}{\tau_{n}}-f(u_{h}^{n-1})\frac{t_{n}-t}{\tau_{n}})\right),e_{u}\right)+\left(\frac{1}{\varepsilon}\left(f(u)-f(u_{h})\right),e_{u}\right),
\end{align*}
for all $t\in (t_{n-1},t_{n}]$ and each $n=1,2,\ldots,N$. Then integrate with respect to $t$, we get

\begin{flalign}
&\frac{1}{2}\|\nabla\Delta^{-1}e_{u}^{N}\|^{2}+\int_{0}^{T}\varepsilon\|\nabla e_{u}\|^{2}dt
\nonumber\\=
&\frac{1}{2}\|\nabla\Delta^{-1}e_{u}^{0}\|^{2}+\int_{0}^{T}\left(\frac{A^{n}w_{h}^{n}-A^{n-1} w_{h}^{n-1}}{2},-\Delta^{-1}e_{u}\right) dt\nonumber\\&+\int_{0}^{T}\left(w_{h}^{n}-R^{n}w_{h}^{n},-\Delta^{-1}e_{u}\right)dt+\int_{0}^{T}  \left((t-t_{n-\frac{1}{2}})\partial_{n}^{2}u_{h},-\Delta^{-1}e_{u}\right) dt \nonumber\\
&+\int_{0}^{T}\left(\mathcal{A}(q(t)-q^{n})+q(t)-q^{n},-\Delta^{-1}e_{u}\right)dt+\int_{0}^{T}-\varepsilon a\left(\epsilon_{u},e_{u}\right)dt\nonumber\\
&+\int_{0}^{T}\left(q^{n}-q(t)-(w_{h}^{n}-w_{h}),-\Delta^{-1}e_{u}\right)dt+\int_{0}^{T}\left(w_{h}^{n}-w_{h},-\Delta^{-1}e_{u}\right)dt\nonumber\\
&+\int_{0}^{T}\varepsilon\left(\frac{A^{n}u_{h}^{n}-A^{n-1}u_{h}^{n-1}}{2},e_{u}\right) dt+\int_{0}^{T}\frac{1}{\varepsilon}\left(u_{h}^{n}-p^{n},e_{u}\right)dt\nonumber\\
&+\int_{0}^{T}\left(\frac{P^{n}f(u_{h}^{n})-P^{n-1}f(u_{h}^{n-1})}{2\varepsilon},e_{u}\right) dt +\int_{0}^{T}\left(-\frac{w_{h}^{n}-w_{h}^{n-1}}{2},e_{u}\right) dt\nonumber\\
&+\int_{0}^{T}\Big(\Big(\varepsilon\mathcal{A}\Big(p(t)-p^{n}\Big)-\frac{1}{\varepsilon}\Big(h(p^{n})\frac{t_{n}-t}{\tau_{n}}-h(p^{n-1})\frac{t_{n}-t}{\tau_{n}}\nonumber\\&
-\frac{1}{2}(t-t_{n-1})(t-t_{n})
\frac{\frac{h(p^{n})-h(p^{n-1})}{\tau_{n}}-\frac{h(p^{n-1})-h(p^{n-2})}{\tau_{n-1}}}{\frac{\tau_{n}+\tau_{n-1}}{2}}\Big)-\Big(q(t)-q^{n}\Big)\Big),e_{u}\Big) dt\nonumber\\
&+\int_{0}^{T}\Big(\frac{1}{\varepsilon}\left(p^{n}-p^{n-1}-\left(u_{h}^{n}-u_{h}^{n-1}\right)\right)\frac{t_{n}-t}{\tau_{n}}+\frac{1}{2\varepsilon}(t-t_{n-1})(t-t_{n})\nonumber\\
&\quad\Big(\frac{\frac{p^{n}-p^{n-1}}{\tau_{n}}-\frac{p^{n-1}-p^{n-2}}{\tau_{n-1}}}{\frac{\tau_{n}+\tau_{n-1}}{2}}-\frac{\frac{u_{h}^{n}-u_{h}^{n-1}}{\tau_{n}}-\frac{u_{h}^{n-1}-u_{h}^{n-2}}{\tau_{n-1}}}{\frac{\tau_{n}+\tau_{n-1}}{2}}\Big),e_{u}\Big) dt+\int_{0}^{T}\left(\frac{1}{\varepsilon}\left(u_{h}^{n}-u_{h}^{n-1}\right)\frac{t_{n}-t}{\tau_{n}},e_{u}\right) dt\nonumber\\
&\quad\int_{0}^{T}\Big(\frac{1}{\varepsilon}\left(f(u_{h}^{n})\frac{t-t_{n-1}}{\tau_{n}}+f(u_{h}^{n-1})\frac{t_{n}-t}{\tau_{n}}-f(p^{n})+f(p^{n})\frac{t_{n}-t}{\tau_{n}}-f(p^{n-1})\frac{t_{n}-t}{\tau_{n}}\right)\nonumber\\&
+\frac{1}{2\varepsilon}(t-t_{n-1})(t-t_{n})
\left(\frac{\frac{f(u_{h}^{n})-f(u_{h}^{n-1})}{\tau_{n}}-\frac{f(u_{h}^{n-1})-f(u_{h}^{n-2})}{\tau_{n-1}}}{\frac{\tau_{n}+\tau_{n-1}}{2}}-\frac{\frac{f(p^{n})-f(p^{n-1})}{\tau_{n}}-\frac{f(p^{n-1})-f(p^{n-2})}{\tau_{n-1}}}{\frac{\tau_{n}+\tau_{n-1}}{2}}\right),e_{u}\Big) dt\nonumber\\
&\quad\int_{0}^{T}\Big(\frac{1}{\varepsilon}\Big(f(u_{h})-f(u_{h}^{n})\frac{t-t_{n-1}}{\tau_{n}}-f(u_{h}^{n-1})\frac{t_{n}-t}{\tau_{n}}-\frac{1}{2}(t-t_{n-1})(t-t_{n})\nonumber\\
&\quad\frac{\frac{f(u_{h}^{n})-f(u_{h}^{n-1})}{\tau_{n}}-\frac{f(u_{h}^{n-1})-f(u_{h}^{n-2})}{\tau_{n-1}}}{\frac{\tau_{n}+\tau_{n-1}}{2}}\Big),e_{u}\Big) dt\nonumber\\
&+\int_{0}^{T}\left(\frac{1}{\varepsilon}\left(f(u)-f(u_{h})\right),e_{u}\right) dt\nonumber\\
:=&\frac{1}{2}\|\nabla\Delta^{-1}e_{u}^{0}\|^{2}+\mathcal{B}_{1}+\cdots+\mathcal{B}_{17},\label{2e15}
\end{flalign}
where
\begin{align*}
\mathcal{B}_{1}&:=\int_{0}^{T}\left(w_{h}^{n}-R^{n}w_{h}^{n},-\Delta^{-1}e_{u}\right)dt;\\
\mathcal{B}_{2}&:=\int_{0}^{T}\left(\frac{A^{n}w_{h}^{n}-A^{n-1} w_{h}^{n-1}}{2},-\Delta^{-1}e_{u}\right) dt;\\
\mathcal{B}_{3}&:=\int_{0}^{T}  \left((t-t_{n-\frac{1}{2}})\partial_{n}^{2}u_{h},-\Delta^{-1}e_{u}\right) dt;\\
\mathcal{B}_{4}&:=\int_{0}^{T}\left(\mathcal{A}(q(t)-q^{n})+q(t)-q^{n},-\Delta^{-1}e_{u}\right)dt;\\
\mathcal{B}_{5}&:=\int_{0}^{T}-\varepsilon a\left(\epsilon_{u},e_{u}\right)dt;\\
\mathcal{B}_{6}&:=\int_{0}^{T}\left(q^{n}-q(t)-(w_{h}^{n}-w_{h}),-\Delta^{-1}e_{u}\right)dt;\\
\mathcal{B}_{7}&:=\int_{0}^{T}\left(w_{h}^{n}-w_{h},-\Delta^{-1}e_{u}\right)dt;\\
\mathcal{B}_{8}&:=\int_{0}^{T}\varepsilon\left(\frac{A^{n}u_{h}^{n}-A^{n-1} u_{h}^{n-1}}{2},e_{u}\right) dt;\\
\mathcal{B}_{9}&:=\int_{0}^{T}\frac{1}{\varepsilon}\left(u_{h}^{n}-p^{n},e_{u}\right)dt;\\
\mathcal{B}_{10}&:=\int_{0}^{T}\left(\frac{P^{n}f(u_{h}^{n})-P^{n-1}f(u_{h}^{n-1})}{2\varepsilon},e_{u}\right) dt;\\
\mathcal{B}_{11}&:=-\int_{0}^{T}\left(\frac{ w_{h}^{n}-w_{h}^{n-1}}{2},e_{u}\right) dt;\\
\mathcal{B}_{12}&:=\int_{0}^{T}\Big(\Big(\varepsilon\mathcal{A}\Big(p(t)-p^{n}\Big)-\frac{1}{\varepsilon}\Big(h(p^{n})\frac{t_{n}-t}{\tau_{n}}-h(p^{n-1})\frac{t_{n}-t}{\tau_{n}}\nonumber\\
&\quad-\frac{1}{2}(t-t_{n-1})(t-t_{n})
\frac{\frac{h(p^{n})-h(p^{n-1})}{\tau_{n}}-\frac{h(p^{n-1})-h(p^{n-2})}{\tau_{n-1}}}{\frac{\tau_{n}+\tau_{n-1}}{2}}\Big)-\Big(q(t)-q^{n}\Big)\Big),e_{u}\Big) dt;\\
\mathcal{B}_{13}&:=\int_{0}^{T}\Big(\frac{1}{\varepsilon}\left(p^{n}-p^{n-1}-\left(u_{h}^{n}-u_{h}^{n-1}\right)\right)\frac{t_{n}-t}{\tau_{n}}+\frac{1}{2\varepsilon}(t-t_{n-1})(t-t_{n})\nonumber\\
&\quad\Big(\frac{\frac{p^{n}-p^{n-1}}{\tau_{n}}-\frac{p^{n-1}-p^{n-2}}{\tau_{n-1}}}{\frac{\tau_{n}+\tau_{n-1}}{2}}-\frac{\frac{u_{h}^{n}-u_{h}^{n-1}}{\tau_{n}}-\frac{u_{h}^{n-1}-u_{h}^{n-2}}{\tau_{n-1}}}{\frac{\tau_{n}+\tau_{n-1}}{2}}\Big),e_{u}\Big) dt;\\
\mathcal{B}_{14}&:=\int_{0}^{T}\left(\frac{1}{\varepsilon}\left(u_{h}^{n}-u_{h}^{n-1}\right)\frac{t_{n}-t}{\tau_{n}},e_{u}\right) dt;\\
\mathcal{B}_{15}&:=\int_{0}^{T}\Big(\frac{1}{\varepsilon}\left(f(u_{h}^{n})\frac{t-t_{n-1}}{\tau_{n}}+f(u_{h}^{n-1})\frac{t_{n}-t}{\tau_{n}}-f(p^{n})+f(p^{n})\frac{t_{n}-t}{\tau_{n}}-f(p^{n-1})\frac{t_{n}-t}{\tau_{n}}\right)\nonumber\\&
\quad+\frac{1}{2\varepsilon}(t-t_{n-1})(t-t_{n})
\left(\frac{\frac{f(u_{h}^{n})-f(u_{h}^{n-1})}{\tau_{n}}-\frac{f(u_{h}^{n-1})-f(u_{h}^{n-2})}{\tau_{n-1}}}{\frac{\tau_{n}+\tau_{n-1}}{2}}-\frac{\frac{f(p^{n})-f(p^{n-1})}{\tau_{n}}-\frac{f(p^{n-1})-f(p^{n-2})}{\tau_{n-1}}}{\frac{\tau_{n}+\tau_{n-1}}{2}}\right),e_{u}\Big) dt;\\
\mathcal{B}_{16}&:=\int_{0}^{T}\Big(\frac{1}{\varepsilon}\Big(f(u_{h})-f(u_{h}^{n})\frac{t-t_{n-1}}{\tau_{n}}-f(u_{h}^{n-1})\frac{t_{n}-t}{\tau_{n}}\nonumber\\
&\quad-\frac{1}{2}(t-t_{n-1})(t-t_{n})\frac{\frac{f(u_{h}^{n})-f(u_{h}^{n-1})}{\tau_{n}}-\frac{f(u_{h}^{n-1})-f(u_{h}^{n-2})}{\tau_{n-1}}}{\frac{\tau_{n}+\tau_{n-1}}{2}}\Big),e_{u}\Big) dt;\\
\mathcal{B}_{17}&:=\int_{0}^{T}\left(\frac{1}{\varepsilon}\left(f(u)-f(u_{h})\right),e_{u}\right) dt.
\end{align*}
Next we estimate each of the terms $\{\mathcal{B}_{j}\}_{j=1,\ldots,17}$, separately.

Step $2$: First, the term $\mathcal{B}_{1}$, which contains a spatial discretization error term, is bounded by using Schwarz inequality
\begin{flalign}
\left|\mathcal{B}_{1}\right|=&\left|\int_{0}^{T}\left(w_{h}^{n}-R^{n}w_{h}^{n},-\Delta^{-1}e_{u}\right)dt\right|\nonumber\\
=&\left|\sum_{n=1}^{N}\int_{t_{n-1}}^{t_{n}}\left(w_{h}^{n}-R^{n}w_{h}^{n},-\Delta^{-1}e_{u}\right)dt\right|\nonumber\\
\leq&\sum_{n=1}^{N}\int_{t_{n-1}}^{t_{n}}\left\|w_{h}^{n}-R^{n}w_{h}^{n}\right\|\cdot\left\|\Delta^{-1}e_{u}\right\|dt\nonumber\\
\leq&\sum_{n=1}^{N}\int_{t_{n-1}}^{t_{n}}C_{1}h\left\|w_{h}^{n}-R^{n}w_{h}^{n}\right\|_{1,\Omega}\left\|\Delta^{-1}e_{u}\right\|dt.\label{2e16e}
\end{flalign}
Owing to Remark \ref{ass2e10}, it can be ignored while $h$ is small enough. In the same way, the term $\mathcal{B}_{9}$ can be also ignored. 

Similarly, the time discretization terms $\mathcal{B}_{2}$ and $\mathcal{B}_{3}$ can be estimated as
\begin{flalign}
\left|\mathcal{B}_{2}\right|&=\left|\int_{0}^{T}\left(\frac{A^{n}w_{h}^{n}-A^{n-1}w_{h}^{n-1}}{2},-\Delta^{-1}e_{u}\right) dt\right|\nonumber\\
&\leq\sum_{n=1}^{N}\int_{t_{n-1}}^{t_{n}}\left\|\frac{A^{n}w_{h}^{n}-A^{n-1} w_{h}^{n-1}}{2}\right\|_{-1}\cdot\left\|\nabla\Delta^{-1}e_{u}\right\|dt\nonumber\\
&:=\sum_{n=1}^{N}\int_{t_{n-1}}^{t_{n}}\gamma_{w}^{n}\left\|\nabla\Delta^{-1}e_{u}\right\|dt.\label{2e16ab}
\end{flalign}
\begin{flalign}
\left|\mathcal{B}_{3}\right|&=\left|\int_{0}^{T}  \left((t-t_{n-\frac{1}{2}})\partial_{n}^{2}u_{h},-\Delta^{-1}e_{u}\right) dt\right|\nonumber\\
&\leq\sum_{n=1}^{N}\int_{t_{n-1}}^{t_{n}}\left\|\frac{\tau_{n}^{2}}{8}\cdot\partial_{n}^{2}u_{h}\right\|_{-1}\cdot\left\|\nabla\Delta^{-1}e_{u}\right\|dt\nonumber\\
&:=\sum_{n=1}^{N}\int_{t_{n-1}}^{t_{n}}\beta_{u}^{n}\left\|\nabla\Delta^{-1}e_{u}\right\|dt.\label{2e16b}
\end{flalign}

Step $3$: For the term $\mathcal{B}_{4}$, based on the definition of elliptic reconstruction, we have
\begin{flalign}
|\mathcal{B}_{4}|&=\left|\int_{0}^{T}\left(\mathcal{A}(q(t)-q^{n})+q(t)-q^{n},-\Delta^{-1}e_{u}\right)dt\right|\nonumber\\
&=\Big|\sum_{n=1}^{N}\int_{t_{n-1}}^{t_{n}}\Big(\mathcal{A}\Big(\frac{t-t_{n-1}}{\tau_{n}}q^{n}+\frac{t_{n}-t}{\tau_{n}}q^{n-1}+\frac{1}{2}(t-t_{n-1})(t-t_{n})\partial_{n}^{2}q-q^{n}\Big)\nonumber\\&
\quad+\Big(\frac{t-t_{n-1}}{\tau_{n}}q^{n}+\frac{t_{n}-t}{\tau_{n}}q^{n-1}+\frac{1}{2}(t-t_{n-1})(t-t_{n})\partial_{n}^{2}q-q^{n}\Big),-\Delta^{-1}e_{u}\Big) dt\Big|\nonumber\\
&=\Big|\sum_{n=1}^{N}\int_{t_{n-1}}^{t_{n}}\Big((((\mathcal{A}R^{n-1}w_{h}^{n-1}+R^{n-1}w_{h}^{n-1})-(\mathcal{A}R^{n}w_{h}^{n}+R^{n}w_{h}^{n}))\frac{t_{n}-t}{\tau_{n}}\nonumber\\&
\quad+\frac{1}{2}(t-t_{n-1})(t-t_{n})(\partial_{n}^{2}\mathcal{A}q+\partial_{n}^{2}q)),-\Delta^{-1}e_{u}\Big) dt\Big|\nonumber\\
&\leq\sum_{n=1}^{N}\int_{t_{n-1}}^{t_{n}}\left\|(\mathcal{A}R^{n-1}w_{h}^{n-1}+R^{n-1}w_{h}^{n-1}-(\mathcal{A}R^{n}w_{h}^{n}+R^{n}w_{h}^{n}))\frac{t_{n}-t}{\tau_{n}}+\frac{\tau_{n}^{2}}{8}\partial_{n}^{2}(\mathcal{A}Rw_{h}+Rw_{h})\right\|_{-1}\cdot\nonumber\\
&\quad\left\|\nabla\Delta^{-1}e_{u}\right\|dt\nonumber\\
&\leq\sum_{n=1}^{N}\int_{t_{n-1}}^{t_{n}}\Big(\left\|A^{n-1}w_{h}^{n-1}+w_{h}^{n-1}-(A^{n}w_{h}^{n}+w_{h}^{n})\right\|_{-1}+\left\|\frac{\tau_{n}^{2}}{8}\partial_{n}^{2}(Aw_{h}+w_{h})\right\|_{-1}\Big)\cdot\nonumber\\
&\quad\left\|\nabla\Delta^{-1}e_{u}\right\|dt\nonumber\\
:&=\sum_{n=1}^{N}\int_{t_{n-1}}^{t_{n}}\eta_{w}^{n}\cdot\left\|\nabla\Delta^{-1}e_{u}\right\|dt.\label{2e16c}
\end{flalign}

Step $4$: The term $\mathcal{B}_{5}$ yields the spatial discretization error, which is bounded as follows
\begin{flalign}
|\mathcal{B}_{5}|&=\left|-\int_{0}^{T}\varepsilon a(\epsilon_{u},e_{u})dt\right|\nonumber\\
&=\left|-\sum_{n=1}^{N}\int_{t_{n-1}}^{t_{n}}\varepsilon a(\epsilon_{u},e_{u})dt\right|\nonumber\\
&\leq\left|\sum_{n=1}^{N}\int_{t_{n-1}}^{t_{n}}\varepsilon\|\nabla\epsilon_{u}\|\cdot\|\nabla e_{u}\|dt\right|\nonumber\\
&\leq\left|\sum_{n=1}^{N}\int_{t_{n-1}}^{t_{n}}\left(2\|\nabla\epsilon_{u}\|^{2}+\frac{\varepsilon^{2}}{8}\|\nabla e_{u}\|^{2}\right)dt\right|\nonumber\\
&\leq\sum_{n=1}^{N}\int_{t_{n-1}}^{t_{n}} 2\left\|\nabla\epsilon_{u}\right\|^{2}dt+\sum_{n=1}^{N}\int_{t_{n-1}}^{t_{n}}\frac{\varepsilon^{2}}{8}\left\|\nabla e_{u}\right\|^{2}dt,\label{2e16}
\end{flalign}
and in view of the triangle inequality and the linearity of the operators $G$ and $\nabla$, we get
\begin{flalign}
\int_{t_{n-1}}^{t_{n}}\|\nabla\epsilon_{u}\|^{2}dt
=&\int_{t_{n-1}}^{t_{n}}\|\nabla(p-u_{h})\|^{2}dt\nonumber\\
=&\int_{t_{n-1}}^{t_{n}}\Big\|\nabla\Big(\Big(R^{n}u_{h}^{n}-u_{h}^{n}\Big)\frac{t-t_{n-1}}{\tau_{n}}+\Big(R^{n-1}u_{h}^{n-1}-u_{h}^{n-1}\Big)\frac{t_{n}-t}{\tau_{n}}\nonumber\\&
+\frac{1}{2}(t-t_{n-1})(t-t_{n})(\partial_{n}^{2}p-\partial_{n}^{2}u_{h})\Big)\Big\|^{2}dt\nonumber\\
\leq&\int_{t_{n-1}}^{t_{n}}\bigg(\big\|\nabla\big(R^{n}u_{h}^{n}-u_{h}^{n}\big)\big\|^{2}\Big(\frac{t-t_{n-1}}{\tau_{n}}\Big)^{2}+\big\|\nabla\big(R^{n-1}u_{h}^{n-1}-u_{h}^{n-1}\big)\big\|^{2}\Big(\frac{t_{n}-t}{\tau_{n}}\Big)^{2}\nonumber
\\&+\frac{(t-t_{n-1})^{2}(t-t_{n})^{2}}{4}\|\partial_{n}^{2}\nabla p-\partial_{n}^{2}\nabla u_{h}\|^{2}\nonumber
\\&+2\big\|\nabla\big(R^{n}u_{h}^{n}-u_{h}^{n}\big)\big\|\cdot \big\|\nabla\big(R^{n-1}u_{h}^{n-1}-u_{h}^{n-1}\big)\big\|\frac{(t-t_{n-1})(t_{n}-t)}{\tau_{n}^{2}}\nonumber\\&
+2\big\|\nabla\big(R^{n}u_{h}^{n}-u_{h}^{n}\big)\big\|\cdot \|\partial_{n}^{2}\nabla p-\partial_{n}^{2}\nabla u_{h}\|\frac{(t-t_{n-1})^{2}(t-t_{n})}{2\tau_{n}}
\nonumber\\&
+2\big\|\nabla\big(R^{n-1}u_{h}^{n-1}-u_{h}^{n-1}\big)\big\|\cdot \|\partial_{n}^{2}\nabla p-\partial_{n}^{2}\nabla u_{h}\|\frac{(t-t_{n-1})(t-t_{n})^{2}}{2\tau_{n}}
\bigg)dt\nonumber\\
\leq&C_{0}^{2}\left((\mathcal{E}_{u}^{n})^{2}\frac{(t-t_{n-1})^{3}}{3(t_{n}-t_{n-1})^{2}}\left.\right|_{t_{n-1}}^{t_{n}}-(\mathcal{E}_{u}^{n-1})^{2}\frac{(t_{n}-t)^{3}}{3(t_{n}-t_{n-1})^{2}}\left.\right|_{t_{n-1}}^{t_{n}}\right.\nonumber\\
&\left.+2\mathcal{E}_{u}^{n}\mathcal{E}_{u}^{n-1}\frac{\frac{1}{2}(t_{n}+t_{n-1})t^{2}-\frac{1}{3}t^{3}-t_{n}t_{n-1}\cdot t}{(t_{n}-t_{n-1})^{2}}\left.\right|_{t_{n-1}}^{t_{n}}\right)\nonumber\\&
+\frac{\tau_{n}^{3}}{12}C_{0}\|\partial_{n}^{2}\nabla p-\partial_{n}^{2}\nabla u_{h}\|\Big(\mathcal{E}_{u}^{n}+\mathcal{E}_{u}^{n-1}\Big)+\frac{\tau_{n}^{5}}{120}\|\partial_{n}^{2}\nabla p-\partial_{n}^{2}\nabla u_{h}\|^{2}\nonumber\\
\leq&\frac{C_{0}^{2}}{3}\tau_{n}\left((\mathcal{E}_{u}^{n})^{2}+(\mathcal{E}_{u}^{n-1})^{2}+\mathcal{E}_{u}^{n}\mathcal{E}_{u}^{n-1}\right)\nonumber\\&
+C_{0}^{2}\frac{\tau_{n}^{2}\tau_{n-1}\Big(\mathcal{E}_{u}^{n}+\mathcal{E}_{u}^{n-1}\Big)+\tau_{n}^{3}\Big(\mathcal{E}_{u-1}^{n}+\mathcal{E}_{u}^{n-2}\Big)}{6\tau_{n-1}(\tau_{n}+\tau_{n-1})}\Big(\mathcal{E}_{u}^{n}+\mathcal{E}_{u}^{n-1}\Big)\nonumber\\
&
+C_{0}^{2}\tau_{n}^{3}\frac{\Big(\tau_{n-1}\Big(\mathcal{E}_{u}^{n}+\mathcal{E}_{u}^{n-1}\Big)+\tau_{n}\Big(\mathcal{E}_{u-1}^{n}+\mathcal{E}_{u}^{n-2}\Big)\Big)^{2}}{30\tau_{n-1}^{2}(\tau_{n}+\tau_{n-1})^{2}}\nonumber\\
:=&\widetilde{\mathcal{E}_{u}^{n}}^{2},\label{2e17}
\end{flalign}
then taking \eqref{2e17} into \eqref{2e16}, we get
\bq\label{2e18}
\begin{aligned}
|\mathcal{B}_{6}|\leq\sum_{n=1}^{N}2\widetilde{\mathcal{E}_{u}^{n}}^{2}+\sum_{n=1}^{N}\int_{t_{n-1}}^{t_{n}}\frac{\varepsilon^{2}}{8}\left\|\nabla e_{u}\right\|^{2}dt.
\end{aligned}
\eq
Similarly, the spatial discretization error term $\mathcal{B}_{6}$ is estimated as follows
\begin{flalign}
\left|\mathcal{B}_{6}\right|
&:=\left|\int_{0}^{T}\left(q^{n}-q(t)-(w_{h}^{n}-w_{h}),-\Delta^{-1}e_{u}\right)dt\right|\nonumber\\
&=\Big|\sum_{n=1}^{N}\int_{t_{n-1}}^{t_{n}}\Big(q^{n}-(q^{n}\frac{t-t_{n-1}}{\tau_{n}}+q^{n-1}\frac{t_{n}-t}{\tau_{n}}+\frac{1}{2}(t-t_{n-1})(t-t_{n})\partial_{n}^{2}q)\nonumber\\&
\quad-(w_{h}^{n}-(w_{h}^{n}\frac{t-t_{n-1}}{\tau_{n}}+w_{h}^{n-1}\frac{t_{n}-t}{\tau_{n}}+\frac{1}{2}(t-t_{n-1})(t-t_{n})\partial_{n}^{2}w_{h})),-\Delta^{-1}e_{u}\Big)dt\Big|\nonumber\\
&=\Big|\sum_{n=1}^{N}\int_{t_{n-1}}^{t_{n}}\Big(\Big(q^{n}-w_{h}^{n}\Big)\frac{t_{n}-t}{\tau_{n}}-\Big(q^{n-1}-w_{h}^{n-1}\Big)\frac{t_{n}-t}{\tau_{n}}\nonumber\\ &\quad-\frac{(t-t_{n-1})(t-t_{n})}{2}(\partial_{n}^{2}q-\partial_{n}^{2}w_{h}),-\Delta^{-1}e_{u}\Big)dt\Big|\nonumber\\
&\leq\sum_{n=1}^{N}\int_{t_{n-1}}^{t_{n}}\Big(\|q^{n}-w_{h}^{n}\|\cdot\|\Delta^{-1}e_{u}\|+\|q^{n-1}-w_{h}^{n-1}\|\cdot\|\Delta^{-1}e_{u}\|\nonumber\\ &\quad
+\|\frac{\tau_{n}^{2}}{8}\partial_{n}^{2}(q-w_{h})\|\cdot\|\Delta^{-1}e_{u}\|\Big)dt\nonumber\\
&\leq\sum_{n=1}^{N}\int_{t_{n-1}}^{t_{n}}\Big(Ch\left\|q^{n}-w_{h}^{n}\right\|_{1,\Omega}+Ch\left\|q^{n-1}-w_{h}^{n-1}\right\|_{1,\Omega}\nonumber\\ &\quad+Ch\|\frac{\tau_{n}^{2}}{8}\partial_{n}^{2}(q-w_{h})\|_{1,\Omega}\Big)\left\|\Delta^{-1}e_{u}\right\|dt.\label{2e18a}
\end{flalign}
According to Remark \ref{ass2e10}, it can be ignored while $h$ is small enough. In the same way, the term $\mathcal{B}_{13}$ can be also ignored.

Step $5$: The time discretization term $\mathcal{B}_{7}$ is bounded as follows
\begin{flalign}
|\mathcal{B}_{7}|&=\left|\int_{0}^{T}\left(w_{h}^{n}-w_{h},-\Delta^{-1}e_{u}\right)dt\right|\nonumber\\
&=\left|\sum_{n=1}^{N}\int_{t_{n-1}}^{t_{n}}\left((w_{h}^{n}-w_{h}^{n-1})\frac{t_{n}-t}{\tau_{n}}-\frac{1}{2}(t-t_{n-1})(t-t_{n})\partial_{n}^{2}w_{h},-\Delta^{-1}e_{u}\right)dt\right|\nonumber\\
&\leq\sum_{n=1}^{N}\int_{t_{n-1}}^{t_{n}}\left(\left\|w_{h}^{n}-w_{h}^{n-1}\right\|_{-1}+\left\|\frac{\tau_{n}^{2}}{8}\partial_{n}^{2}w_{h}\right\|_{-1}\right)\left\|\nabla\Delta^{-1}e_{u}\right\|dt\nonumber\\
&:=\sum_{n=1}^{N}\int_{t_{n-1}}^{t_{n}}\delta_{w}^{n}\cdot\left\|\nabla\Delta^{-1}e_{u}\right\|dt.\label{2e18b}
\end{flalign}
Similarly, the terms $\mathcal{B}_{14}$, $\mathcal{B}_{8}$, $\mathcal{B}_{10}$, $\mathcal{B}_{11}$ are also the time discretization terms, and they are estimated as follows
\begin{flalign}
|\mathcal{B}_{14}|&=\left|\int_{0}^{T}\left(\frac{1}{\varepsilon}\left(u_{h}^{n}-u_{h}^{n-1}\right)\frac{t_{n}-t}{\tau_{n}},e_{u}\right)\right|\nonumber\\
&\leq\sum_{n=1}^{N}\int_{t_{n-1}}^{t_{n}}\left\|\frac{u_{h}^{n}-u_{h}^{n-1}}{\varepsilon}\right\|\left\|e_{u}\right\|dt\nonumber\\
&:=\sum_{n=1}^{N}\int_{t_{n-1}}^{t_{n}}\delta_{u}^{n}\cdot\left\|e_{u}\right\|dt,\label{2e18c}
\end{flalign}

\begin{flalign}
|\mathcal{B}_{8}|&=\left|\int_{0}^{T}\varepsilon\left(\frac{A^{n}u_{h}^{n}-A^{n-1}u_{h}^{n-1}}{2},e_{u}\right) dt\right|\nonumber\\
&\leq\sum_{n=1}^{N}\int_{t_{n-1}}^{t_{n}}\varepsilon\left\|\frac{A^{n}u_{h}^{n}-A^{n-1} u_{h}^{n-1}}{2}\right\|\cdot\left\|e_{u}\right\|dt\nonumber\\
:&=\sum_{n=1}^{N}\int_{t_{n-1}}^{t_{n}}\gamma_{u}^{n}\left\|e_{u}\right\|dt,\label{2e19}
\end{flalign}


\begin{flalign}
|\mathcal{B}_{10}|=&\left|\int_{0}^{T}\left(\frac{P^{n}f(u_{h}^{n})-P^{n-1}f(u_{h}^{n-1})}{2\varepsilon},e_{u}\right) dt\right|\nonumber\\
\leq&\sum_{n=1}^{N}\int_{t_{n-1}}^{t_{n}}\left\|\frac{P^{n}f(u_{h}^{n})-P^{n-1}f(u_{h}^{n-1})}{2\varepsilon}\right\|\left\|e_{u}\right\|dt\nonumber\\
:=&\sum_{n=1}^{N}\int_{t_{n-1}}^{t_{n}}\xi_{u}^{n}\cdot\left\|e_{u}\right\|dt,\label{2e20a}\\
|\mathcal{B}_{11}|=&\left|-\int_{0}^{T}\left(\frac{w_{h}^{n}-w_{h}^{n-1}}{2},e_{u}\right) dt\right|\nonumber\\
\leq&\sum_{n=1}^{N}\int_{t_{n-1}}^{t_{n}}\left\|\frac{w_{h}^{n}-w_{h}^{n-1}}{2}\right\|\left\|e_{u}\right\|dt\nonumber\\
:=&\sum_{n=1}^{N}\int_{t_{n-1}}^{t_{n}}\beta_{w}^{n}\cdot\left\|e_{u}\right\|dt.\label{2e20}
\end{flalign}

Step $6$: The term $\mathcal{B}_{12}$, which also contains a time discretization term, is estimated by using the definitions of $w^{n}$, $w$ and $R^{n}$,
\begin{align}
|\mathcal{B}_{12}|&=\Big|\int_{0}^{T}\Big(\Big(\varepsilon\mathcal{A}\Big(p(t)-p^{n}\Big)-\frac{1}{\varepsilon}\Big(h(p^{n})\frac{t_{n}-t}{\tau_{n}}-h(p^{n-1})\frac{t_{n}-t}{\tau_{n}}\nonumber\\
&\quad-\frac{1}{2}(t-t_{n-1})(t-t_{n})
\frac{\frac{h(p^{n})-h(p^{n-1})}{\tau_{n}}-\frac{h(p^{n-1})-h(p^{n-2})}{\tau_{n-1}}}{\frac{\tau_{n}+\tau_{n-1}}{2}}\Big)-\Big(q(t)-q^{n}\Big)\Big),e_{u}\Big) dt\Big|\nonumber\\
&=\Big|\sum_{n=1}^{N}\int_{t_{n-1}}^{t_{n}}\Big(\Big(\varepsilon\mathcal{A}\Big(p^{n}\frac{t-t_{n-1}}{\tau_{n}}+p^{n-1}\frac{t_{n}-t}{\tau_{n}}+\frac{1}{2}(t-t_{n-1})(t-t_{n})\partial_{n}^{2}p-p^{n}\Big)\nonumber\\&
\quad-\frac{1}{\varepsilon}\Big(h(p^{n})-h(p^{n-1})\Big)\frac{t_{n}-t}{\tau_{n}}+\frac{1}{2\varepsilon}(t-t_{n-1})(t-t_{n})\frac{\frac{h(p^{n})-h(p^{n-1})}{\tau_{n}}-\frac{h(p^{n-1})-h(p^{n-2})}{\tau_{n-1}}}{\frac{\tau_{n}+\tau_{n-1}}{2}}\nonumber\\&
\quad-\Big(q^{n}\frac{t-t_{n-1}}{\tau_{n}}+q^{n-1}\frac{t_{n}-t}{\tau_{n}}+\frac{1}{2}(t-t_{n-1})(t-t_{n})\partial_{n}^{2}q-q^{n}\Big)\Big),e_{u}\Big) dt\Big|\nonumber\\
&=\Big|\sum_{n=1}^{N}\int_{t_{n-1}}^{t_{n}}\Big(\Big(\big(\varepsilon\mathcal{A}R^{n-1}u_{h}^{n-1}+\frac{1}{\varepsilon}h(R^{n-1}u_{h}^{n-1})-q^{n-1}\big)-\big(\varepsilon\mathcal{A}R^{n}u_{h}^{n}+\frac{1}{\varepsilon}h(R^{n}u_{h}^{n})-q^{n}\big)\Big)\frac{t_{n}-t}{\tau_{n}}\nonumber\\
&\quad +\frac{\frac{(\varepsilon\mathcal{A}p^{n}+\frac{ h(p^{n})}{\varepsilon}-q^{n})-(\varepsilon\mathcal{A}p^{n-1}+\frac{ h(p^{n-1})}{\varepsilon}-q^{n-1})}{\tau_{n}}-\frac{(\varepsilon\mathcal{A}p^{n-1}+\frac{ h(p^{n-1})}{\varepsilon}-q^{n-1})-(\varepsilon\mathcal{A}p^{n-2}+\frac{ h(p^{n-2})}{\varepsilon}-q^{n-2})}{\tau_{n-1}}}{\frac{\tau_{n}+\tau_{n-1}}{2}}\nonumber\\&\quad\frac{(t-t_{n-1})(t-t_{n})}{2},e_{u}\Big) dt\Big|\nonumber\\
&\leq\sum_{n=1}^{N}\int_{t_{n-1}}^{t_{n}} \Big(2\Big\|\Big(\varepsilon A^{n-1}u_{h}^{n-1}+\frac{1}{\varepsilon}h(u_{h}^{n-1})-w_{h}^{n-1}\Big)-\Big(\varepsilon A^{n}u_{h}^{n}+\frac{1}{\varepsilon}h(u_{h}^{n})-w_{h}^{n}\Big)\Big\|\nonumber\\
&\quad+\Big\|\Big(\varepsilon A^{n-1}u_{h}^{n-1}+\frac{1}{\varepsilon}h(u_{h}^{n-1})-w_{h}^{n-1}\Big)-\Big(\varepsilon A^{n-2}u_{h}^{n-2}+\frac{1}{\varepsilon}h(u_{h}^{n-2})-w_{h}^{n-2}\Big) \Big\|\Big)\cdot \|e_{u}\|dt\nonumber\\
&=\sum_{n=1}^{N}\int_{t_{n-1}}^{t_{n}} \theta_{u}^{n}\cdot\|e_{u}\|dt.\label{2e21}
\end{align}

Step $7$: The term $\mathcal{B}_{15}$ also yields a time discretization error, which is estimated by using Lagrange mean value theorem and embedding theorem.
\begin{align}\label{2e22}
\left|\mathcal{B}_{16}\right|=&\Big|\int_{0}^{T}\Big(\frac{1}{\varepsilon}\left(f(u_{h}^{n})\frac{t-t_{n-1}}{\tau_{n}}+f(u_{h}^{n-1})\frac{t_{n}-t}{\tau_{n}}-f(p^{n})+f(p^{n})\frac{t_{n}-t}{\tau_{n}}-f(p^{n-1})\frac{t_{n}-t}{\tau_{n}}\right)\nonumber\\&
+\frac{1}{2\varepsilon}(t-t_{n-1})(t-t_{n})
\left(\frac{\frac{f(u_{h}^{n})-f(u_{h}^{n-1})}{\tau_{n}}-\frac{f(u_{h}^{n-1})-f(u_{h}^{n-2})}{\tau_{n-1}}}{\frac{\tau_{n}+\tau_{n-1}}{2}}-\frac{\frac{f(p^{n})-f(p^{n-1})}{\tau_{n}}-\frac{f(p^{n-1})-f(p^{n-2})}{\tau_{n-1}}}{\frac{\tau_{n}+\tau_{n-1}}{2}}\right),e_{u}\Big) dt\Big|\nonumber\\
\leq&\left|\sum_{n=1}^{N}\int_{t_{n-1}}^{t_{n}}\left(\frac{t-t_{n-1}}{\tau_{n}\varepsilon}\left(f(u_{h}^{n})-f(p^{n})\right),e_{u}\right) dt\right|+\left|\sum_{n=1}^{N}\int_{t_{n-1}}^{t_{n}}\left(\frac{t_{n}-t}{\tau_{n}\varepsilon}\left(f(u_{h}^{n-1})-f(p^{n-1})\right),e_{u}\right) dt\right|\nonumber\\
&+\left|\sum_{n=1}^{N}\int_{t_{n-1}}^{t_{n}}\left(\frac{1}{2\varepsilon}(t-t_{n-1})(t-t_{n})\frac{\frac{f(u_{h}^{n})-f(p^{n})}{\tau_{n}}}{\frac{\tau_{n}+\tau_{n-1}}{2}},e_{u}\right) dt\right|\nonumber\\&+\left|\sum_{n=1}^{N}\int_{t_{n-1}}^{t_{n}}\left(\frac{1}{\varepsilon}(t-t_{n-1})(t-t_{n})\frac{\frac{f(u_{h}^{n-1})-f(p^{n-1})}{\tau_{n}}}{\frac{\tau_{n}+\tau_{n-1}}{2}},e_{u}\right) dt\right|\nonumber\\
&+\left|\sum_{n=1}^{N}\int_{t_{n-1}}^{t_{n}}\left(\frac{1}{2\varepsilon}(t-t_{n-1})(t-t_{n})\frac{\frac{f(u_{h}^{n-2})-f(p^{n-2})}{\tau_{n}}}{\frac{\tau_{n}+\tau_{n-1}}{2}},e_{u}\right) dt\right|\nonumber\\
\leq&\sum_{n=1}^{N}\int_{t_{n-1}}^{t_{n}}\frac{3}{2\varepsilon}\left\|f'(\xi_{1})\right\|_{0,3,\Omega}\cdot\left\|u_{h}^{n}-p^{n}\right\|_{0,6,\Omega}\cdot\left\|e_{u}\right\|dt\nonumber\\
&+\sum_{n=1}^{N}\int_{t_{n-1}}^{t_{n}}\frac{2}{\varepsilon}\left\|f'(\xi_{2})\right\|_{0,3,\Omega}\cdot\left\|u_{h}^{n-1}-p^{n-1}\right\|_{0,6,\Omega}\cdot\left\|e_{u}\right\|dt\nonumber\\
&+\sum_{n=1}^{N}\int_{t_{n-1}}^{t_{n}}\frac{1}{2\varepsilon}\left\|f'(\xi_{3})\right\|_{0,3,\Omega}\cdot\left\|u_{h}^{n-2}-p^{n-2}\right\|_{0,6,\Omega}\cdot\left\|e_{u}\right\|dt\nonumber\\
\leq&\sum_{n=1}^{N}\int_{t_{n-1}}^{t_{n}}\frac{3}{2\varepsilon}\left\|f'(\xi_{1})\right\|_{1,\Omega}\cdot\left\|u_{h}^{n}-p^{n}\right\|_{1,\Omega}\cdot\left\|e_{u}\right\|dt\nonumber\\
&+\sum_{n=1}^{N}\int_{t_{n-1}}^{t_{n}}\frac{2}{\varepsilon}\left\|f'(\xi_{2})\right\|_{1,\Omega}\cdot\left\|u_{h}^{n-1}-p^{n-1}\right\|_{1,\Omega}\cdot\left\|e_{u}\right\|dt\nonumber\\
&+\sum_{n=1}^{N}\int_{t_{n-1}}^{t_{n}}\frac{1}{2\varepsilon}\left\|f'(\xi_{3})\right\|_{1,3,\Omega}\cdot\left\|u_{h}^{n-2}-p^{n-2}\right\|_{1,6,\Omega}\cdot\left\|e_{u}\right\|dt\nonumber\\
\leq&\sum_{n=1}^{N}\int_{t_{n-1}}^{t_{n}}\frac{1}{\varepsilon}C\mathcal{E}_{u}^{n}\left\|e_{u}\right\|dt+\sum_{n=1}^{N}\int_{t_{n-1}}^{t_{n}}\frac{1}{\varepsilon}C\mathcal{E}_{u}^{n-1}\left\|e_{u}\right\|dt+\sum_{n=1}^{N}\int_{t_{n-1}}^{t_{n}}\frac{1}{\varepsilon}C\mathcal{E}_{u}^{n-2}\left\|e_{u}\right\|dt\nonumber\\
:=&\sum_{n=1}^{N}\int_{t_{n-1}}^{t_{n}}\alpha_{u}^{n}\left\|e_{u}\right\|dt.
\end{align}

Step $8$: In order to estimate the term $\mathcal{B}_{16}$, which also yields a time discretization error, we first simplify the following formula
\begin{flalign}
&f(u_{h})-f(u_{h}^{n})\frac{t-t_{n-1}}{\tau_{n}}-f(u_{h}^{n-1})\frac{t_{n}-t}{\tau_{n}}-\frac{1}{2}(t-t_{n-1})(t-t_{n})
\frac{\frac{f(u_{h}^{n})-f(u_{h}^{n-1})}{\tau_{n}}-\frac{f(u_{h}^{n-1})-f(u_{h}^{n-2})}{\tau_{n-1}}}{\frac{\tau_{n}+\tau_{n-1}}{2}}\nonumber\\
=&\big[3(u_{h}^{n})^{2}u_{h}^{n-1}-2(u_{h}^{n})^{3}-(u_{h}^{n-1})^{3}\big]\cdot \left(\frac{t-t_{n-1}}{\tau_{n}}\right)^{2}\frac{t_{n}-t}{\tau_{n}}\nonumber\\
&+\big[3u_{h}^{n}(u_{h}^{n-1})^{2}-(u_{h}^{n})^{3}-2(u_{h}^{n-1})^{3}\big]\cdot \frac{t-t_{n-1}}{\tau_{n}}\left(\frac{t_{n}-t}{\tau_{n}}\right)^{2}\nonumber\\
&+\big[3(u_{h}^{n})^{2}\partial_{n}^{2}u_{h}-3u_{h}^{n}u_{h}^{n-1}\partial_{n}^{2}u_{h}\big]\cdot\frac{1}{2}(t-t_{n-1})(t-t_{n})\left(\frac{t-t_{n-1}}{\tau_{n}}\right)^{2}\nonumber\\
&+\big[3(u_{h}^{n-1})^{2}\partial_{n}^{2}u_{h}-3u_{h}^{n}u_{h}^{n-1}\partial_{n}^{2}u_{h}\big]\cdot\frac{1}{2}(t-t_{n-1})(t-t_{n})\left(\frac{t_{n}-t}{\tau_{n}}\right)^{2}\nonumber\\
&+\big[3u_{h}^{n}u_{h}^{n-1}\partial_{n}^{2}u_{h}-\frac{\frac{(u_{h}^{n})^{3}-(u_{h}^{n-1})^{3}}{\tau_{n}}-\frac{(u_{h}^{n-1})^{3}-(u_{h}^{n-2})^{3}}{\tau_{n-1}}}{\frac{\tau_{n}+\tau_{n-1}}{2}}\big]\cdot\frac{1}{2}(t-t_{n-1})(t-t_{n})\nonumber\\
&+\Big[\big(3u_{h}^{n}(\partial_{n}^{2}u_{h})^{2}-3u_{h}^{n-1}(\partial_{n}^{2}u_{h})^{2}\big)\cdot \left(\frac{t-t_{n-1}}{\tau_{n}}\right)+3u_{h}^{n-1}(\partial_{n}^{2}u_{h})^{2}\Big]\cdot\Big(\frac{1}{2}(t-t_{n-1})(t-t_{n})\Big)^{2}\nonumber\\
&+\Big[\Big(\frac{1}{2}(t-t_{n-1})(t-t_{n})\Big)^{3}(\partial_{n}^{2}u_{h})^{3}\Big],\label{2e23}
\end{flalign}
thus we have
\begin{flalign}
\left|\mathcal{B}_{16}\right|
=&\Big|\int_{0}^{T}\Big(\frac{1}{\varepsilon}\Big(f(u_{h})-f(u_{h}^{n})\frac{t-t_{n-1}}{\tau_{n}}-f(u_{h}^{n-1})\frac{t_{n}-t}{\tau_{n}}\nonumber\\&-\frac{1}{2}(t-t_{n-1})(t-t_{n})\frac{\frac{f(u_{h}^{n})-f(u_{h}^{n-1})}{\tau_{n}}-\frac{f(u_{h}^{n-1})-f(u_{h}^{n-2})}{\tau_{n-1}}}{\frac{\tau_{n}+\tau_{n-1}}{2}}\Big),e_{u}\Big) dt\Big|\nonumber\\
=&\Big|\sum_{n=1}^{N}\int_{t_{n-1}}^{t_{n}}\frac{1}{\varepsilon} \Big( \big[3(u_{h}^{n})^{2}u_{h}^{n-1}-2(u_{h}^{n})^{3}-(u_{h}^{n-1})^{3}\big]\cdot \Big(\frac{t-t_{n-1}}{\tau_{n}}\Big)^{2}\frac{t_{n}-t}{\tau_{n}}\nonumber\\
&+\big[3u_{h}^{n}(u_{h}^{n-1})^{2}-(u_{h}^{n})^{3}-2(u_{h}^{n-1})^{3}\big]\cdot \frac{t-t_{n-1}}{\tau_{n}}\Big(\frac{t_{n}-t}{\tau_{n}}\Big)^{2}\nonumber\\
&+\big[3(u_{h}^{n})^{2}\partial_{n}^{2}u_{h}-3u_{h}^{n}u_{h}^{n-1}\partial_{n}^{2}u_{h}\big]\cdot\frac{1}{2}(t-t_{n-1})(t-t_{n})\left(\frac{t-t_{n-1}}{\tau_{n}}\right)^{2}\nonumber\\
&+\big[3(u_{h}^{n-1})^{2}\partial_{n}^{2}u_{h}-3u_{h}^{n}u_{h}^{n-1}\partial_{n}^{2}u_{h}\big]\cdot\frac{1}{2}(t-t_{n-1})(t-t_{n})\left(\frac{t_{n}-t}{\tau_{n}}\right)^{2}\nonumber\\
&+\big[3u_{h}^{n}u_{h}^{n-1}\partial_{n}^{2}u_{h}--\frac{\frac{(u_{h}^{n})^{3}-(u_{h}^{n-1})^{3}}{\tau_{n}}-\frac{(u_{h}^{n-1})^{3}-(u_{h}^{n-2})^{3}}{\tau_{n-1}}}{\frac{\tau_{n}+\tau_{n-1}}{2}}\big]\cdot\frac{1}{2}(t-t_{n-1})(t-t_{n})\nonumber\\
&+\Big[\big(3u_{h}^{n}(\partial_{n}^{2}u_{h})^{2}-3u_{h}^{n-1}(\partial_{n}^{2}u_{h})^{2}\big)\cdot \left(\frac{t-t_{n-1}}{\tau_{n}}\right)+3u_{h}^{n-1}(\partial_{n}^{2}u_{h})^{2}\Big]\cdot\Big(\frac{1}{2}(t-t_{n-1})(t-t_{n})\Big)^{2}\nonumber\\
&+\Big[\Big(\frac{1}{2}(t-t_{n-1})(t-t_{n})\Big)^{3}(\partial_{n}^{2}u_{h})^{3}\Big]
,e_{u}\Big) dt \Big|\nonumber\\
\leq&\sum_{n=1}^{N}\int_{t_{n-1}}^{t_{n}}\left\|\frac{3(u_{h}^{n})^{2}u_{h}^{n-1}-2(u_{h}^{n})^{3}-(u_{h}^{n-1})^{3}}{\varepsilon}\right\|\cdot\left\|e_{u}\right\|dt\nonumber\\
&+ \sum_{n=1}^{N}\int_{t_{n-1}}^{t_{n}}\left\|\frac{3u_{h}^{n}(u_{h}^{n-1})^{2}-2(u_{h}^{n-1})^{3}-(u_{h}^{n})^{3}}{\varepsilon} \right\|\cdot\left\|e_{u}\right\|dt\nonumber\\
&+\sum_{n=1}^{N}\int_{t_{n-1}}^{t_{n}}\left\|\frac{\big(3(u_{h}^{n})^{2}\partial_{n}^{2}u_{h}-3u_{h}^{n}u_{h}^{n-1}\partial_{n}^{2}u_{h}\big)\tau_{n}^{2}}{8\varepsilon} \right\|\cdot\left\|e_{u}\right\|dt\nonumber\\
&+\sum_{n=1}^{N}\int_{t_{n-1}}^{t_{n}}\left\|\frac{\big(3(u_{h}^{n-1})^{2}\partial_{n}^{2}u_{h}-3u_{h}^{n}u_{h}^{n-1}\partial_{n}^{2}u_{h}\big)\tau_{n}^{2}}{8\varepsilon} \right\|\cdot\left\|e_{u}\right\|dt\nonumber\\
&+\sum_{n=1}^{N}\int_{t_{n-1}}^{t_{n}}\left\|\frac{\big(3u_{h}^{n}u_{h}^{n-1}\partial_{n}^{2}u_{h}-\frac{\frac{(u_{h}^{n})^{3}-(u_{h}^{n-1})^{3}}{\tau_{n}}-\frac{(u_{h}^{n-1})^{3}-(u_{h}^{n-2})^{3}}{\tau_{n-1}}}{\frac{\tau_{n}+\tau_{n-1}}{2}}\big)\tau_{n}^{2}}{8\varepsilon} \right\|\cdot\left\|e_{u}\right\|dt\nonumber\\
&+\sum_{n=1}^{N}\int_{t_{n-1}}^{t_{n}}\left(\left\|\frac{\big(3u_{h}^{n}(\partial_{n}^{2}u_{h})^{2}-3u_{h}^{n-1}(\partial_{n}^{2}u_{h})^{2}\big)\tau_{n}^{4}}{64\varepsilon} \right\|+\left\|\frac{\big(3u_{h}^{n-1}(\partial_{n}^{2}u_{h})^{2}\big)\tau_{n}^{4}}{64\varepsilon}\right\|\right)\cdot\left\|e_{u}\right\|dt\nonumber\\
&+\sum_{n=1}^{N}\int_{t_{n-1}}^{t_{n}}\left(\left\|\frac{(\partial_{n}^{2}u_{h})^{3}\tau_{n}^{6}}{512\varepsilon} \right\|\right)\cdot\left\|e_{u}\right\|dt\nonumber\\
:=&\sum_{n=1}^{N}\int_{t_{n-1}}^{t_{n}}\zeta_{u}^{n}\left\|e_{u}\right\|dt.\label{2e23a}
\end{flalign}

Step $9$: Grouping together \eqref{2e16c}, \eqref{2e18a} and \eqref{2e18b}, we have
\begin{flalign}
\left|\mathcal{B}_{1}\right|+\cdots&+\left|\mathcal{B}_{4}\right|+\left|\mathcal{B}_{6}\right|+\left|\mathcal{B}_{7}\right|\nonumber\\
\leq&\sum_{n=1}^{N}\int_{t_{n-1}}^{t_{n}}\gamma_{w}^{n}\left\|\nabla\Delta^{-1}e_{u}\right\|dt+\sum_{n=1}^{N}\int_{t_{n-1}}^{t_{n}}\beta_{u}^{n}\left\|\nabla\Delta^{-1}e_{u}\right\|dt\nonumber\\
&+\sum_{n=1}^{N}\int_{t_{n-1}}^{t_{n}}\delta_{w}^{n}\cdot\left\|\nabla\Delta^{-1}e_{u}\right\|dt+\sum_{n=1}^{N}\int_{t_{n-1}}^{t_{n}}\eta_{w}^{n}\cdot\left\|\nabla\Delta^{-1}e_{u}\right\|dt\nonumber\\
:=&\sum_{n=1}^{N}\int_{t_{n-1}}^{t_{n}}\eta_{0}\left\|\nabla\Delta^{-1}e_{u}\right\|dt\nonumber\\
\leq&\sum_{n=1}^{N}\int_{t_{n-1}}^{t_{n}}\frac{1}{2}\eta_{0}^{2}dt+\sum_{n=1}^{N}\int_{t_{n-1}}^{t_{n}}\frac{1}{2}\left\|\nabla\Delta^{-1}e_{u}\right\|^{2}dt.\label{2e24a}
\end{flalign}
Summing up \eqref{2e18c}-\eqref{2e23a}, it holds that
\begin{align}
\left|\mathcal{B}_{8}\right|+\cdots+\left|\mathcal{B}_{16}\right|\leq&\sum_{n=1}^{N}\int_{t_{n-1}}^{t_{n}}\gamma_{u}^{n}\left\|e_{u}\right\|dt+\sum_{n=1}^{N}\int_{t_{n-1}}^{t_{n}}\xi_{u}^{n}\left\|e_{u}\right\|dt+\sum_{n=1}^{N}\int_{t_{n-1}}^{t_{n}}\beta_{w}^{n}\left\|e_{u}\right\|dt\nonumber\\
&+\sum_{n=1}^{N}\int_{t_{n-1}}^{t_{n}} \theta_{u}^{n}\cdot\|e_{u}\|dt+\sum_{n=1}^{N}\int_{t_{n-1}}^{t_{n}}\delta_{u}^{n}\cdot\left\|e_{u}\right\|dt\nonumber\\&+\sum_{n=1}^{N}\int_{t_{n-1}}^{t_{n}}\alpha_{u}^{n}\left\|e_{u}\right\|dt+\sum_{n=1}^{N}\int_{t_{n-1}}^{t_{n}}\zeta_{u}^{n}\left\|e_{u}\right\|dt\nonumber\\
:=&\sum_{n=1}^{N}\int_{t_{n-1}}^{t_{n}}\eta_{1}\left\|e_{u}\right\|dt\nonumber\\
\leq&\sum_{n=1}^{N}\int_{t_{n-1}}^{t_{n}}\frac{1}{2}\eta_{1}^{2}dt+\sum_{n=1}^{N}\int_{t_{n-1}}^{t_{n}}\frac{1}{2}\left\|e_{u}\right\|^{2}dt.\label{2e24}
\end{align}

Step $10$: As for estimation of the term $\mathcal{B}_{17}$, according to the Remark \ref{rem1e1}, the spectrum  estimate \cite{AF1993,C1994},
 and the fact that
\[\left(f(a)-f(b)-f'(b)(a-b)\right)\left(a-b\right)\geq -\tilde{f}(b)\mid a-b\mid^{3}\]
with $\tilde{f}(b)=3|b|$, we obtain
\begin{align}\label{2e25}
\left|\mathcal{B}_{17}\right|=&\left|\int_{0}^{T}\left(\frac{1}{\varepsilon}\left(f(u)-f(u_{h})\right),e_{u}\right) dt\right|\nonumber\\
\leq&\left|\int_{0}^{T}\left(-\frac{1}{\varepsilon}\left(f'(u_{h})e_{u},e_{u}\right)+\frac{1}{\varepsilon}\left(\tilde{f}(u_{h}),|u-u_{h}|^{3}\right)\right)dt\right|\nonumber\\
\leq&\left|\int_{0}^{T}\left(-\frac{1-\varepsilon}{\varepsilon}\left( f'(u_{h})e_{u},e_{u}\right)-\left(f'(u_{h})e_{u},e_{u}\right)+\frac{1}{\varepsilon}\|\tilde{f}(u_{h})\|_{L^{\infty}(\Omega)}\|e_{u}\|_{L^{3}}^{3}\right)dt\right|\nonumber\\
\leq&\left|\int_{0}^{T}\left(\left(1-\varepsilon\right)\overline{\Lambda}_{CH}(t)\|\nabla\Delta^{-1} e_{u}\|^{2} +\varepsilon\left(1-\varepsilon\right)\|\nabla e_{u}\|^{2}+2\|e_{u}\|^{2}+\frac{1}{\varepsilon}\mu_{g}\|e_{u}\|_{L^{3}}^{3}\right)dt\right|.
\end{align}

Step $11$: Taking \eqref{2e18}, \eqref{2e24a}-\eqref{2e25} into \eqref{2e15}, we have
\begin{flalign}
\frac{1}{2}\|\nabla\Delta^{-1}e_{u}^{N}\|^{2}+\int_{0}^{T}\varepsilon\|\nabla e_{u}\|^{2}dt
\leq&\frac{1}{2}\|\nabla\Delta^{-1}e_{u}^{0}\|^{2}+\sum_{n=1}^{N}2\widetilde{\mathcal{E}_{u}^{n}}^{2}+\sum_{n=1}^{N}\int_{t_{n-1}}^{t_{n}}\frac{\varepsilon^{2}}{8}\left\|\nabla e_{u}\right\|^{2}dt\nonumber\\
&+\sum_{n=1}^{N}\int_{t_{n-1}}^{t_{n}}\frac{1}{2}\eta_{0}^{2}dt+\sum_{n=1}^{N}\int_{t_{n-1}}^{t_{n}}\frac{1}{2}\left\|\nabla\Delta^{-1}e_{u}\right\|^{2}dt\nonumber\\
&+\sum_{n=1}^{N}\int_{t_{n-1}}^{t_{n}}\frac{1}{2}\eta_{1}^{2}dt+\sum_{n=1}^{N}\int_{t_{n-1}}^{t_{n}}\frac{1}{2}\left\|e_{u}\right\|^{2}dt\nonumber\\
&+\int_{0}^{T}\left(1-\varepsilon\right)\overline{\Lambda}_{CH}(t)\|\nabla\Delta^{-1}e_{u}\|^{2}dt\nonumber\\&+\sum_{n=1}^{N}\int_{t_{n-1}}^{t_{n}}2\left\|e_{u}\right\|^{2}dt
+\int_{0}^{T}\varepsilon\left(1-\varepsilon\right)\|\nabla e_{u}\|^{2}dt\nonumber\\&+\int_{0}^{T}\frac{1}{\varepsilon}\mu_{g}\|e_{u}\|_{L^{3}}^{3}dt.\label{2e27a}
\end{flalign}
Note that
\begin{flalign}
\left\|e_{u}\right\|^{2}&=\left(\nabla(-\Delta^{-1}e_{u}),\nabla e_{u}\right)\nonumber\\
&\leq\|\nabla\Delta^{-1}e_{u}\|\|\nabla e_{u}\|\nonumber\\
&\leq\frac{1}{\varepsilon^{2}}\|\nabla\Delta^{-1}e_{u}\|^{2}+\frac{\varepsilon^{2}}{4}\|\nabla e_{u}\|^{2},\label{2e27}
\end{flalign}
plugging \eqref{2e27} into \eqref{2e27a}, and further simplification, then we obtain
\begin{flalign}\label{2e26}
\|\nabla\Delta^{-1}e_{u}^{N}\|^{2}+&\int_{0}^{T}\frac{\varepsilon^{2}}{2}\|\nabla e_{u}\|^{2}dt
\nonumber\\
\leq&\|\nabla\Delta^{-1}e_{u}^{0}\|^{2}+\sum_{n=1}^{N}4\widetilde{\mathcal{E}_{u}^{n}}^{2}+\sum_{n=1}^{N}\left(\eta_{0}^{2}+\eta_{1}^{2}\right)\tau_{n}\nonumber\\
&+\int_{0}^{T}\left(1+\frac{5}{2\varepsilon^{2}}+2\left(1-\varepsilon\right)\overline{\Lambda}_{CH}(t)\right)\|\nabla\Delta^{-1}e_{u}\|^{2}dt\nonumber\\
&+\int_{0}^{T}\frac{2}{\varepsilon}\mu_{g}\|e_{u}\|_{L^{3}}^{3}dt.
\end{flalign}

According to Lemma \ref{lem3e1} and assume that $\|e_{u}\|_{L^{\infty}}\leq C$, then it holds that
\[\begin{aligned}
\int_{0}^{T}\|e_{u}\|_{L^{3}}^{3}dt\leq&\int_{0}^{T}C_{I}\|e_{u}\|_{L^{\infty}(\Omega)}^{1-\sigma}\|\nabla\Delta^{-1}e_{u}\|^{\sigma}\|\nabla e_{u}\|^{2}dt\\
\leq&\int_{0}^{T}C_{I}\|e_{u}\|_{L^{\infty}(\Omega)}^{1-\sigma}\|\nabla\Delta^{-1}e_{u}\|^{\sigma}\|\nabla e_{u}\|^{2}dt\\
\leq& C_{S}\left(\sup_{t\in(0,T)}\|\nabla\Delta^{-1}e_{u}\|^{\sigma}\right)\int_{0}^{T}\|\nabla e_{u}\|^{2}dt.
\end{aligned}
\]
Setting
\begin{align*}
y_{1}(t)&:=\|\nabla\Delta^{-1}e_{u}\|^{2}, \qquad\qquad y_{2}(t):=\frac{\varepsilon^{2}}{2}\|\nabla e_{u}\|^{2},\qquad\qquad\qquad
y_{3}(t):=2\varepsilon^{-1}\mu_{g}\|e_{u}\|_{L^{3}}^{3},\\
B&:=2\varepsilon^{-1}\mu_{g}C_{S}, \qquad\qquad\qquad E:=\exp\left(\int_{0}^{T}a(t)dt\right),\qquad \qquad \beta:=\sigma,
\end{align*}
then by Lemma \ref{lem3e2}, we have
\[
\sup_{t\in[0,T]}\|\nabla\Delta^{-1}e_{u}\|^{2}+\int_{0}^{T}\frac{\varepsilon^{2}}{2}\|\nabla e_{u}\|^{2}dt\leq8\eta^{2}\exp\left(\int_{0}^{T}a(t)dt\right).
\]
\end{proof}

\section*{Acknowledgments}
Chen's research was supported by NSFC Project (12201010), Natural Science Research Project of Higher Education in Anhui Province (2022AH040027). Huang's research was partially supported by NSFC Project (11971410) and China's National Key R\&D Programs  (2020YFA0713500). Yi's research was partially supported by NSFC Project (12071400, 12261131501) and Hunan Provincial NSF Project (2021JJ40189). 




\end{document}